\documentclass[11pt,a4paper]{amsart}

\usepackage[top=2.5cm,bottom=2.5cm,left=2.5cm,right=2.5cm]{geometry}
\usepackage[T1]{fontenc}
\usepackage[utf8]{inputenc}
\usepackage[english]{babel}
\usepackage{float}
\usepackage{xspace}
\usepackage{lmodern}
\DeclareFontFamily{U}{wncyr}{}
\DeclareFontShape{U}{wncyr}{m}{n}{<->wncyr10}{}
\DeclareFontShape{U}{wncyr}{m}{it}{<->wncyi10}{}
\DeclareFontShape{U}{wncyr}{m}{sc}{<->wncysc10}{}
\DeclareFontShape{U}{wncyr}{b}{n}{<->wncyb10}{}
\DeclareTextCommand{\guillemotleft}{T1}{%
  {\fontencoding{U}\fontfamily{wncyr}\selectfont\symbol{"3C}}%
}
\DeclareTextCommand{\guillemotright}{T1}{%
  {\fontencoding{U}\fontfamily{wncyr}\selectfont\symbol{"3E}}%
}

\usepackage{amsmath,amsfonts,amsthm,amssymb}
\usepackage{tikz-cd}
\usepackage{tikz,pgfplots}
\definecolor{PineGreen}{rgb}{0.0, 0.47, 0.44}
\definecolor{BrickRed}{rgb}{0.8, 0.25, 0.33}
\definecolor{bole}{rgb}{0.47, 0.27, 0.23}
\definecolor{amber}{rgb}{1.0, 0.75, 0.0}
\usepackage{color}

\usepackage{graphicx}
\usepackage{ifthen}
\usetikzlibrary{arrows,calc,intersections,patterns}
\usepackage{rotating} 

\usepackage{caption}
\usepackage{subcaption}

\usepackage{enumitem}

\usepackage{cite}

\allowdisplaybreaks[3] 
\usepackage{mathtools}
\usepackage{hyperref}
\usepackage{xcolor}

\newtheorem{theorem}{Theorem}[section]
\newtheorem{proposition}[theorem]{Proposition}
\newtheorem{corollary}[theorem]{Corollary}
\newtheorem{lemma}[theorem]{Lemma}

\newtheorem{theorembis}{Theorem}

\newtheorem{theorembisH}{Theorem}

\theoremstyle{definition}
\newtheorem{definition}[theorem]{Definition}

\theoremstyle{remark}
\newtheorem{remark}[theorem]{Remark}

\numberwithin{equation}{section}
\DeclareRobustCommand{\SkipTocEntry}[5]{}

\renewcommand{\d}{\,\mathrm{d}}

\renewcommand{\leq}{\leqslant}

\renewcommand{\geq}{\geqslant}

\usepackage[]{hyperref}
\hypersetup{
    colorlinks=true,       
    linkcolor=blue,          
    citecolor=magenta,        
    filecolor=magenta,      
    urlcolor=cyan           
}
\begin{document}

\title[Grushin-Schrödinger equation]{Probabilistic Local Well-Posedness for the Schrödinger Equation Posed for the Grushin Laplacian}

\author{Louise Gassot}
\address{Département de Mathématiques et Applications, Ecole Normale Supérieure -- PSL Research University, 45 rue d'Ulm 75005 Paris, France \&
 Université Paris-Saclay, CNRS, Laboratoire de mathématiques d’Orsay, 91405 Orsay, France}
\email{louise.gassot@universite-paris-saclay.fr}
\author{Mickaël Latocca}
\address{Département de Mathématiques et Applications, Ecole Normale Supérieure -- PSL Research University, 45 rue d'Ulm 75005 Paris, France}
\email{mickael.latocca@ens.fr}

\date{\today}

\begin{abstract}
\begin{sloppypar}
We study the local well-posedness of the nonlinear Schr{\"o}dinger equation associated to the Grushin operator with random initial data. To the best of our knowledge, no well-posedness result is known in the Sobolev spaces $H^k$ when $k\leq\frac{3}{2}$. In this article, we prove that there exists a large family of initial data such that, with respect to a suitable randomization in $H^k$, $k \in (1,\frac{3}{2}]$, almost-sure local well-posedness holds. The proof relies on bilinear and trilinear estimates.  
\end{sloppypar}
\end{abstract}

\maketitle

\tableofcontents
\section{Introduction}

\subsection{The Schrödinger equation on the Heisenberg group and the Grushin equation}

We consider the Grushin-Schrödinger equation
\begin{equation}
    \label{eq:NLSG}
    \tag{NLS-G}
    i\partial_tu - \Delta_Gu = |u|^2u\,,
\end{equation}
where $(t,x,y)\in\mathbb{R}\times\mathbb{R}^2$ and $\Delta_G=\partial_{xx}+x^2\partial_{yy}$. The natural associated Sobolev spaces in this case are the Grushin Sobolev spaces $H^k_G$ on $\mathbb{R}^2$, defined by replacing powers of  the usual operator $\sqrt{-\Delta}$ by powers of $\sqrt{-\Delta_G}$.

This equation is a simplification of the semilinear Schrödinger equation on the Heisenberg group in the radial case
\begin{equation}
    \label{eq:NLSH}
     \tag{NLS-$\mathbb{H}^1$}
    i\partial_tu - \Delta_{\mathbb{H}^1}u = |u|^2u\,,
\end{equation}
where $(t,x,y,s)\in \mathbb{R}\times\mathbb{H}^1$. In the radial case, the solution $u$ only depends on $t,|x+iy|$ and $s$ and the sub-Laplacian is written $\Delta_{\mathbb{H}^1}=\frac{1}{4}(\partial_{xx}+\partial_{yy})+(x^2+y^2)\partial_{ss}$.
Our simplification of this equation consists in removing one of the two variables $x,y$ since they play the same role, leading to~\eqref{eq:NLSG}.

When $k>k_C$ where $k_C=2$ (resp. $k_C=\frac{3}{2}$ for~\eqref{eq:NLSG}), one can use the algebra property of the spaces $H^{k}(\mathbb{H}^1)$ (resp. $H^k_G)$) and solve the Cauchy problem associated to \eqref{eq:NLSH} (resp. \eqref{eq:NLSG}) locally in time, see Appendix~\ref{appendix:B} for details. However, the conservation of energy only controls the $H^k$ norm when $k=1$, and since no conservation law is known  for $k>1$, we have no information about global existence of maximal solutions in the range of Sobolev exponents $k>k_C$.

For Sobolev exponents below the critical exponent $k_C$, 
existence and uniqueness of general weak solutions is an open problem. To go further, the Schrödinger equation on the Heisenberg group displays a total lack of dispersion~\cite{BahouriGerardXu2000}, implying that the flow map for~\eqref{eq:NLSH} (resp.~\eqref{eq:NLSG}) cannot be smooth in the Sobolev spaces $H^k$ when $k<k_C$. We refer to the introduction of~\cite{GerardGrellier2010} and Remark~2.12 in~\cite{BurqGerardTzvetkov2004} for details. 

Finally, note that equation~\eqref{eq:NLSG} enjoys the $H^{1/2}_G$-critical scaling invariance
\[
u\mapsto u_\lambda(t,x,y)=\lambda u(\lambda^2t,\lambda x,\lambda^2 y)\,.
\]
As a consequence, the equation can be shown to be ill-posed in $H^s_G$ when $s<1/2$, see for instance~\cite{ChristCollianderTao2003norm}. 
Recently, Sun and Tzvetkov~\cite{SunTzvetkov2020pathological} showed that for nonlinear wave equations, when the initial data are regularized by convolution, ill-posedness under the form of norm inflation of the solution in $H^s$ for arbitrarily small times holds on a pathological set which contains a dense $G_{\delta}$ set. This result is complementary to the existing probabilistic well-posedness theory below this exponent, see their Theorem 1 for more details. According to the extension of this work to the nonlinear Schrödinger equation in~\cite{CampsGassot2022}, one can see that the existence of a dense pathological set would also hold for equations~\eqref{eq:NLSG} and~\eqref{eq:NLSH} below the scaling-critical exponent $1/2$. However, it is an open problem whether local well-posedness holds between the exponents $1/2$ and  $k_C$.

The situation for equation~\eqref{eq:NLSG} can be summarized by the following diagram.
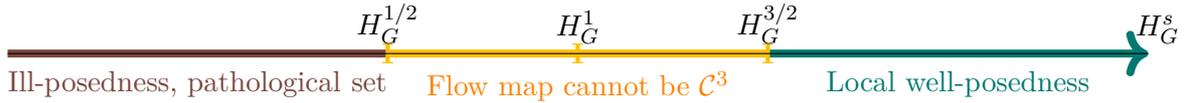
\begin{center}
\begin{tikzpicture}
\draw[line width=3, bole] (-6,0)--(-1,0);
\draw[line width=3, amber] (-1,0)--(4,0);
\draw[bole] (-3.5,-0.1) node[below] {Ill-posedness, pathological set};
\draw[orange] (1.5,-0.1) node[below] {Flow map cannot be $\mathcal{C}^3$};
\draw[PineGreen] (6.5,-0.1) node[below] {Local well-posedness};
\draw[->, line width=3, PineGreen] (4,0)--(9,0);
\draw (8.7,0)  node[above right]{$H^s_G$};
\draw (-1,0) node{{\color{amber} \bf{I}}} node[above]{ $H^{1/2}_G$};
\draw (1.5,0) node{{\color{amber}\bf I}} node[above]{$H^{1}_G$};
\draw (4,0) node{{\color{amber} \bf I}} node[above]{$H^{3/2}_G$};
\draw (-6,0)--(9,0);
\end{tikzpicture}
\end{center}

\subsection{Main results}

Our main result is the following.
\begin{theorembis}[Local Cauchy Theory for~\eqref{eq:NLSG}]\label{thm:main_short}
Let $k\in(1,\frac 32]$ and $\ell \in (\frac{3}{2},k+\frac{1}{2})$. There exists a dense set $A_{k,\ell}$ in $H^k_G$ such that for every $u_0\in A_{k,\ell}$, there exists $T>0$ and a unique local solution with initial data $u_0$ to~\eqref{eq:NLSG} in the space
\[
    e^{it\Delta_G}u_0+ \mathcal{C}^0([0,T),H^{\ell}_G) \subset \mathcal{C}^0([0,T),H^{k}_G)\,.
\]
\end{theorembis}

In this subsection we will only introduce the needed notations to precise to our main result using random initial data. We refer to Section~\ref{sec:preliminaries} for more precise definitions.


Fix $u_0$ in some Sobolev space $H^k_G$ for $k> 0$. Then $u$ decomposes as a sum
\begin{equation}\label{eq:decompo_Grushin}
u_0=\sum_{(I,m)\in 2^{\mathbb{Z}}\times\mathbb{N}}u_{I,m}\,,
\end{equation}
where the $u_{I,m}$ will be defined by~\eqref{def:BuildingBlock},

Let $(\Omega,\mathcal{A}, \mathbb{P})$ be a probability space. 
We consider a sequence $(X_{I,m})_{(I,m)\in 2^{\mathbb{Z}}\times\mathbb{N}}$ of independent and identically distributed Gaussian random variables and define the measure $\mu_{u_0}$ as the image measure of $\mathbb{P}$ under the \textit{randomization map}
\begin{equation}
    \label{eq:randomizationMap}
    \omega \in\Omega\mapsto u_0^{\omega} := \sum_{(I,m)\in 2^{\mathbb{Z}}\times\mathbb{N}} X_{I,m}(\omega)u_{I,m}\,.
\end{equation}

For $k\geq 0$ and $\rho\geqslant 0$, we introduce the subspace $\mathcal{X}^k_\rho$ of $H^k_G$ in which we will prove almost-sure local well-posedness. Denoting $\langle x\rangle=\sqrt{1+x^2}$, the space $\mathcal{X}^k_\rho$ corresponds to the norm
\begin{equation}
    \label{def:Xk_rho}
    \|u_0\|_{\mathcal{X}^k_{\rho}} ^2:= \sum_{(I,m)\in 2^{\mathbb{Z}}\times\mathbb{N}} (1+(2m+1)I)^k \left\langle I \right\rangle ^{\rho} \|u_{I,m}\|^2_{L^2_G}\,.
\end{equation}
The powers $(1+(2m+1)I)^k$ refer to the Sobolev regularity: for instance, when $\rho =0$, then $\|u_0\|_{\mathcal{X}^{k}_0}\sim\|u_0\|_{H^k_G}$. However, the powers $\left\langle I \right\rangle ^{\rho}$ only corresponds to partial regularity with respect to the last variable, see the precise definition of decomposition \eqref{eq:decompo_Grushin} in Section~\ref{sec:preliminaries} and Remark~\ref{rem:Xk_rho} for details.

In order to establish Theorem~\ref{thm:main_short}, we prove the following more specific theorem.
\begin{theorembisH}[Local Cauchy Theory for~\eqref{eq:NLSG}, precised]\label{th.main} 
Let $k\in(1,\frac{3}{2}]$.  
\begin{enumerate}[label=(\textit{\roman*})]
\item  Let $u_0\in \mathcal{X}^k_{1}\subset H^k_G$. For any $\ell \in (\frac{3}{2},k+\frac{1}{2})$, for almost-every $\omega\in\Omega$, there exists $T>0$ and a unique local solution with initial data $u_0^{\omega}$ to~\eqref{eq:NLSG} in the space
\[
    e^{it\Delta_G}u_0^{\omega}+ \mathcal{C}^0([0,T),H^{\ell}_G) \subset \mathcal{C}^0([0,T),H^{k}_G)\,.
\]
More precisely, there exists $c>0$ such that for all $R\geq 1$, outside a set of probability at most $e^{-cR^2}$, one can choose $T\geq (R\|u_0\|_{\mathcal{X}^k_1})^{-2}$.
 
\item (Non-smoothing under randomization) If $u_0\in H^k_G \setminus (\bigcup_{\varepsilon >0}H^{k+\varepsilon}_G)$, then
\begin{equation*}
    \operatorname{supp}(\mu_{u_0})\subset H^k_G \setminus (\bigcup_{\varepsilon >0}H^{k+\varepsilon}_G)\,.
\end{equation*}
\item
 (Density of measures with rough potentials) Let $u_0\in H^k_G$ and $\varepsilon>0$, then there exists $v_0\in \mathcal{X}^k_1 \setminus (\bigcup_{\varepsilon' >0}H^{k+\varepsilon'}_G)$ such that
\begin{equation*}
    \operatorname{supp}(\mu_{v_0}) \cap  B_{H^k_G}(u_0,\varepsilon) \neq \varnothing\,.
\end{equation*}
\end{enumerate}
\end{theorembisH}

\begin{remark}[Consequences of (\textit{i})]
For every $k>1+2\varepsilon>1$ and $u_0\in \mathcal{X}^k_{1+2\varepsilon+1-k}$, the continuous embedding $\mathcal{X}^k_{1+2\varepsilon+1-k} \hookrightarrow \mathcal{X}^{1+2\varepsilon}_{1}$ implies that for almost-every $\omega\in\Omega$, the initial data $u_0^{\omega}$ gives rise to a unique local solution 
\[
    u \in e^{it \Delta _G}u_0^{\omega} + \mathcal{C}^0([0,T),H_G^{\frac{3}{2}+\varepsilon})\,,
\]
where we check that $\ell=\frac{3}{2}+\varepsilon\in(\frac{3}{2},1+2\varepsilon+\frac{1}{2})$.

Therefore, in the case $k=\frac{3}{2}$, we observe that the limiting almost-sure well-posedness space is $\bigcap_{\varepsilon>0} \mathcal{X}^{\frac{3}{2}}_{\frac{1}{2}+\varepsilon}$. We recall that for $k= \frac{3}{2}+\varepsilon$, local well-posedness is known to hold in $H_G^{\frac{3}{2}+\varepsilon}=\mathcal{X}_{0}^{\frac{3}{2}+\varepsilon}$. It is interesting to note that our approach looses an exponent $\rho=\frac{1}{2}$, since we do not recover the same limit space in the limit $k \to \frac{3}{2}$. 
\end{remark}

\begin{remark}[Decomposition of the solution]
In (\textit{i}), we claim that it is possible to construct local solutions to~\eqref{eq:NLSG} in the space $\mathcal{C}([0,T),H^k_G)$ for small values of $k$. However, uniqueness holds only on a smaller subset, as a consequence of an \textit{a priori} decomposition of the solution as sum of the solution to the linear equation \eqref{eq:LS-G} with initial data $u_0^{\omega}$ and a smoother part. This decomposition 
can be interpreted as a more \textit{nonlinear} decomposition of the solution than seeking $v \in H^k_G$, as we seek for $u$ in an affine space $e^{it\Delta_G}u_0^{\omega} + H^{\ell}_G$ instead of a vector space. Such a decomposition is the simplest nonlinear decomposition, akin to~\cite{burqTzvetkov1}, and is a key feature in most random data well-posedness works. Such a decomposition can be traced back to the seminal paper~\cite{bourgain2}, whose method we follow. More nonlinear decompositions for the solutions to some dispersive equations have been exhibited in~\cite{ohTzvetkovWang20}, and more recently in~\cite{bringmann21} and~\cite{dengNahmodYue19,dengNahmodYue21} (see also~\cite{gubinelliImkellerPerkowski15,hairer13,hairer14} in the context of stochastic equations), and this would be an interesting problem whether these approaches could give even more insight on the probabilistic well-posedness theory in our context.
\end{remark}

\begin{remark}[Regularity and density of the measures]
Parts (\textit{ii}) and (\textit{iii}) give regularity properties of the measures $\mu_{u_0}$. In fact, (\textit{ii}) ensures that the measure $\mu_{u_0}$ does not charge solutions which are more regular than $u_0$. 
Actually, the measure $\mu_{u_0}$ charges solutions that have regularity $W^{k+1/4,4}_G$ on $L^p$ based Sobolev spaces, but no better regularity bound is expected to hold (see Proposition~\ref{prop.random.integrability.improvement} and Remark~\ref{rk:wkp}), so this estimate alone is not enough to establish Theorem~\ref{th.main}.

Statement (\textit{iii}) goes even further since 
we prove that
\[
\bigcup\Big\{ \operatorname{supp}(\mu_{v_0})\mid v_0\in \mathcal{X}^k_1 \setminus (\bigcup_{\varepsilon >0}H^{k+\varepsilon}_G)\Big\}
\] is dense in $H^k_G$. This result is related to the support density in the Euclidean case. Indeed, for the nonlinear Schrödinger equation on the torus, probabilistic local well-posedness holds with respect to measures which are dense in Sobolev spaces, see for instance Appendix~B of~\cite{burqTzvetkov1} and~\cite{burqTzvetkov}. 
\end{remark}

\begin{remark}[Admissible initial data]
When $u_0$ only has a finite number of modes $m$, the assumption $u_0\in \mathcal{X}^k_1$ is equivalent to the condition $u_0\in H^{k+1}_G$, but since $k+1>k_C=\frac{3}{2}$, the result is void. For this reason, our result does not extend to the nonlinear half-wave equation $\partial_tu\pm\sqrt{-\Delta}u=|u|^2u$, which also admits a similar decomposition to~\eqref{eq:decompo_Grushin} but only with a finite number of modes $m$. However, we will see in Remark~\ref{rem:Xk_rho} that the condition $u_0\in \mathcal{X}^k_{1}$ still allows a general set of low regularity initial data in our context.
\end{remark}

\begin{remark}[Defocusing case] One can replace equation~\eqref{eq:NLSG} by its defocusing variant and get the exact same local well-posedness theory. Indeed, we only address local well-posedness, which mainly depends on the order of magnitude of the nonlinearity and not its sign.
\end{remark}

\begin{remark}[Randomization] The measures $\mu_{u_0}$ are defined in~\eqref{eq:randomizationMap} with Gaussian random variables. However, most of the results in this article are stated for more general subgaussian random variables (see Definition~\ref{def:subgaussianity}), except when using the Wiener chaos estimates from Corollary~\ref{coro:WienerChaos}, which is stated only for Gaussian random variables.

Note that the randomization along on a unit scale in the variable $m$ is quite classical, as this is the variable along which we establish our bilinear estimates. However, the variable $I$ plays a different role which allows us to only take a randomization on a dyadic scale. One could compare this choice with the construction of adapted dilated cubes in~\cite{benyiOhPocovnicu15a}.
\end{remark}

\subsection{Deterministic and probabilistic Cauchy theory for~\eqref{eq:NLSH}}

As mentioned at the beginning of this introduction, the nonlinear Schrödinger equation on the Heisenberg group lacks dispersion, therefore the dispersive paradigm cannot be applied for lowering the critical well-posedness exponent below $k_C=2$ (resp. $k_C=\frac{3}{2}$ for~\eqref{eq:NLSG}) given by the Sobolev embedding. More precisely, the lack of dispersion precludes the usual way in which Strichartz estimates are proven, that is a combination of a dispersive estimate and a duality $TT^*$ argument. The result in~\cite{BahouriGerardXu2000} goes even further, as the non smoothness of the flow map for~\eqref{eq:NLSH} in $H^k$ for $k<k_C$ makes it impossible to implement a fixed point argument. 

The lack of dispersion and the lack of Strichartz estimates for the Schrödinger equation on the Heisenberg group have been recently investigated in~\cite{BahouriGallagher2020} and~\cite{bahouriBarilariGallagher}. In these works, the authors prove that there exist anisotropic Strichartz estimates~\cite{bahouriBarilariGallagher}, local in space versions of the dispersive estimates (Theorem~1 in~\cite{BahouriGallagher2020}) and local version of Strichartz estimates (Theorem~3 in~\cite{BahouriGallagher2020}). These results follow the general strategy of Fourier restriction methods for proving Strichartz estimates, dating back to Strichartz~\cite{strichartz77}, and use the Fourier analysis on the Heisenberg group~\cite{bahouriCheminDanchin18, bahouriCheminDanchin19}. We also refer to~\cite{muller90} for restriction theorems on the Heisenberg group. 

Probabilistic methods have proven to be very useful to break the scaling barrier in the context of dispersive equations. Such a study has been pioneered in~\cite{bourgain1}: the purpose is to construct global solutions for nonlinear Schrödinger equations posed on the torus, using invariant measures and a probabilistic local Cauchy theory. In~\cite{burqTzvetkov1,burqTzvetkov2}, the authors extend these results to other dispersive equations, opening the way to a very active area of research and leading to an immense body of results. 

Invariant measure methods mostly reduce their scope to compact spaces, the setting of the torus being used on many works. For non-compact spaces, the probabilistic method of~\cite{burqTzvetkov1} remains largely adaptable through the use of Gaussian random initial data. 
We refer for example to~\cite{benyiOhPocovnicu15a,benyiOhPocovnicu15b} where probabilistic local well-posedness is obtained for the nonlinear Schrödinger equations on $\mathbb{R}^d$, and  to~\cite{ohPocovnicu16,pocovnicu17} for similar results with the wave equation.  

Several works go beyond the Euclidean Laplacian. For instance, in~\cite{btz} the authors replace the standard Laplacian $-\Delta$ with a harmonic oscillator $-\Delta + x^2$ and study the local Cauchy theory for the associated nonlinear Schrödinger equation. Our work is partly inspired from this work, and also subsequent works~\cite{deng12,burqThomann20,latocca20}. Indeed, in our case, rescaled harmonic oscillators parameterized by one of the variables appear when one considers a partial Fourier transform of the equation.

We point out that although no progress had been obtained in the direction of local well-posedness up to now, traveling waves and their stability have been studied in~\cite{gassot19a, gassot19b}. 

\subsection{Outline of the proof and main arguments}

In this section, we briefly review the main ideas leading to the proof of Theorem~\ref{th.main}.

\subsubsection{General strategy for almost-sure local well-posedness} 

We follow  the probabilistic approach to the local well-posedness problem from~\cite{burqTzvetkov1} in order to study~\eqref{eq:NLSG}. We fix $u_0 \in H^k_G$, where $0<k<k_C=\frac{3}{2}$, and consider the randomization $u_0^{\omega}$ defined in \eqref{eq:randomizationMap}. We seek for solutions to~\eqref{eq:NLSG} under the form  
\[
    u(t)=e^{it\Delta_G}u_0^{\omega}+v(t)
\] 
where $v(t)$ belongs to some space $H^{\frac{3}{2}+\varepsilon}_G$, $\varepsilon>0$, on which a deterministic local well-posedness theory is known to hold. Plugging this ansatz in the Duhamel representation of~\eqref{eq:NLSG} leads to
\[
    u(t)=e^{it\Delta_G}u_0^{\omega} - i\int_0^te^{i(t-t')\Delta_G} \left(|e^{it'\Delta_G}u_0^{\omega}+v(t')|^2(e^{it'\Delta_G}u_0^{\omega}+v(t'))\right)\,\mathrm{d}t'\,,
\]
so that with the notation $z^{\omega}(t)=e^{it\Delta_G}u_0^{\omega}$,  we expect that  $v(t)=\Phi v(t)\in H^{\frac{3}{2}+\varepsilon}_G$, where
\[
\Phi v(t) =-i\int_0^te^{i(t-t')\Delta_G} \left(|z^{\omega}(t')+v(t')|^2(z^{\omega}(t')+v(t'))\right)\,\mathrm{d}t' \,.
\]
In view of the lack of Strichartz estimates, the best known bounds on $\Phi v$ are the trivial estimates (we take $t \leqslant T$ and forget time estimates, as we only give heuristic arguments)
\[
    \|\Phi v\|_{H^{\frac{3}{2}+\varepsilon}_G} \lesssim \|(v+z^{\omega})^3\|_{H^{\frac{3}{2}+\varepsilon}_G} \lesssim \|v^3\|_{H^{\frac{3}{2}+\varepsilon}_G} + \|(z^{\omega})^3\|_{H^{\frac{3}{2}+\varepsilon}_G} + \|z^{\omega}v^2\|_{H^{\frac{3}{2}+\varepsilon}_G} + \|(z^{\omega})^2v\|_{H^{\frac{3}{2}+\varepsilon}_G}\,.
\]
The term $v^3$ is handled using the algebra property of the space $ H^{\frac{3}{2}+\varepsilon}_G$, since $v$ has high regularity. The terms involving $z^{\omega}$ are mode difficult because $z^{\omega}\in H^k_G \setminus H^{\frac{3}{2}+\varepsilon}_G$ only has $H^k$ regularity since the randomization does not gain derivatives, as stated in Theorem~\ref{th.main} (\textit{ii}) and proven in Section~\ref{sec:randomLinearEstimates}. A first approach would be to estimate
\[
    \|(z^{\omega})^3\|_{H_G^{\frac{3}{2}+\varepsilon}} \sim \|\langle - \Delta_G \rangle^{\frac{3}{4}+\frac{\varepsilon}{2}} z^{\omega} (z^{\omega})^2\|_{L^2_G} \lesssim \|z^{\omega}\|_{W^{\frac{3}{2},\infty}_G}\|z^{\omega}\|_{L^4_G}^2\,,
\]
thus it is important to study the effect of the randomization on $u_0$ in terms of regularity in $L^p$ based Sobolev spaces. As proven in Proposition~\ref{prop.random.integrability.improvement}, linear estimates in $W^{k, p}_G$ spaces gain up to $\frac{1}{4}$ derivatives. However, the linear estimates alone are not sufficient to gain the $\frac{1}{2}$ derivatives in regularity and therefore deal with low values of $k$ in Theorem~\ref{th.main}.

In order to improve our estimates, we establish bilinear estimates: we prove that given the random solutions $z^{\omega}$ with initial data in $H^k_G$ to the linear Schrödinger equation 
\begin{equation}
    \label{eq:LS-G}
    \tag{LS-G}
    i\partial_tz - \Delta_G z = 0\,,
\end{equation}
then almost-surely we have $(z^{\omega})^2 \in H^{k+\frac{1}{2}}_G$.  In this article, we prove the following bilinear and trilinear estimates in the spaces $\mathcal{X}^k_{\rho}$, which could be of independent interest.

\begin{theorembis}[Bilinear and trilinear estimates for random solutions]\label{th.bilinear_estimate}
Let $k\in(1,\frac{3}{2}]$ and $u_0\in \mathcal{X}^k_{1}\subset H^k_G$ (see~\eqref{def:Xk_rho}). 
Let $u_0^{\omega}$ as in \eqref{eq:randomizationMap}, and let us denote by $z^{\omega}=e^{it\Delta_G}u_0^{\omega}\in\mathcal{C}^0(\mathbb{R},H^k_G)$ the solution to~\eqref{eq:LS-G} associated to $u_0^{\omega}$. Then there exists $c>0$ such that the following statements hold.  Fix $q\in[2, \infty)$. For $T>0$, denote $L^q_T:= L^q([0,T])$.
\begin{enumerate}[label=(\textit{\roman*})]
    \item For all $R\geq 1$ and $T>0$, outside a set of probability at most $e^{-cR^2}$, one has 
    \begin{equation}
        \label{eq:Bilin_zz}
        \|(z^{\omega})^2\|_{L^q_TH^{k+\frac{1}{2}}_G} 
        +\||z^{\omega}|^2\|_{L^q_TH^{k+\frac{1}{2}}_G} \leqslant R^2T^{\frac{1}{q}} \|u_0\|_{\mathcal{X}^k_1}^2\,,
    \end{equation}
    \begin{equation}
        \label{eq:Bilin_zzz*}
        \||z^{\omega}|^2z^{\omega}\|_{L^q_TH^{k+\frac{1}{2}}_G} \leqslant R^3T^{\frac{1}{q}} \|u_0\|_{\mathcal{X}^k_1}^3\,.
    \end{equation}
    \item
     We further require that $u_0 \in \mathcal{X}^{k}_{1+\varepsilon_0}$ for some $\varepsilon_0>0$. 
      Let $\ell<k+\frac{1}{2}$. For all $R\geq 1$ and $T>0$, there exists a set $E_{R,T}$ of probability at least $1-e^{-cR^2}$ such that the following holds. Fix $\omega\in E_{R,T}$ and $v,w \in L^{\infty}_TH^{\ell}_G$, then
    \begin{equation}
        \label{eq:Bilin_zv}
        \|z^{\omega}vw\|_{L^q_TH^{\ell}_G} \leqslant RT^{\frac{1}{q}} \|u_0\|_{\mathcal{X}^k_{1+\varepsilon_0}}\|v\|_{L^{\infty}_TH^{\ell}_G}\|w\|_{L^{\infty}_TH^{\ell}_G}\,.
    \end{equation}
    Note that $v$ and $w$ may depend on $\omega$.
\end{enumerate}
\end{theorembis}

\begin{remark}
The time variable does not play an important role. Indeed, in the course of the proof, we establish deterministic pointwise estimates in the time variable, and the $L^q_T$ norm instead of $L^{\infty}_T$ only appears in order to apply the Khinchine inequality and Wiener chaos estimates from part~\ref{subsec:probabilistic_preliminaries}. In comparison, bilinear smoothing estimates for the nonlinear Schrödinger equation on $\mathbb{R}^2$ crucially exploit the time variable, as the smoothing occurs in time averages.
\end{remark}



The heuristic explained above for proving Theorem~\ref{th.main} by using Theorem~\ref{th.bilinear_estimate} is implemented rigorously in Section~\ref{sec:localWellPosedness}.

\subsubsection{Multilinear random estimates}

The bulk of this paper aims at establishing Theorem~\ref{th.bilinear_estimate}, thus we briefly outline the main aspects of the proof. 

First, because of the random nature of the $z^{\omega}$, we use random decoupling in order to reduce the estimates to ``building block estimates'', that is estimating products $\|u_av_bw_c\|_{H^{\frac{3}{2}+\varepsilon}_G}$ for $u_a$, $v_b$ and $w_c$ obtained by restricting the Sobolev frequencies of $u$, $v$ and $w$ around the values $a$, $b$ and $c$. This reduction is a consequence of Corollary~\ref{coro:WienerChaos}. 

In the Euclidean setting, the Bernstein estimates and the Littlewood-Paley decomposition would justify the heuristics $\nabla (u_a v_b w_c) \simeq \nabla (u_a) v_bw_c+u_a\nabla(v_b)w_c +u_av_b\nabla (w_c)$ and thus reduce the analysis of $\|u_av_bw_c\|_{H^{\frac{3}{2}+\varepsilon}_G}$ to that of $\|u_av_bw_c\|_{L^2_G}$. In our case, this results still holds, but rigorous justification is more intricate and is the content of Section~\ref{sec:LaplaceAction}. 

The purpose of Section~\ref{sec:deterministicBilinear} is to prove ``building-block'' estimates of the form
\[
\|u_av_b\|_{L^2_G}
	\leqslant \frac{C}{\max\{a,b\}^{\frac{1}{2}}}\|u_a\|_{L^2_G}\|v_b\|_{L^2_G}\,.
\]
To give the main idea of the proof, it is instructing to see that in partial Fourier transform along the last variable, we can think of $u_a$ and $v_b$ as
\[
    \mathcal{F}_{y \to \eta}(u_a)(x,\eta)=f(\eta) h_m(|\eta|^{\frac{1}{2}}x) \text{ and } \mathcal{F}_{y \to \eta}(v_b)(x,\eta)=g(\eta)h_n (|\eta|^{\frac{1}{2}}x)\,,
\]
where $h_m,h_n$ are Hermite functions. Thus we can see that estimating $u_au_b$ involves estimating convolution products of rescaled Hermite functions.  The key fact is now that Hermite functions, due to their localization and normalization, enjoy bilinear estimates that are better than trivial Hölder bounds, see~\cite{btz}, relying on pointwise estimates from~\cite{kt2005}, see also~\cite{ktz}.

\subsubsection{Bilinear random-deterministic estimates} 

It turns out to be more difficult to prove a multilinear estimate on probabilistic-detemrinistic products such as $u^{\omega}vw$, where $v$ and $w$ are deterministic, which is the content of~\eqref{eq:Bilin_zv}. In this case, one should pay attention that the required set of $\omega$ constructed in Theorem~\ref{th.bilinear_estimate} (\textit{ii}) does not depend on $v$ and $w$. This precludes a direct use of decoupling offered by Corollary~\ref{coro:WienerChaos} as exploited in the proof of Theorem~B~(\textit{i}). To understand the difference between the treatment of $|z^{\omega}|^2z^{\omega}$ and $z^{\omega}vw$, remark that for example in~\eqref{eq:probaDecouplingTrilinear} the set of $\omega$ which is removed depends on the $z_{I,m}$, that is on $z$, which is fixed in Theorem~\ref{th.bilinear_estimate}. In the case of $z^{\omega}v^2$ this would remove a set of $\omega$ depending on $v$ and $w$. 

In order to circumvent such a difficulty, the idea is to apply probabilistic decoupling only on terms involving the random part $z^{\omega}$. The implementation of this strategy is carried out in Section~\ref{sec:proba-deterministic} and relies on a preparatory step introduced in Section~\ref{subsec:preparatory} aiming at splitting the analysis of $z^{\omega}vw$ into a deterministic part $\mathbf{K}\coloneqq\mathbf{K}(v,w)$ and a probabilistic part $\mathbf{J}^{\omega}\coloneqq\mathbf{J}(z^{\omega})$, which are treated in Section~\ref{subsec:deterministic} and Section~\ref{subsec:probabilistic}.    


\subsection{Further work} 
We identify two directions for future works related to the use of probabilistic methods in the study of the Cauchy problem to~\eqref{eq:NLSG} and~\eqref{eq:NLSH}. 

First, as our motivation was originally to study~\eqref{eq:NLSH}, it is natural to ask what can be said about probabilistic construction of solutions to the Schrödinger equation~\eqref{eq:NLSH} on the Heisenberg group. For the radial Heisenberg sub-laplacian and $u_0\in H^k(\mathbb{H}^1)$, there holds a decomposition similar to~\eqref{eq:decompo_Grushin} and we believe that assuming that $u_0$ belongs to a space similar to~\eqref{def:Xk_rho} should imply that the local well-posedness theory for the randomized Cauchy problem holds for almost every $u_0^{\omega}\in H^k(\mathbb{H}^1)$, $k\in(\frac{3}{2},2]$ (we recall that $k_C=2$) in the space $e^{it\Delta_{\mathbb{H}^1}}u_0^{\omega}+ \mathcal{C}^0([0,T),H^{\ell}(\mathbb{H}^1)) \subset \mathcal{C}^0([0,T),H^{k}(\mathbb{H}^1))$ with $\ell > 2$. 

The gain of $1/2$ derivative only in both the Grushin and the Heisenberg case is due to the bilinear random estimate, which is responsible for the trilinear random interaction estimate. Intuitively, in the Grushin case, the key argument lies in Lemma~\ref{lem.bilinearRescalingHermite}, in which we gain a factor $m^{1/2}$ when we estimate the product of rescaled Hermite functions $h_mh_n(\alpha\cdot)$ in $L^2$ for $n\lesssim m$. In the Heisenberg case, we would need to estimate the products of $h_{m_1}h_{n_1}(\alpha_1\cdot)$ and of  $h_{m_2}h_{n_2}(\alpha_2\cdot)$ in $L^2$, under some conditions $m_1+m_2=m$ and $n_1+n_2=n$, for $n\lesssim m$. As a consequence, two applications of Lemma~\ref{lem.bilinearRescalingHermite} only lead to a gain of $m^{1/2}$, indeed we can only ensure either $m_1\geq\frac  m2$ or $m_2\geq\frac m2$. The total gain is therefore similar to the Grushin case. 

Note that for equation~\eqref{eq:NLSH}, a randomization along every mode of the decomposition along products of Hermite functions presented for instance in~\cite{gassot19a} would not preserve the radial property. In the non radial case, one has to tackle the additional terms in the expression of the sub-Laplacian on the Heisenberg group $\mathcal{L}=\Delta_{\mathbb{H}^1}+(y\partial_x-x\partial_y)\partial_s$.

Second, as this work was partly motivated by~\cite{btz} in which the authors construct global solutions using a Gibbs invariant measure, it is natural to ask if our construction of local solutions for~\eqref{eq:NLSG} can be used to construct a global dynamics. We believe that this is a non-trivial task and most likely requires new ideas. In order to construct global solutions using a Gibbs one needs first to be able to define this Gibbs measure at least formally. In the setting of the torus~\cite{bourgain1} this is made possible by the fact that the Laplacian on $\mathbb{T}^2$ has countably many eigenfunctions, and this does not work on~$\mathbb{R}^d$ (see Proposition~3.2 in~\cite{burqThomann20}). In the setting of~\cite{btz}, such a construction is possible because the harmonic oscillator also has countably many eigenfunctions. In both cases this is the discrete decomposition of solutions on an basis of eigenfunctions that makes the Gibbs measure construction possible. In the case~\eqref{eq:NLSG}, the frequency space is $\mathbb{N}_m \times \mathbb{R}_{\eta}$ (see~\eqref{eq.strangyNorm}), which is uncountable. Moreover, even when taking a periodic variable $y\in\mathbb{T}$ so that we get a countable basis of eigenfunctions, we expect the regularity of the Gibbs measure to be really low (below $L^2$) compared to the currently reachable exponents by probabilistic methods ($k>1$ in our context).

Finally, we mention that in~\cite{dengNahmodYue19}, the authors introduce the notion of {\it probabilistic scaling exponent}, below which the first Duhamel iterate does not converge in $H^s$, implying that any iterative probabilistic method diverges below this threshold. It would be interesting to ask whether such a notion can be extended to more general geometry, in particular non compact ones. However, the lack of a universal method of randomization for the initial data in our context leaves this question as a widely open problem. In particular, the choice of randomization in the frequency variable $I$ corresponding to the space variable $y$  impacts the probabilistic result in Theorem~\ref{thm:main_short} that we would get, and our choice of dyadic distribution is by no means the only possible choice.

\subsection*{Acknowledgments} 
The authors are grateful to the referees for their constructive remarks on
the manuscript and useful suggestions. They would like to warmly thank Hajer Bahouri, Nicolas Burq, Isabelle Gallagher, Patrick Gérard and Nikolay Tzvetkov for interesting discussions during the course of this work. They also thank Nicolas Camps for an important remark.


\section{Notation and preliminary estimates}
\label{sec:preliminaries}

The purpose of this section is twofold. First we introduce decomposition~\eqref{eq:decompo_Grushin}, which is due to the structure of the Grushin operator. 
Then, we  recall some useful estimates, such as Sobolev embeddings, product laws, eigenfunction estimates and probabilistic decoupling estimates.

We will use the notation $f\lesssim g$ to denote that there exists $C>0$ such that $f\leqslant Cg$.

\subsection{Decomposition along Hermite functions for the Grushin operator}
In this subsection we give an explicit description of the Grushin operator $\Delta_G= \partial_x^2 + x^2\partial_y^2$, acting on $L^2(\mathbb{R}^2)$.

Let us consider the orthonormal basis of $L^2(\mathbb{R})$ given by the Hermite functions $(h_m)_{m\geqslant 0}$. By definition, the Hermite functions are eigenfunctions of the harmonic oscillator: for all $m\geqslant 0$, we have
\[(-\partial _x^2+x^2)h_m=(2m+1)h_m\,.\]
Taking  the Fourier transform in $y$, with Fourier variable $\eta$, we observe that for all $\eta\in\mathbb{R}$, we have
\[(-\partial _x^2 + x^2\eta ^2)h_m(|\eta|^{\frac{1}{2}}x)=(2m+1)|\eta|h_m(|\eta|^{\frac{1}{2}}x)\,.\]
Therefore, one can decompose the Fourier transform $\mathcal{F}_{y \to \eta}(u)(\cdot,\eta)$ of $u\in H^k_G$ along the basis $(h_m(|\eta|^{\frac{1}{2}}\cdot))_{m\geqslant 0}$, so $u$ becomes a sum
\[
\mathcal{F}_{y \to \eta}(u)(x,\eta) =\sum_{m\geqslant 0}f_m(\eta)h_m(|\eta|^{\frac 12}x)\,.
\] 
Moreover, this decomposition is invariant by the action of the Grushin operator $-\Delta_G$, so that we can explicitly write the $H^k_G$ norm as
\begin{equation}
    \label{eq.strangyNorm}
    \|u\|^2_{H^k_G} := \|(\operatorname{Id}-\Delta_G)^{\frac{k}{2}}u\|_{L^2_G}^2=\sum_{m\geqslant0} \int_{\mathbb{R}} (1+(2m+1)|\eta|)^k|f_m(\eta)|^2|\eta|^{-\frac{1}{2}}\,\mathrm{d}\eta\,.
\end{equation}

\begin{remark} The quantity $(2m+1)|\eta|$ plays the role of ``taking two'' derivatives, that is a similar role as the Fourier variable $|\xi|^2$ in the context of the Euclidean Sobolev spaces. Keep in mind that however this quantity mixes the Hermite modes $m$ with the partial Fourier variable $\eta$.   
\end{remark}

\begin{remark} In~\eqref{eq.strangyNorm}, the extra factor $|\eta|^{-\frac{1}{2}}$ should be understood as a normalization factor. Indeed, the  $L^2_x$ norm of the function $x \mapsto h_m(|\eta|^{\frac{1}{2}}x)$ is $|\eta|^{-\frac 14}$.
\end{remark}

In order to deal with the Sobolev norms, we further decompose $u$ according to its regularity in the $y$ variable as~\eqref{eq:decompo_Grushin}
\[
  u=\sum_{(I,m)\in 2^{\mathbb{Z}}\times\mathbb{N}}u_{I,m}\,.
\]
The definition for the $u_{I,m}$ is the following (see also~\cite{martiniSikora12} for an analogous decomposition). Taking the Fourier transform along the $y$ variable, the support of $\mathcal{F}_{y\to \eta}(u_{I,m})(x,\eta)$ satisfies the condition $|\eta|\in[I,2I]$ for some dyadic relative integer $I\in 2^{\mathbb{Z}}$:
\begin{equation}
    \label{def:BuildingBlock}
    \mathcal{F}_{y\to \eta}(u_{I,m})(x,\eta)=f_m(\eta)h_m(|\eta|^{\frac{1}{2}}x){\bf 1}_{|\eta|\in[I,2I]}\,.
\end{equation}

When using the bilinear estimates, it will be useful to regroup the global frequencies $1+(2m+1)|\eta|$ in dyadic blocs $1+(2m+1)I\in [A,2A]$ where $A\in 2^{\mathbb{N}}$ is a dyadic integer. 
For the shortness of notation, we will write $(m+1)I\sim A$ instead of $1+(2m+1)I\in [A,2A]$. Therefore, we denote
\begin{equation}\label{def:dyadic}
    u_A:=\sum_{\substack{(I,m)\in 2^{\mathbb{Z}}\times\mathbb{N}\\(m+1)I\sim A}}u_{I,m}
\end{equation}
so that
\begin{equation*}
    u=\sum_{A\in2^{\mathbb{N}}}u_A\,.
\end{equation*}
 It is useful to note that, writing $\langle I\rangle=\sqrt{1+I^2}$, we have $\frac{1}{2m+1}\lesssim \frac{\langle I\rangle}{1+(2m+1)I}$. 

Because of the orthogonality of the $h_{m}$, we have the following useful identities, which we will refer to as using \textit{orthogonality}. 
\begin{lemma}
For all $k\in\mathbb{R}$, there holds
\[
    \|u\|_{H^k_G}^2 
     = \sum_{(I,m)\in 2^{\mathbb{Z}}\times\mathbb{N}}\|u_{I,m}\|^2_{H^k_G}
      = \sum_{A\in 2^{\mathbb{N}}} \|u_{A}\|^2_{H^k_G} \,.
\] 
\end{lemma}

\begin{proof} Using~\eqref{eq.strangyNorm} we have
\begin{align*}
    \|u\|^2_{H^k_G}& =\sum_{m\geqslant 0}\int_{\mathbb{R}}(1+(2m+1)|\eta|)^k|f_m(\eta)|^2|\eta|^{-\frac{1}{2}}\,\mathrm{d}\eta \\
    &= \sum_{(I,m)\in 2^{\mathbb{Z}}\times\mathbb{N}} \int_{\mathbb{R}} \left(\frac{(1+(2m+1)|\eta|)^{\frac{k}{2}}|f_m(\eta)|}{|\eta|^{\frac{1}{4}}}\mathbf{1}_{\eta \in [I,2I]}\right)^2\,\mathrm{d}\eta\,,
\end{align*}
and since $\||\eta|^{\frac{1}{4}}h_m(|\eta|^{\frac{1}{2}} \cdot )\|_{L^2_x}^2=1$ we can write 
\begin{multline*}
    \int_{\mathbb{R}} \left(\frac{(1+(2m+1)|\eta|)^{\frac{k}{2}}|f_m(\eta)|}{|\eta|^{\frac{1}{4}}}\mathbf{1}_{\eta \in [I,2I]}\right)^2\,\mathrm{d}\eta \\ 
    = \int_{\mathbb{R}^2}\left((1+(2m+1)|\eta|)^{\frac{k}{2}}|f_m(\eta)|\mathbf{1}_{\eta \in [I,2I]}h_m(|\eta|^{\frac{1}{2}}x)\right)^2\,\mathrm{d}x\mathrm{d}\eta\,.
\end{multline*}
Thanks to the definition of the $H^k_G$ norm and the Fourier-Plancherel theorem we finally arrive at $\|u\|_{H^k_G}^2 = \sum_{(I,m)\in 2^{\mathbb{Z}}\times\mathbb{N}}\|u_{I,m}\|^2_{H^k_G}$. 

The proof of the equality $\|u_A\|_{H^k_G}^2=\sum_{(m+1)I\sim A} \|u_{m,I}\|^2_{H^k_G}$ is similar.  
\end{proof}

Observe that on the support of $u_{I,m}$, we have 
\[1+(2m+1)|\eta|\in[1+(2m+1)I,1+(2m+1)2I]\subset[A,4A]\,,
\]
so that
\[
    \|u_{I,m}\|_{H^k_G} \sim (1+(2m+1)I)^{\frac{k}{2}}\|u_{I,m}\|_{L^2_G}\sim A^{\frac{k}{2}}\|u_{I,m}\|_{L^2_G}\,.
\]
Using orthogonality, one also has: for any $A \in 2^{\mathbb{N}}$,
\[
    \|u_{A}\|_{H^k_G} \sim A^{\frac{k}{2}}\|u_A\|_{L^2_G}\,.
\]

\begin{remark}\label{rem:Xk_rho}
With the above notation, one can interpret the norm~\eqref{def:Xk_rho} in $\mathcal{X}^k_\rho$
\begin{equation*}
    \|u\|_{\mathcal{X}^k_{\rho}} ^2:= \sum_{(I,m)\in 2^{\mathbb{Z}}\times\mathbb{N}} (1+(2m+1)I)^k \left\langle I \right\rangle ^{\rho} \|u_{I,m}\|^2_{L^2_G}
\end{equation*}
as
\begin{equation*}
    \|u\|^2_{\mathcal{X}^k_\rho} \sim \sum_{m\geqslant0} \int_{\mathbb{R}} (1+(2m+1)|\eta|)^{k}(1+|\eta|)^{\rho} |f_m(\eta)|^2|\eta|^{-\frac{1}{2}}\,\mathrm{d}\eta\,.
\end{equation*}
In particular, every function $u\in H^k_G$ with additional partial regularity of $\frac{\rho}{2}$ in the $y$ variable belongs to $\mathcal{X}^k_\rho$.
\end{remark}

\subsection{Hermite functions}

Let us first recall some pointwise bounds for the Hermite functions $(h_m)_{m \geqslant 0}$. We denote by $\lambda_m=\sqrt{2m+1}$ the square root of the $m$-th eigenvalue for the harmonic oscillator.

\begin{theorem}[Pointwise estimates for Hermite functions~\cite{kt2005}, Lemma~5.1]\label{th.hermiteBounds} For any $m\geqslant 0$ and $x\in\mathbb{R}$, there holds
\[|h_m(x)| \lesssim \left\{
\begin{array}{ccc}
    \frac{1}{|\lambda_m^2-x^2|^{1/4}} & \text{ if }& |x|<\lambda_m-\lambda_m^{-1/3} \\
    \lambda_m^{-1/6} & \text{ if }& ||x|-\lambda_m| \leq \lambda_m^{-1/3}\\
    \frac{e^{-s_m(x)}}{|\lambda_m^2-x^2|^{1/4}}  &\text{ if }&  |x|>\lambda_m+\lambda_m^{-1/3}\,,
\end{array}
\right.\]
where
\[
s_m(x)=\int_{\lambda_m}^x\sqrt{t^2-\lambda_m^2}\mathrm{d}t\,.
\]
\end{theorem}  

\begin{remark} In order to understand how bilinear estimates on $h_mh_n$ are proven 
(see~\cite{btz}), one may roughly picture $h_m$ to be concentrated on $[-\sqrt{2m+1},\sqrt{2m+1}]$, and work with models of the form 
\[
    h_m(x) \sim \frac{1}{\sqrt{2}}(2m+1)^{-\frac{1}{4}}\mathbf{1}_{[-\sqrt{2m+1},\sqrt{2m+1}]}(x)\,,
\] 
as long as one does not take pointwise estimates, and as long as one does not consider $L^p$ norms for $p$ too big. 
\end{remark}

The pointwise estimates imply the following lemma on the  $L^p$ norm of the Hermite functions.
\begin{lemma}[$L^p$ norms for Hermite functions~\cite{kt2005}, Corollary~3.2]\label{lem.Hermite_Lp} For any $p \geqslant 2$ there holds uniformly in $m$
\[
\|h_m\|_{L^p(\mathbb{R})}
	\lesssim  \frac{1}{\lambda_m^{\zeta(p)}},
\]
 where $\lambda_m=\sqrt{2m+1}$ and
\begin{equation}\label{def.rho}
\zeta(p)=
\begin{cases}
\frac{1}{2}-\frac{1}{p} &\text{ if } 2\leq p\leq 4
\\
\frac{1}{6}+\frac{1}{3p} &\text{ if } 4<p\leq +\infty.
\end{cases}
\end{equation}
\end{lemma}

However, these $L^p$ norm estimates will not be sufficient for our purpose, and we will rather make use of the following simplified pointwise estimates.

\begin{corollary}[Rough pointwise estimates for Hermite functions]\label{coro.pointwiseBoundsHermite} There exists $c>0$ such that for any $m\geqslant 0$ and $x\in\mathbb{R}$,
\[
|h_m(x)| \lesssim \left\{
\begin{array}{ccc}
    \lambda_m^{-\frac{1}{2}} & \text{ if } & |x| \leqslant \frac{\lambda_m}{2} \\
     \left(\lambda_m^{\frac{2}{3}} + |x^2-\lambda_m^2|\right)^{-\frac{1}{4}}&\text{ if } & \frac{\lambda_m}{2}\leqslant |x| \leqslant 2 \lambda_m \\
    e^{-\frac{1}{8}x^2}& \text { if } & |x| \geqslant 2\lambda_m.
\end{array}
\right. 
\]
\end{corollary}

For a proof of Corollary~\ref{coro.pointwiseBoundsHermite} based on Theorem~\ref{th.hermiteBounds}, see Appendix \ref{appendix:A}.

\subsection{Sobolev spaces}

We will use on several occasions Sobolev embeddings for the Grushin operator, which correspond to the Folland-Stein embedding for the Sobolev spaces on the Heisenberg group~\cite{Folland1974}.

\begin{theorem}[Folland-Stein embedding] Let $p \in[2,\infty)$.  
\begin{enumerate}[label=(\textit{\roman*})]
    \item For $k>\frac{3}{2}$
    , then
$        H^k_G \hookrightarrow L^{\infty}_G $.
    \item For $k\leqslant \frac{3}{2}$
    , if  $\frac{1}{p} \geqslant \frac{1}{2}-\frac{k}{3}$
    , one has 
$       H^k_G \hookrightarrow L^p_G$.
\end{enumerate}
\end{theorem}

\begin{proposition}\label{lem.productLaw} Let $k>0$ and $u, v \in H^k_G$. Then 
\begin{enumerate}[label=(\textit{\roman*})]
    \item (Product rule) $\|uv\|_{H^k_G} \lesssim \|u\|_{H^k_G}\|v\|_{L^{\infty}_G} + \|u\|_{L^{\infty}_G}\|v\|_{H^k_G}$;
    \item (Algebra property) if $k>\frac{3}{2}$, $\|uv\|_{H^k_G} \lesssim \|u\|_{H^k_G}\|v\|_{H^k_G}$;
    \item (Chain rule) if $k>\frac{3}{2}$, for every $p\in\mathbb{N}^*$, $\|u^p\|_{H^k_{G}} \lesssim \|u\|^p_{H^k_G}$. 
\end{enumerate}
\end{proposition} 
More details about the proof of Proposition~\ref{lem.productLaw} can be found in Appendix~\ref{appendix:B}.


\subsection{Probabilistic preliminaries}\label{subsec:probabilistic_preliminaries}

Only basic probability notions will be used in this article. Recall that  we have fixed once and for all a probability space $(\Omega, \mathcal{A}, \mathbb{P})$, and denote $\omega \in \Omega$.    

Our main probabilistic tool is the Khinchine inequality for subgaussian random variables, and a is a multilinear version for Gaussian random variables.

\begin{definition}[Subgaussian random variables]\label{def:subgaussianity} We say that the family of independent and identically distributed complex-valued random variables $(X_{I,m})_{(I,m)\in 2^{\mathbb{Z}}\times\mathbb{N}}$ is {\it subgaussian} if there exists $c>0$ such that for all $\gamma >0$,
\begin{equation*}
    \mathbb{E}\left[e^{\gamma X_{I,m}}\right] \leqslant e^{c\gamma ^2}\,.
\end{equation*}
\end{definition}

\begin{theorem}[Khinchine inequality / Kolmogorov-Paley-Zygmund~\cite{burqTzvetkov1}, Lemma~3.1]\label{th:khinchine} Let $\mathcal{I}$ be a countable set, and let $(X_n)_{n \in \mathcal{I}}$ be a sequence of independent and identically distributed complex-valued subgaussian random variables. 
Then there exists $C>0$ such that for every complex-valued sequence $(\Psi_n)_{n \in \mathcal{I}} \in \ell^2(\mathcal{I})$ and all $r \in [2,\infty)$, one has:
    \begin{equation*}
    \left \| \sum_{n \in \mathcal{I}} \Psi_n X_n \right \|_{L^r_\Omega} \leq  C\sqrt{r}\left(\sum_{n \in \mathcal{I}} |\Psi_n|^2\right)^{\frac{1}{2}}\,.
    \end{equation*}
\end{theorem}


\begin{corollary}[Probabilistic decoupling]\label{cor:khinchine}
Let $\mathcal{I}$ be a countable set, and let $(X_n)_{n \in \mathcal{I}}$ be a sequence of independent and identically distributed complex-valued subgaussian random variables. Then there exists $c>0$ such that the following holds.

Fix a sequence $(\Psi_n)_{n\in\mathcal{I}}$ of functions of the variable $\psi\in\mathbb{R}^d$ in $L^p_{\psi}$, $p=(p_1,\dots,p_d) \in [2,\infty)^d$. Fix a countable set $\mathcal{P}$ and a partition of $\mathcal{I}$ denoted $(\mathcal{I}_k)_{k\in\mathcal{P}}$. Then there exists $R_0(p)$ large enough such that for $R\geq R_0(p)$, outside a set of probability less that $e^{-cR^2}$, there holds:
\[
    \sum_{k\in\mathcal{P}}\left\|\sum_{n\in\mathcal{I}_k}\Psi_nX_n(\omega)\right\|_{L^p_{\psi}}^2 \leqslant R^2\sum_{n \in \mathcal{I}} \|\Psi_n\|_{L^p_{\psi}}^2\,.
\]
\end{corollary}

\begin{proof} First, from the triangle inequality, we have
\[
\left\| \sum_{k\in\mathcal{P}}\left\|\sum_{n\in\mathcal{I}_k}\Psi_nX_n(\omega)\right\|_{L^p_{\psi}}^2\right\|_{L^r_{\Omega}}
	\leq \sum_{k\in\mathcal{P}}\left\|\sum_{n\in\mathcal{I}_k}\Psi_nX_n(\omega)\right\|_{L^{2r}_{\Omega}L^p_{\psi}}^2\,.
\]
Now, by the Minkowski inequality and Theorem~\ref{th:khinchine} we have for $2r \geqslant \max\{p_1,\dots,p_d\}$, for all $k\in\mathcal{P}$,
\[\left\|\sum_{n\in\mathcal{I}_k} \Psi_nX_n\right\|_{L^{2r}_{\Omega}L^p_{\psi}} \leqslant C \sqrt{r} \left\|\left(\sum_{n\in\mathcal{I}_k}|\Psi_n|^2\right)^{\frac{1}{2}}\right\|_{L^p_{\psi}}\,,\]
and applying the Minkowski inequality again,
\[\left\|\sum_{n\in\mathcal{I}_k} \Psi_nX_n\right\|_{L^{2r}_{\Omega}L^p_{\psi}} \leqslant C \sqrt{r} \left(\sum_{n\in\mathcal{I}_k}\|\Psi_n\|_{L^p_{\psi}}^2\right)^{\frac{1}{2}}\,.\]
Thus by the Markov inequality, there exists $C>0$ such that
\[\mathbb{P}\left(\sum_{k\in\mathcal{P}}\left\|\sum_{n\in\mathcal{I}_k} \Psi_nX_n\right\|_{L^p_{\psi}}^2 > R^2\sum_{n\in\mathcal{I}}\|\Psi_n\|_{L^p_{\psi}}^2\right) \leqslant \left(\frac{C\sqrt{r}}{R}\right)^r\,,\]
and the conclusion follows by optimizing in $r$, which leads to the choice $r=\frac{R^2}{4C}$.
\end{proof}

We have the following multilinear version for Gaussian variables. 

\begin{theorem}[Wiener Chaos estimates~\cite{ohThomann18}, Lemma~2.6 and~\cite{simon74}, Lemma I.18-I.22]\label{thm:WienerChaos} Let $\mathcal{I}$ be a countable set, $ \ell \geqslant 1$ an integer, and let $\Psi : \mathcal{I}^{\ell} \to \mathbb{R}$. Let $(g_n)_{n \in \mathcal{I}}$ be independent and identically distributed standard real-valued Gaussian variables. Then there exists $C(\ell)$ such that the following holds. Let 
\[
    F^{\omega}:= \sum_{(n_1, \dots, n_{\ell})\in\mathcal{I}^{\ell}}\Psi_{n_1, \dots, n_{\ell}} g_1^{\omega} \cdots g_{\ell}^{\omega}\,,
\]
and assume that $F^{\omega} \in L^2_{\Omega}$. Then one has that for any $r \geqslant 2$, 
\[
    \|F^{\omega}\|_{L^r_{\Omega}} \leq C(\ell) r^{\frac{\ell}{2}} \|F^{\omega}\|_{L^2_\Omega}\,.
\]
\end{theorem}
We now state the main consequences of this theorem and of the Markov inequality that we will use in this article. 

\begin{corollary}[Probabilistic decoupling]\label{coro:WienerChaos} Let $\mathcal{I}$ be a countable set. Let $X_n = g_n +ih_n$ complex Gaussian random variables, where the $\{g_n, h_n\}_{n\in \mathcal{I}}$ are independent and identically distributed real-valued Gaussian variables. There exists $c>0$ such that the following holds.  Let $p=(p_1,\dots,p_d) \in [2,\infty)^d$ and $(\Psi_{n,n'})_{n,n'\in\mathcal{I}}$ and $(\Psi_{n,n',n''})_{n,n',n''\in\mathcal{I}}$ be functions of the variable $\psi\in\mathbb{R}^d$ belonging to $L^p_{\psi}=L^{p_1}_{\psi_1}\dots L^{p_d}_{\psi_d}$.   Then for $R\geq R_0(p)$ large enough the following holds. 
\begin{enumerate}[label=(\textit{\roman*})]
    \item Outside a set of probability at most $e^{-cR^2}$, 
    \[\left\|\sum_{n,n' \in \mathcal{I}} \Psi_{n,n'}X_{n}\bar X_{n'}\right\|_{L^p_{\psi}}^2 \leqslant R^4 \sum_{n,n'\in \mathcal{I}} \left\|\Psi_{n,n'}\right\|_{L^p_{\psi}}^2+R^4 \left(\sum_{n\in\mathcal{I}}\|\Psi_{n,n}\|_{L^p_{\psi}}\right)^2\]
\[\left\|\sum_{n,n' \in \mathcal{I}} \Psi_{n,n'}X_{n} X_{n'}\right\|_{L^p_{\psi}}^2 \leqslant R^4 \sum_{n,n'\in \mathcal{I}} \left\|\Psi_{n,n'}\right\|_{L^p_{\psi}}^2\,.\]    
    \item We assume that $\Psi_{n,n',n''}=\Psi_{n',n,n''}$ for every $n,n',n''\in\mathcal{I}$. Outside a set of probability at most $e^{-cR^2}$,
    \[\left\|\sum_{n,n',n'' \in \mathcal{I}} \Psi_{n,n',n''}X_{n}X_{n'}\bar X_{n''}\right\|_{L^p_{\psi}}^2
     \leqslant R^6\sum_{n,n',n'' \in \mathcal{I}}\|\Psi_{n,n',n''}\|^2_{L^p_{\psi}}
    +R^6\sum_{n\in\mathcal{I}} \left(\sum_{n'\in\mathcal{I}}\|\Psi_{n,n,n'}\|_{L^p_{\psi}}\right)^2.\]
    \item We assume that $\psi=(\psi_-,\psi_+)$, $p=(p_-,p_+)$, and we relax the assumption on $p$ as $p_-=(p_1,\dots,p_{d_-}) \in [1,\infty)^{d_-}$ and $p_+=(p_{d_-+1},\dots,p_d) \in [2,\infty)^{d-d_-}$. Then outside a set of probability at most $e^{-cR^2}$,
\[
\left\|\sum_{n,n' \in \mathcal{I}} \Psi_{n,n'}X_{n}\bar X_{n'}\right\|_{L^p_{\psi}}
	\leqslant R^2\left\|\left( \sum_{n,n'\in \mathcal{I}} \left\|\Psi_{n,n'}\right\|_{L^{p_+}_{\psi_+}}^2\right)^{1/2}+\sum_{n\in\mathcal{I}}\|\Psi_{n,n}\|_{L^{p_+}_{\psi_+}}\right\|_{L^{p_-}_{\psi_-}}\,.
	\]
\end{enumerate}
\end{corollary}


\begin{proof} The proof follows from Theorem~\ref{thm:WienerChaos} by expansion, writing that $\Psi_{n,n'}=b_{n,n'}+ic_{n,n'}$ and using the independence of $g_n$ from $h_n$. We fix $r\geq \max\{p_1,\dots,p_d\}$.
 
(\textit{i}) Applying the Minkowski inequality, Theorem~\ref{thm:WienerChaos} and the Markov inequality, we get that outside a set of probability at most $e^{-cR^2}$ there holds
\[
\left\|\sum_{n,n' \in \mathcal{I}} \Psi_{n,n'}X_{n}\bar X_{n'}\right\|_{L^p_{\psi}}
 	\leqslant R^2 \left\|\sum_{n,n'\in \mathcal{I}} \Psi_{n,n'} X_n \bar X_{n'}\right\|_{L^p_{\psi}L^{2}_{\Omega}}\,.
\]
But for every $x\in\mathbb{R}^d$, we expand
\begin{align*}
\left\|\sum_{n,n'\in \mathcal{I}} \Psi_{n,n'}(\psi) X_n\bar X_{n'}\right\|_{L^2_\Omega}^2
 	= & \mathbb{E}\left[\left|\sum_{n,n'\in \mathcal{I}} \Psi_{n,n'}(\psi) X_n \bar X_{n'}\right|^2\right]\\
     = & \sum_{n_1, n_2,n_1',n_2'\in \mathcal{I} } \mathbb{E}[X_{n_1}\bar{X}_{n_1'}\bar{X}_{n_2}X_{n'_2}] \Psi_{n_1,n'_1}(\psi)\overline{ \Psi_{n_2,n'_2}(\psi)}\,.
\end{align*}
The crucial observation is then the following,
\(
    \mathbb{E}[X_{n_1}\bar X_{n_1'}\bar{X}_{n_2}X_{n'_2}]=0
\)
unless $\{n_1,n'_2\}=\{n_1',n_2\}$. Indeed, since $X\coloneqq X_{I,m}$ is a complex gaussian, we have $0=\mathbb{E}[X]=\mathbb{E}[X^2]$, and the same holds for $\bar X$. Therefore we are left with:
\[
\left\|\sum_{n,n'\in \mathcal{I}} \Psi_{n,n'}(\psi) X_n \bar X_{n'}\right\|_{L^2_\Omega}^2
    = C\sum_{n,n'\in\mathcal{I}} |\Psi_{n,n'}(\psi)|^2+C' \left(\sum_{n\in \mathcal{I}} |\Psi_{n,n}(\psi)|\right)^2 \,.
\]
Taking the norm $L^p_{\psi}$ of the square root of this term and using the Minkowski inequality lead to the first inequality. For the second inequality, the proof is the same except for the remark that \(
    \mathbb{E}[X_{n_1}{X}_{n'_1}\bar X_{n_2}\bar X_{n'_2}]=0
\)
unless $\{n_1,n_1'\}=\{n_2,n'_2\}$.

(\textit{ii}) Applying the Minkowski inequality, Theorem~\ref{thm:WienerChaos} and the Markov inequality, we get that outside a set of probability at most $e^{-cR^2}$ there holds
\[\left\|\sum_{n,n',n'' \in \mathcal{I}} \Psi_{n,n',n''}X_{n}X_{n'}\bar X_{n''}\right\|_{L^p_{\psi}} \leqslant R^3 \left\|\sum_{n,n',n'' \in \mathcal{I}} \Psi_{n,n',n''}X_{n}X_{n'}\bar X_{n''}\right\|_{L^p_{\psi}L^{2}_{\Omega}}\,.
\]
We first fix $x\in\mathbb{R}^d$ and estimate  $\left\|\sum_{n,n',n'' \in \mathcal{I}} \Psi_{n,n',n''}(\psi)X_{n}X_{n'}\bar X_{n''}\right\|_{L^{2}_{\Omega}}$ by expanding
\begin{multline*}
\left\|\sum_{n,n',n'' \in \mathcal{I}} \Psi_{n,n',n''}(\psi)X_{n}X_{n'}\bar X_{n''}\right\|_{L^{2}_{\Omega}}^2\\
     = \sum_{n_1,n'_1,n''_1,n_2,n'_2,n''_2 \in \mathcal{I}} \mathbb{E}[X_{n_1}X_{n'_1}\bar X_{n''_1}\bar X_{n_2}\bar X_{n'_2} X_{n''_2}]  \Psi_{n_1,n'_1,n''_1}(\psi)\overline{\Psi_{n_2,n'_2,n''_2}(\psi)}\,.
\end{multline*}
Since $X:= X_{n}$ is a complex gaussian, we have $0=\mathbb{E}[X]=\mathbb{E}[X^2]=\mathbb{E}[X^3]$, and the same holds for $\bar X$. Therefore, if 
\[ \mathbb{E}[X_{n_1}X_{n'_1}\bar X_{n''_1}\bar X_{n_2}\bar X_{n'_2} X_{n''_2}]\neq0,\]
then one can see by double inclusion that 
\[\{n_1,n'_1,n''_2\}= \{n''_1,n_2,n'_2\}\,.\]
More precisely, using that $\mathbb{E}[X^2\bar X]=0$, we even have that the indices $n_1,n'_1,n''_2$, are in bijection with the indices $n''_1,n_2,n'_2$. This implies that up to permutation of the pairs of indices playing the same role, namely the pair of indices $n_1$ and $n'_1$ and the pair of indices $n_2$ and $n'_2$, we are in one of the following two cases:
\begin{itemize}
\item $n_1=n_2$, $n'_1=n'_2$, $n''_1=n''_2$;
\item $n_1=n''_1$, $n'_1=n'_2$ and $n_2=n''_2$.
\end{itemize}
Therefore, outside a set of probability at most $e^{-cR^2}$, there holds 
\begin{equation*}
\left\|\sum_{n,n',n'' \in \mathcal{I}} \Psi_{n,n',n''}X_{n}X_{n'}\bar X_{n''}\right\|_{L^p_{\psi}}^2
    \leqslant R^6 \left\|\left(\sum_{n,n',n'' \in \mathcal{I}}|\Psi_{n,n',n''}|^2\right)^{1/2}\right\|^2_{L^p_{\psi}}
    +R^6\left\|\left(\sum_{n' \in \mathcal{I}} \left(\sum_{n \in \mathcal{I}}|\Psi_{n,n',n}|\right)^2\right)^{1/2}\right\|_{L^{p}_{\psi}}^2\,,
\end{equation*}
so that from the Minkowski inequality on the sum over $n,n',n''$ resp.\@ $n'$,
\begin{equation*}
\left\|\sum_{n,n',n'' \in \mathcal{I}} \Psi_{n,n',n''}X_{n}{X}_{n'}\bar X_{n''}\right\|_{L^p_{\psi}}^2
    \leqslant R^6 \sum_{n,n',n'' \in \mathcal{I}}\|\Psi_{n,n',n''}\|^2_{L^p_{\psi}}\\
    +R^6\sum_{n' \in \mathcal{I}} \left\|\sum_{n \in \mathcal{I}}|\Psi_{n,n',n}|\right\|_{L^p_{\psi}}^2\,,
\end{equation*}
and from the triangle inequality on the sum over $n$,
\begin{equation*}
\left\|\sum_{n,n',n'' \in \mathcal{I}} \Psi_{n,n',n''}X_{n}{X}_{n'}\bar X_{n''}\right\|_{L^p_{\psi}}^2
    \leqslant R^6 \sum_{n,n',n'' \in \mathcal{I}}\|\Psi_{n,n',n''}\|^2_{L^p_{\psi}}\\
    +R^6\sum_{n' \in \mathcal{I}} \left(\sum_{n \in \mathcal{I}}\|\Psi_{n,n',n}\|_{L^p_{\psi}}\right)^2\,.
\end{equation*}

(\textit{iii})  Applying the Minkowski inequality, Theorem~\ref{thm:WienerChaos} and the Markow inequality, we have 
that outside a set of probability at most $e^{-cR^2}$ there holds
\[
\left\|\sum_{n,n' \in \mathcal{I}} \Psi_{n,n'}X_{n}\bar X_{n'}\right\|_{L^p_{\psi}}
	\leqslant R^2\left\|\sum_{n,n' \in \mathcal{I}} \Psi_{n,n'}X_{n}\bar X_{n'}\right\|_{L^p_{\psi}L^2_{\Omega}}\,.
\]
For $\psi\in\mathbb{R}^d$, we expand
\begin{align*}
\left\|\sum_{n,n' \in \mathcal{I}} \Psi_{n,n'}(\psi)X_{n}\bar X_{n'}\right\|_{L^2_\Omega}^2
	= &\mathbb{E}\left[\left\vert\sum_{n,n' \in \mathcal{I}} \Psi_{n,n'}(\psi)X_{n}\bar X_{n'}\right\vert^2\right]\\
     =& \sum_{n_1,n_2,n'_1,n'_2\in\mathcal{I}} \mathbb{E}[X_{n_1}\bar X_{n'_1}\bar X_{n_2} X_{n'_2}] \Psi_{n_1,n'_1}(\psi)\overline{\Psi_{n_2,n'_2}}(\psi)\,.
\end{align*}
Since $X:= X_{n}$ is a complex gaussian, we have $0=\mathbb{E}[X]=\mathbb{E}[X^2]$, and the same holds for $\bar X$. Therefore, if 
\(
\mathbb{E}[X_{n_1}\bar X_{n'_1}\bar X_{n_2} X_{n'_2}] \neq 0\,,
\)
then one can see that either $\Big(n_1=n'_1$ and $n_2=n'_2\Big)$ or $\Big(n_1=n_2$ and $n'_1=n'_2\Big)$.
Therefore, we have 
\[
\left\|\sum_{n,n' \in \mathcal{I}} \Psi_{n,n'}X_{n}\bar X_{n'}\right\|_{L^p_{\psi}}
    \leqslant R^2\left\|\left(\sum_{n,n'\in\mathcal{I}} |\Psi_{n,n'}|^2\right)^{1/2}
	+   \sum_{n\in\mathcal{I}}|\Psi_{n,n}|\right\|_{L^{p}_{\psi}}\,.
\]
Moreover, applying the Minkowski inequality for the admissible exponents $p_+\in [2,\infty)^{d-d_-}$, we obtain that outside a set of probability at most $e^{-cR^2}$ there holds
\[
\left\|\sum_{n,n' \in \mathcal{I}} \Psi_{n,n'}X_{n}\bar X_{n'}\right\|_{L^p_{\psi}}
	\leqslant R^2\left\|\left( \sum_{n,n'\in \mathcal{I}} \left\|\Psi_{n,n'}\right\|_{L^{p_+}_{\psi_+}}^2\right)^{1/2}+\sum_{n\in\mathcal{I}}\|\Psi_{n,n}\|_{L^{p_+}_{\psi_+}}\right\|_{L^{p_-}_{\psi_-}}\,.\qedhere
\]
\end{proof}





\section{Random data and linear random estimates}\label{sec:randomLinearEstimates}

In this section, we study in detail how the multiplication of each mode by independent normal distributions improves the integrability of the potential in the $L^p$ spaces, without changing the Sobolev regularity.


Let us recall the construction of the random linear solutions associated to some fixed initial data $u_0\in\mathcal{X}^k_\rho\subset H^k_G$. From~\eqref{def:Xk_rho}, the $u_{I,m}$ defined by~\eqref{def:BuildingBlock} satisfy 
\begin{equation*}
    \|u_0\|_{\mathcal{X}_{\rho}^k}^2=\sum_{(I,m)\in 2^{\mathbb{Z}}\times\mathbb{N}} (1+(2m+1)I)^k\langle I\rangle^{\rho}\|u_{I,m}\|^2_{L^2_G} \,.
\end{equation*}

Given a family of Gaussian independent and identically distributed random variables denoted $(X_{I,m})_{(I,m)\in 2^{\mathbb{Z}}\times\mathbb{N}}$, the probability measure $\mu_{u_0}$ is the push-forward of $\mathbb{P}$ under the map~\eqref{eq:randomizationMap}
\begin{equation*}
    \omega \mapsto u_0^{\omega} := \sum_{(I,m)\in 2^{\mathbb{Z}}\times\mathbb{N}} X_{I,m}(\omega)u_{I,m}\,.
\end{equation*}
We denote $z^{\omega}(t)=e^{it\Delta_G}u_0^{\omega}$ the solution to the linear equation~\eqref{eq:LS-G} associated to the initial data $u_0^{\omega}$. In particular, we have
\begin{equation}\label{eq:decomposition_Grushin_random}
z^{\omega}=\sum_{(I,m)\in 2^{\mathbb{Z}}\times\mathbb{N}} X_{I,m}(\omega)z_{I,m}\,,
\end{equation}
with 
$z_{I,m}(t)=e^{it\Delta_G}u_{I,m}$. 


\subsection{Probabilistic integrability and smoothing estimates}


In this part, we prove that the randomized potential $z^{\omega}$ belongs to the space $L^q_TL^p_G$ for every $p,q\in [2,\infty)$, and even to the space $L^q_TW^{k+\zeta(p),p}_G$, $\zeta(p)\leq \frac{1}{4}$ being defined in~\eqref{def.rho}. 

\begin{proposition}[Integrability improvement]\label{prop.random.integrability.improvement} 
There exists $c>0$ such that the following holds. Let $k\in\mathbb{R}$. Fix $u_0 \in \mathcal{X}^{k}_{\zeta (p)+\frac{1}{2}-\frac{1}{p}}$, denote $u_0^{\omega}$ its randomization from~\eqref{eq:randomizationMap}, and write $z^{\omega}=e^{it\Delta_G}u_0^{\omega}$.
Let $p,q \in [2,\infty)$ and $\varepsilon>0$. Then outside a set of probability at most $e^{-cR^2}$, there holds
\begin{equation}
    \label{eq.largeDeviation}
    \|z^{\omega}\|_{L^q_TW_G^{k+\zeta(p),p}}\leqslant RT^{\frac{1}{q}}\|u_0\|_{\mathcal{X}^{k}_{\zeta(p)+\frac{3}{2}-\frac{3}{p}}}\,,
\end{equation}
\begin{equation}
    \label{eq.largeDeviation_dyadic}
    \sum_{A\in 2^{\mathbb{N}}}A^{k+\zeta(p)-\varepsilon}\|z^{\omega}_A\|_{L^q_TW^{\varepsilon,p}_G}^2\leqslant R^2T^{\frac{2}{q}}\|u_0\|_{\mathcal{X}^{k}_{\zeta(p)+\frac{3}{2}-\frac{3}{p}}}^2\,.
\end{equation}
\end{proposition}
Note that we have $\zeta(p)+\frac{3}{2}-\frac{3}{p}=2-\frac{4}{p}$ when $p\leq 4$ and $\zeta(p)+\frac{3}{2}-\frac{3}{2p}=\frac{5}{3}-\frac{8}{3p}$ when $p>4$. The derivative gain is maximal when $\zeta (p)$ is, that is when $p=4$ leading to a $\frac{1}{4}$ gain of derivatives. 


\begin{remark}[Case $p=\infty$]
Since $\frac{4}{p} \times p = 4>3$ one can use the embedding $W_G^{\frac{4}{p},p} \hookrightarrow L^{\infty}_G$ and obtain that oustide a set of probability at most $e^{-cR^2}$, there holds
\begin{equation*}
    \|z^{\omega}\|_{L^q_TW_G^{\varepsilon,\infty}}\lesssim_{\varepsilon} RT^{\frac{1}{q}}\|u_0\|_{\mathcal{X}^{\frac{4}{p}-\zeta (p) + \varepsilon}_{\zeta (p)+\frac{3}{2}-\frac{3}{p}}}\,,
\end{equation*} 
which is obtained by taking $k=\frac{4}{p}-\zeta (p)+\varepsilon$ in~\eqref{eq.largeDeviation}. Observe that $\mathcal{X}^1_1 \hookrightarrow  \mathcal{X}^{\frac{4}{p}-\zeta (p) + \varepsilon}_{\zeta (p)+\frac{3}{2}-\frac{3}{p}}$ when $\varepsilon + \frac{1}{2} +\frac{1}{p}\leqslant 1$, satisfied as soon as $\varepsilon < \frac{1}{2}$ by choosing $p$ large enough. Therefore, for $\varepsilon \in [0,\frac{1}{2})$, outside a set of probability at most $e^{-cR^2}$, there holds
\[
\|z^{\omega}\|_{L^q_TW_G^{\varepsilon, \infty}} \lesssim_{\varepsilon} RT^{\frac{1}{q}}\|u_0\|_{\mathcal{X}^1_1}\,.
\]
Similarly, one obtains that for $\varepsilon \in (0,\frac{1}{2})$, outside a set of probability at most $e^{-cR^2}$, there holds  
\begin{equation}
    \label{eq.largeDeviation.v2}
    \sum_{A\in 2^{\mathbb{N}}}A^{\varepsilon}\|z^{\omega}_{A}\|^2_{L^q_TL^{\infty}_G}\lesssim_{\varepsilon} R^2T^{\frac{2}{q}} \|u_0\|^2_{\mathcal{X}^1_{1}}\,. 
\end{equation}
\end{remark}

The proof of Proposition~\ref{prop.random.integrability.improvement} relies on the following deterministic estimates.

\begin{lemma}\label{cor:random_integrability_improvement}
With the notation form Proposition \ref{prop.random.integrability.improvement}, let $p,q\in[2,\infty)$. 
Then 
 for all $T>0$, writing $L^q_T=L^q([0,T])$,
\begin{equation}
    \label{eq.largeDeviation_v0}
    \sum_{(I,m)\in 2^{\mathbb{Z}}\times\mathbb{N}} (1+(2m+1)I)^{k+\zeta(p)}\|z_{I,m}\|^2_{L^q_TL^p_G}\lesssim T^{\frac{2}{q}}\|u_0\|^2_{\mathcal{X}^k_{\zeta(p)+\frac{3}{2}-\frac{3}{p}}}\,.
\end{equation}
\end{lemma}


\begin{proof}[Proof of Lemma \ref{cor:random_integrability_improvement}] 
 Let $p\in [2,\infty)$. 
We apply the Hausdorff-Young inequality in the $y$ variable, so that if we denote $p'$ the conjugate exponent of $p$, we obtain that for all $t\in[0,T]$,
\begin{align*}
    \|z_{I,m}(t)\|_{L^p_G} &\lesssim \| e^{it(2m+1)|\eta|}f_{m}(\eta)h_m(|\eta|^{\frac 12}x){\bf 1}_{|\eta|\in[I,2I]}\|_{L^p_xL^{p'}_\eta} \\
        & = \|f_{m}(\eta)h_m(|\eta|^{\frac 12}x){\bf 1}_{|\eta|\in[I,2I]}\|_{L^p_xL^{p'}_\eta }\,.
\end{align*}
Moreover, using the Minkowski inequality, since $p\geqslant p'$, we have
\begin{equation*}
    \|z_{I,m}(t)\|_{L^p_G}      
     \lesssim \|f_{m}(\eta)h_m(|\eta|^{1/2}x){\bf 1}_{|\eta|\in[I,2I]}\|_{L^{p'}_\eta L^p_x}\,.
\end{equation*}
Since  $z_{I,m}(0)=u_{I,m}$, this implies
\[
\|z_{I,m}(t)\|_{L^p_G} \lesssim
	 \|u_{I,m}\|_{L^{p'}_\eta L^p_x}\,.
\]

Therefore, estimate~\eqref{eq.largeDeviation_v0} is a consequence of the inequality
\begin{equation}
    \label{eq.largeDeviation_v0noTime_proof}
    \sum_{(I,m)\in 2^{\mathbb{Z}}\times\mathbb{N}} (1+(2m+1)I)^{k+\zeta(p)}\|u_{I,m}\|^2_{L^{p'}_\eta L^p_x}\lesssim \|u_0\|^2_{\mathcal{X}^k_{\zeta(p)+\frac{3}{2}-\frac{3}{p}}}\,,
\end{equation}
which we will now establish.

Using the upper bounds in $L^p$ on the Hermite functions $h_m$ from Lemma~\ref{lem.Hermite_Lp}, and the notation~\eqref{def.rho} for the exponent $\zeta(p)$, we have
\begin{align*}
\|u_{I,m}\|_{L^{p'}_\eta L^p_x}
	&= \|f_{m}(\eta)h_m(|\eta|^{\frac 12}x){\bf 1}_{|\eta|\in[I,2I]}\|_{L^{p'}_\eta L^p_x}\\
	&\lesssim (2m+1)^{-\frac{\zeta(p)}{2}}\left\|f_m(\eta)\eta^{-\frac{1}{4}}{\bf 1}_{|\eta|\in[I,2I]}\eta^{\frac{1}{4}-\frac{1}{2p}}\right\|_{L^{p'}_\eta}\\
	&\lesssim (2m+1)^{-\frac{\zeta(p)}{2}} I^{\frac{1}{4}-\frac{1}{2p}}\|f_m(\eta)\eta^{-\frac{1}{4}}{\bf 1}_{|\eta|\in[I,2I]}\|_{L^{p'}_\eta}\,.
\end{align*}
From Hölder's inequality in $\eta$ (and because the interval $[I,2I]$ has length $I$), we get
\begin{align*}
\|u_{I,m}\|_{L^p_G}
	&\lesssim (2m+1)^{-\frac{\zeta(p)}{2}} I^{\frac{1}{4}-\frac{1}{2p}-\frac{1}{2}+\frac{1}{p'}}\|f_m(\eta)\eta^{-\frac{1}{4}}{\bf 1}_{|\eta|\in[I,2I]}\|_{L^2_\eta}\\
	&\lesssim (1+(2m+1)I)^{-\frac{\zeta(p)}{2}} \langle I\rangle^{\frac{\zeta (p)}{2}+\frac{1}{4}-\frac{1}{2p}-\frac{1}{2}+\frac{1}{p'}}\|u_{I,m}\|_{L^2_G}\,,
\end{align*}
which leads to~\eqref{eq.largeDeviation_v0noTime_proof} with $\zeta (p)+\frac{1}{2}-\frac{1}{p}-1+\frac{2}{p'}=\zeta(p)+\frac{3}{2}-\frac{3}{p}$.
\end{proof} 

\begin{proof}[Proof of Proposition \ref{prop.random.integrability.improvement}] 

Let us first establish~\eqref{eq.largeDeviation}. We start with an application of the probabilistic decoupling from Corollary~\ref{cor:khinchine}: outside a set or probability at most $e^{-cR ^2}$, we have
\begin{equation*}
\|z^{\omega}\|_{L^q_TW^{k+\zeta(p),p}_G}^2
	 \lesssim  R^2\sum_{I,m} \|z_{I,m}\|_{L^q_TW^{k+\zeta(p),p}_G}^2 \,.
\end{equation*}
Since for all $t\in\mathbb{R}$ and all $k'\in\mathbb{R}$, we have $\|(-\Delta_G)^{k'/2}z_{I,m}(t)\|_{L^p_G}^2\lesssim (1+(2m+1)I)^{k'/2}\|z_{I,m}(t)\|_{L^p_G}^2$, we deduce that
\begin{equation*}
\|z^{\omega}\|_{L^q_TW^{k+\zeta(p),p}_G}
	\lesssim R\left(\sum_{I,m} (1+(2m+1)I)^{k+\zeta(p)}\| z_{I,m}\|_{L^q_TL^p_G}^2\right)^{\frac{1}{2}}\,.
\end{equation*}
Inequality~\eqref{eq.largeDeviation} is now a consequence of~\eqref{eq.largeDeviation_v0}.

Similarly, fix $\varepsilon>0$. For $(I,m)\in 2^{\mathbb{Z}}\times\mathbb{N}$, we denote by $A$ the dyadic integer such that $(m+1)I\sim A$. Then from Corollary~\ref{cor:khinchine} applied to the partition given by the $\{(I,m) \in  2^{\mathbb{Z}}\times\mathbb{N}\mid (2m+1)I\sim A\}$ for $A\in 2^{\mathbb{N}}$. Therefore, outside a set or probability at most $e^{-cR ^2}$, we have
\begin{align*}
   \sum_{A\in 2^{\mathbb{N}}}A^{k+\zeta(p)-\varepsilon}\|z^{\omega}_A\|_{L^q_TW^{\varepsilon,p}_G}^2
   &\lesssim R^2   \sum_{I,m} A^{k+\zeta(p)-\varepsilon}\|z_{I,m}\|_{L^q_TW^{\varepsilon,p}_G}^2\\
      &\lesssim R^2   \sum_{I,m} (1+(2m+1)I)^{k+\zeta(p)}\|z_{I,m}\|_{L^q_TL^{p}_G}^2\,,
\end{align*}
enabling us to conclude thanks to~\eqref{eq.largeDeviation_v0} again.
\end{proof} 

\subsection{Non-smoothing properties of the randomization}

We now prove Theorem \ref{th.main} (\textit{ii}), stating that the randomization does not improve the Sobolev regularity in the $H^k$ spaces, in a similar spirit as in~\cite{burqTzvetkov1}. In this subsection we relax the assumption that the $(X_{I,m})_{(I,m)\in 2^{\mathbb{Z}\times\mathbb{N}}}$ are Gaussian independent and identically distributed random variables.

\begin{proposition}[Non-smoothing for random initial data]
Assume that the $(X_{I,m})_{(I,m)\in 2^{\mathbb{Z}\times\mathbb{N}}}$ are independent, identically distributed, have finite variance, zero expectation, and satisfy the condition $\mathbb{P}(X_{I,m}=0)<1$.
Let $u_0\in H^k_G \setminus ( \bigcup_{\varepsilon >0} H_G^{k+\varepsilon})$ and let $u_0^{\omega}$ defined by~\eqref{eq:randomizationMap}. Then, almost surely, we have
\(
u_0^{\omega} \in H^k_G \setminus ( \bigcup_{\varepsilon >0} H_G^{k+\varepsilon}).
\)

\end{proposition}

%

\begin{proof} Without loss of generality, we assume that $k=0$. First, let us remark that $\mathbb{E}[u_0^{\omega}]=0$. Then, since by orthogonality one has
\[
    \left\|u_0^{\omega}\right\|_{L^2_G}^2
	=\sum_{(I,m)\in 2^{\mathbb{Z}}\times\mathbb{N}}|X_{I,m}(\omega)|^2\|u_{I,m}\|_{L^2_G}^2\,,
\]
and since the $X_{I,m}$ have a finite variance, we conclude that 
\[
    \mathbb{E}\left[\left\|u_0^{\omega}\right\|_{L^2_G}^2\right]\lesssim \|u_0\|_{L^{2}_G}^2 <\infty\,.
\]

Let $\varepsilon >0$, it remains to prove that almost surely we have $u_0^{\omega}\not\in H^{\varepsilon}_G$. In fact, one only needs to show that
\begin{equation}
    \label{eq.noReg}
    \mathbb{E}\left[e^{-\|u_0^{\omega}\|^2_{H^{\varepsilon}_G}}\right] = 0\,.
\end{equation}

In order to establish~\eqref{eq.noReg}, we expand decomposition~\eqref{eq:randomizationMap} using the independence of the random variables $X_{I,m}$:
\[
    \mathbb{E}\left[e^{-\|u_0^{\omega}\|^2_{H^{\varepsilon}_G}}\right]
	= \prod_{(I,m)\in 2^{\mathbb{Z}}\times\mathbb{N}}\mathbb{E}\left[e^{-|X_{I,m}|^2\|u_{I,m}\|_{H^{\varepsilon}_G}^2}\right]\,.
\]
Since for all $I,m$, there holds $\|u_{I,m}\|_{H^{\varepsilon}_G}^2\geqslant (1+(2m+1)I)^{\varepsilon}\|u_{I,m}\|_{L^2_G}^2$, we get
\[
\mathbb{E}\left[e^{-\|u_0^{\omega}\|^2_{H^{\varepsilon}_G}}\right]
	\leqslant \prod_{(I,m)\in 2^{\mathbb{Z}}\times\mathbb{N}} \mathbb{E}\left[e^{-|X_{I,m}|^2(1+(2m+1)I)^{\varepsilon}\|u_{I,m}\|_{L^2_G}^2}\right]\,.
\]
Then we have the following alternative.

{\it First case:} Assume that the terms $(1+(2m+1)I)^{\varepsilon}\|u_{I,m}\|_{L^2_G}^2$ do not go to zero as $(I,m) \to \infty$. Then there exists $\gamma>0$ and an infinite set $S$ of indices $(I,m)\in 2^{\mathbb{Z}}\times\mathbb{N}$ satisfying $(1+(2m+1)I)^{\varepsilon} \|u_{I,m}\|_{L^2_G}^2\geqslant \gamma$. This leads to the bound
\begin{equation*}
\mathbb{E}\left[e^{-\|u_0^{\omega}\|^2_{H^{\varepsilon}_G}}\right]
	\leqslant \prod _{(I,m) \in S} \mathbb{E}\left[e^{-|X_{I,m}|^2(1+(2m+1)I)^{\varepsilon}\|u_{I,m}\|_{L^2_G}^2}\right]
	\leqslant\prod _{(I,m) \in S} \mathbb{E}\left[e^{-|X_{I,m}|^2\gamma}\right]
	=0\,,
\end{equation*}
where in the last step we used that $S$ is infinite, the fact that the $X_{I,m}$ are identically distributed and the assumption \(
    \mathbb{P}(X_{I,m} =0)<1.
\)

{\it Second case:} Assume that $(1+(2m+1)I)^{\varepsilon}\|u_{I,m}\|_{L^2_G}^2\to 0$ as $(I,m) \to \infty$. In this case, we fix $R>0$ such that
\(
    \delta_R := \mathbb{P}(|X_{I,m}|>R)>0\,,
\)
which is possible thanks to the assumption
\(
    \mathbb{P}(X_{I,m} =0)<1
\) 
on the $X_{I,m}$. Moreover, $\delta_{R}$ does not depend on $(I,m)$ since the $X_{I,m}$ are identically distributed.

For every $I,m$, we have
\[
    \mathbb{E}\left[e^{-|X_{I,m}|^2(1+(2m+1)I)^{\varepsilon}\|u_{I,m}\|_{L^2_G}^2}\right] 
	\leqslant (1-\delta_R) + \delta_R e^{-R ^2 (1+(2m+1)I)^{\varepsilon}\|u_{I,m}\|_{L^2_G}^2}\,.
\]
Since the sequence $\left((1+(2m+1)I)^{\varepsilon}\|u_{I,m}\|_{L^2_G}^2\right)_{I,m}$ is convergent to zero as $(I,m)\to\infty$, this sequence is bounded by some constant $C>0$. But on the interval $[0,R^2C]$, there holds the inequality $e^{-x}\leq 1-c_R x$ for some $c_R>0$, so that
\[
    \mathbb{E}\left[e^{-|X_{I,m}|^2(1+(2m+1)I)^{\varepsilon}\|u_{I,m}\|_{L^2_G}^2}\right]
	\lesssim 1- \delta_R c_RR^2 (1+(2m+1)I)^{\varepsilon}\|u_{I,m}\|_{L^2_G}^2\,,
\]
and the upper bound is positive. Moreover, since 
\[
    \sum_{I,m} \delta_R c_RR^2 (1+(2m+1)I)^{\varepsilon}\|u_{I,m}\|_{L^2_G}^2 
	\gtrsim \delta_R c_RR^2 \|u_0\|_{H^{\varepsilon}_G}^2=\infty\,,
\]
we conclude that the infinite product with general term $(1- \delta_R c_RR^2 (1+(2m+1)I)^{\varepsilon}\|u_{I,m}\|_{L^2_G}^2)$ is zero, so that
\[
    \prod_{(I,m)\in 2^{\mathbb{Z}}\times\mathbb{N}} \mathbb{E}\left[e^{-|X_{I,m}|^2(1+(2m+1)I)^{\varepsilon}\|u_{I,m}\|_{L^2_G}^2}\right] =0\,.
\]

In any case the above discussion proves~\eqref{eq.noReg}. 
\end{proof}

\begin{remark}\label{rk:wkp}
Note that using the fact that the decay rate for the $L^p$ norms of $h_m$ is optimal in Lemma~\ref{lem.Hermite_Lp} (see Lemma~5.1 in~\cite{kt2005}) and the probabilistic decoupling argument from Corollary~\ref{coro:WienerChaos}, it seems highly unlikely that $u_0^{\omega}$ belongs to $W^{\zeta(p),p}_G$ for $\varepsilon>0$ (see also~\cite{ImekrazRobertThomann2016}).
\end{remark}

\subsection{Density of the measure $\mu_{u_0}$}


In this part, we establish Theorem \ref{th.main} (\textit{iii}). Before we turn to the density properties of the measures $\mu_{u_0}$ associated to rough potentials $u_0$, we briefly justify that we can construct functions $u_0\in \mathcal{X}^k_1$ such that $u_0 \in H^k_G \setminus ( \bigcup_{\varepsilon >0} H_G^{k+\varepsilon})$.

\begin{lemma}[Existence of rough potentials]\label{lem:existencePotentials} Let $k \geqslant 0$ and $\rho \geqslant 0$. There exists a function $u_0 \in \mathcal{X}^{k}_{\rho}\subset H^k_G$ such that for all $\varepsilon >0$, $u_0 \notin H^{k+\varepsilon}_G$.

\end{lemma}

\begin{remark} In fact, this lemma implies that there exists uncountably many such functions. Indeed, we can apply another randomization argument to the potential $u_0$ from the lemma. Take $v_0^{\omega} = \sum_{(I,m)\in 2^{\mathbb{Z}}\times\mathbb{N}} \varepsilon_{I,m}(\omega)u_{I,m}$, where $\varepsilon_{I,m}$ are independent random signs, then the functions $v_0^{\omega}$ almost-surely satisfy the requirements.
\end{remark}

\begin{proof} Let $k, \rho \geqslant 0$. We consider $\displaystyle u_0 \coloneqq \sum_{(I,m)\in 2^{\mathbb{Z}}\times\mathbb{N}} u_{I,m}$ defined in Fourier variable  by $u_{I,m}=0$ if $I< 1$, and if $I\geq 1$, we take $u_{I,m}$ "constant" by parts
\[
    \mathcal{F}_{y\to \eta}(u_{I,m})(x,\eta) = \|u_{I,m}\|_{L^2_G}|\eta|^{\frac{1}{4}}h_m\left(|\eta|^{\frac{1}{2}}x\right)\,.
\]
where for any $(I,m)\in 2^{\mathbb{N}}\times\mathbb{N}$, we choose
\[
    \|u_{I,m}\|^2_{L^2_G}= \frac{1}{(1+(2m+1)I)^{k}\langle I\rangle^{\rho}\log(1+I)^2(m+1)\log(m+2)^2}\,.
\]

First, we observe that for $\varepsilon\in\mathbb{R}$,
\[
    \|u_0\|^2_{H^{k+\varepsilon}_G} \sim\sum_{ I\in 2^{\mathbb{N}} }\frac{1}{\langle I\rangle^{\rho}\log(1+I)^2} \sum_{m\in\mathbb{N}}\frac{(1+(2m+1)I)^{\varepsilon}}{(m+1)\log(m+2)^2}\,.
\]
This series with positive general term is divergent when $\varepsilon>0$. For $\varepsilon\leq 0$, this series is bounded by $C\sum_{ I\in 2^{\mathbb{N}} }\frac{1}{\langle I\rangle^{\rho}\log(1+I)^2}$ which is convergent.

It remains to prove that $u_0 \in \mathcal{X}^k_{\rho}$. In order to do so, we compute:
\begin{equation*}
    \|u_0\|^2_{\mathcal{X}^k_{\rho}}= \sum_{(I,m)\in 2^{\mathbb{N}}\times\mathbb{N}}  \frac{1}{\log(1+I)^2(m+1)\log(m+2)^2}\,,
\end{equation*}
which is indeed finite.
\end{proof}

We now establish density properties of the support of the measures with rough potentials (Theorem \ref{th.main} (\textit{iii})). 

\begin{proposition}[Density of measures with rough potentials] 
Let $k \geqslant 0$. 
Assume that the $(X_{I,m})_{(I,m)\in 2^{\mathbb{Z}\times\mathbb{N}}}$ are independent, identically distributed, have finite variance, and satisfy the condition $\mathbb{P}(|X_{I,m}-1|<r)>0$ for all $r>0$.
Let $u_0 \in H^k_G$ and $\varepsilon >0$. 
 Then there exists $v_0\in \mathcal{X}^k_1\setminus ( \bigcup_{\varepsilon' >0} H_G^{k+\varepsilon'})$ such that 
\[\mu_{v_0} \left(B_{H^k_G}(u_0,\varepsilon)\right)>0\,.\]
\end{proposition}

\begin{proof} Without loss of generality, we assume that $k=0$, and we only deal with the Grushin case. 
Fix $\varepsilon>0$ and $u_0\in L^2_G$. We first construct $v_0\in \mathcal{X}^0_2\setminus ( \bigcup_{\varepsilon' >0} H_G^{\varepsilon'})$ satisfying the condition $\|u_0-v_0\|_{L^2_G} \leqslant \varepsilon$, then we prove that $v_0$ meets the requirements from the proposition.


To construct $v_0$,  let $K_0 >0$ be such that 
\[
    \sum_{|(I,m)|>K} \|u_{I,m}\|^2_{L^2_G} \leqslant \varepsilon^2\,,
\]
where $|(I,m)|=\max\{|\log(1+I)|,m\}$. 
Then, let $\delta>0$ to be determined later, and denote by 
\[
\displaystyle \widetilde{u}_0=\sum_{(I,m)\in 2^{\mathbb{Z}}\times\mathbb{N}}\widetilde{u}_{I,m}\in \mathcal{X}^0_1\setminus ( \bigcup_{\varepsilon' >0} H_G^{\varepsilon'})
\]
the potential constructed in Lemma~\ref{lem:existencePotentials} for $k=0$ and $\rho=1$. Set 
\[
    v_{I,m}\coloneqq
    \begin{cases} 
        u_{I,m} \quad &\text{ if }|(I,m)| \leqslant K_0\\
        \delta\widetilde{u}_{I,m} \quad &\text{ if }|(I,m)| > K_0\,.
    \end{cases}
\] 
Since the coefficients $v_{I,m}$ for $|(I,m)| > K_0$ are the ones from Lemma~\ref{lem:existencePotentials}, we know that $v $ belongs to $ \mathcal{X}^0_{1} \setminus (\bigcup _{\varepsilon' >0}H^{\varepsilon'}_G)$. 
Then we compute that  when $\delta$ is small enough,
\[
    \|u_0-v_0\|_{L^2_G}^2 \leqslant  \sum_{|(I,m)|>K_0} \|u_{I,m}\|^2_{L^2_G}+\delta^2 \sum_{|(I,m)|>K_0} \|\widetilde{u}_{I,m}\|^2_{L^2_G}\leqslant 2\varepsilon^2\,.
\]


In the end of this proof, we establish that 
\begin{equation*}
    \mu_{v_0} \left(w \in L^2_G, \|w-v_0\|_{L^2_G} \leqslant \varepsilon\right) >0\,,
\end{equation*}
as this would imply
\begin{equation*}
    \mu_{v_0} \left(w \in L^2_G, \|w-u_0\|_{L^2_G} \leqslant 2\varepsilon\right) >0\,.
\end{equation*}
For future occurrences, for $K>0$ and $w=\sum_{I,m}w_{I,m}\in L^2_G$, let us denote 
\[
    \Pi_K( w) := \sum_{|(I,m)|\leqslant K} w_{I,m}\,.
\]
Because of the inclusion 
\begin{multline*}
    \left\{ w\in L^2_G\mid \|\Pi_K (w-v_0)\|_{L^2_G}\leqslant \frac{\varepsilon}{2}\right\}\cap \left\{w\in L^2_G\mid \|(\mathrm{Id}-\Pi_K)(w-v_0)\|_{L^2_G}\leqslant \frac{\varepsilon}{2}\right\}\\
	\subset \left\{w \in L^2_G\mid \|w-v_0\|_{L^2_G} \leqslant \varepsilon \right\} \,,
\end{multline*}
we know by independence that
\begin{multline*}
    \mu_{v_0} \left(w \in L^2_G, \|w-v_0\|_{L^2_G} \leqslant \varepsilon\right)
	\geqslant \mu_{v_0}\left(w \in L^2_G, \|\Pi_{K} (w-v_0)\|_{L^2_G}\leqslant \frac{\varepsilon}{2}\right)\\
	\times\mu_{v_0}\left(w \in L^2_G, \|(\mathrm{Id}-\Pi_{K})(w-v_0)\|_{L^2_G}\leqslant \frac{\varepsilon}{2}\right)\,.
\end{multline*}

To handle the second term in the right-hand side, since $\Pi_Kv_0$ tends to $v_0$ in $L^2_G$ as $K\to \infty$, one has that for large enough $K>0$,
\[
\mu_{v_0}\left(w \in L^2_G, \|(\mathrm{Id}-\Pi_K)(w-v_0)\|_{L^2_G}\leqslant\frac{\varepsilon}{2}\right)
	\geqslant \mu\left( w\in L^2_G, \|(\mathrm{Id}-\Pi_K)w\|_{L^2_G}\leqslant\frac{\varepsilon}{4}\right)\,.
\]
But we have
\[
    \mu_{v_0}\left(w \in L^2_G, \|(\mathrm{Id}-\Pi_K)w\|_{L^2_G}> \frac{\varepsilon}{4} \right) \underset{K\to \infty}{\longrightarrow}0\,,
\] 
since from the Markov inequality,
\begin{multline*}
    \mu_{v_0} \left(w \in L^2_G, \|(\mathrm{Id}-\Pi_K)w\|_{L^2_G}> \frac{\varepsilon}{4} \right)
    =\mathbb{P} \left( \|(\mathrm{Id}-\Pi_K)v^{\omega}\|_{L^2_G}> \frac{\varepsilon}{4} \right)\\
     \lesssim \varepsilon^{-2} \mathbb{E}\left[\|(\mathrm{Id}-\Pi_K)v^{\omega}\|^2_{L^2_G}\right]
    \lesssim \varepsilon^{-2} \sum_{|(I,m)|>K} \|v_{I,m}\|^2_{L^2_G}\,,
\end{multline*}
which goes to zero as $K$ goes to infinity. We therefore fix $K>0$ large enough so that
\[
\mu_{v_0}\left(w \in L^2_G, \|(\mathrm{Id}-\Pi_K)(w-v_0)\|_{L^2_G}\leqslant\frac{\varepsilon}{2}\right)
	>0\,.
\]

It only remains to handle the first term in the right-hand side and prove that
\begin{equation}
    \label{eq.density2}
    \mu_{v_0}\left(w \in L^2_G, \|\Pi_K(w-v_0)\|_{L^2_G}\leqslant \frac{\varepsilon}{2} \right)>0\,. 
\end{equation}
Let $c>0$ small enough so that the following inclusion holds:
\begin{multline*}
\bigcap_{|(I,m)|<K} \left\{w\in L^2_G\mid \exists \widetilde{X}_{I,m}\in\mathbb{R}, w_{I,m}=\widetilde{X}_{I,m}v_{I,m}
	 \text{ and } \|(\widetilde{X}_{I,m}-1)v_{I,m}\|_{L^2_G}\leqslant c\varepsilon\right\}\\
	\subset \left\{w \in L^2_G\mid \|\Pi_K(w-v_0)\|_{L^2_G}\leq\frac{\varepsilon}{2} \right\}\,.
\end{multline*}
This implies by independence of the random variables $X_{I,m}$ that
\begin{equation*}
 \mu_{v_0}\left(w \in L^2_G, \|\Pi_K(w-v_0)\|_{L^2_G}\leqslant\frac{\varepsilon}{2} \right)
 	\geqslant \prod_{|(I,m)|\leqslant K} \mathbb{P}\left(\|(X_{I,m}(\omega)-1)v_{I,m}\|_{L^2_G}\leq c\varepsilon \right),
\end{equation*}
which is positive thanks to the assumption $\mathbb{P}(|X_{I,m}-1|<r)>0$ for all $r>0$, thus~\eqref{eq.density2} holds. 
\end{proof}

\section{Action of the Laplace operator}\label{sec:LaplaceAction}

In the remaining of this article, our general strategy is to group the decompositions~\eqref{eq:decompo_Grushin} of $u$ and $v$ in dyadic packets~\eqref{def:dyadic}
\[
u=\sum_{A\in 2^{\mathbb{N}}}u_{A}
\quad
\text{and}
\quad
v=\sum_{B\in 2^{\mathbb{N}}}v_{B}\,,
\]
and write the estimates for the $H^k$ norm of the product $uv$ in terms of the $L^2$ norms of products $u_Av_B$ for dyadic $A$ and $B$.

In this section, we consider two elements $u,v\in H^k_G$ 
 and provide useful estimates for the $H^k$ norm of the products $u_Av_B$ and $u_Av$ in terms of the $L^2$ norms of the products $u_Av_B$ for dyadic $A$ and $B$.

\subsection{The action of $-\Delta_G$  on a product of functions}

In this part, we prove that for $u,v\in H^k_G$, the term $-\Delta_G(uv)$ can be expressed thanks to shifted versions of $u$ and $v$, defined as follows.


\begin{definition}[$\delta$-shifted functions]\label{def:delta_shifted}
Let and $u\in H^k_G$, decomposed as~\eqref{eq:decompo_Grushin}
\(
u=\sum_{(I,m)\in 2^{\mathbb{Z}}\times\mathbb{N}}u_{I,m}
\), where $u_{I,m}$ is defined in~\eqref{def:BuildingBlock} as
\begin{equation*}
    \mathcal{F}_{y\to \eta}(u_{I,m})(x,\eta)=f_{I,m}(\eta)h_m(|\eta|^{\frac{1}{2}}x)\,,
\end{equation*}
and $f_{I,m}(\eta)=f_m(\eta){\bf 1}_{|\eta|\in[I,2I]}$.
Throughout this section we use the notation
\begin{equation}
    \label{eq.D1}
    D_1 \coloneqq \{-1,0,1\} \times \{+,-\}\,,
\end{equation}
\begin{equation}
    \label{eq.D2}
    D_2= D_1 \times D_1\,. 
\end{equation}
For  $\delta\in D_1:=\{-1,0,1\}\times\{+,-\}$ we write $\delta=(\delta_0,\pm)$  and for $m\in\mathbb{N}$, we write $m+\delta$ as a shortcut for $m+\delta:= m+\delta_0$. We introduce the shifted function $u^{\delta}$ from its decomposition~\eqref{eq:decompo_Grushin}
\(
u^{\delta}=\sum_{(I,m)\in 2^{\mathbb{Z}}\times\mathbb{N}}u_{I,m}^{\delta}
\)
as follows: for all $(I,m)\in 2^{\mathbb{Z}}\times\mathbb{N}$, 
\[
    \mathcal{F}_{y\to \eta}(u_{I,m}^{\delta})(x,\eta)=F_{I,m}^{\delta}(\eta) f_{I,m}(\eta)h_{m+\delta}(|\eta|^{\frac{1}{2}}x)\,,
\]
 and if $A\in 2^{\mathbb{N}}$ is the dyadic integer  such that $(m+1)I\sim A$,
\begin{equation}\label{def:F_delta}
F_{I,m}^{\delta}(\eta):=
\begin{cases}
\left(\frac{(2m+1)|\eta|}{4A}\right)^{\frac{1}{2}} &\text{ if }\delta\in\{(-1,+),(1,+)\}
\\
\left(\frac{(2m+1)}{4A}\right)^{\frac{1}{2}}\frac{\eta}{\sqrt{|\eta|}} &\text{ if }\delta\in\{(-1,-),(1,-)\}
\\
\frac{(2m+1)|\eta|}{4A} &\text{ if }\delta=(0,+)
\\
1&\text{ if }\delta=(0,-).
\end{cases}
\end{equation}
\end{definition}

Note that by definition, for all $k\geqslant 0$ and $u\in H^k_G$, we have
\(
\|u^{\delta}\|_{H^k_G}\leqslant \|u\|_{H^k_G}.
\)

\begin{lemma}[Action of $\Delta_G$]\label{lem:shift} Let $A,B\in2^{\mathbb{N}}$, $I,I\in 2^{\mathbb{Z}}$ and $m,n\in\mathbb{N}$ such that $(m+1)I\sim A$ and $(n+1)J\sim B$. 
Then there holds:
\begin{equation*}
(\operatorname{Id}-\Delta_G)(u_{I,m}v_{J,n}) = \sum_{(\delta_1,\delta_2)\in D_2}C_{A,B}(\delta_1,\delta_2)u_{I,m}^{\delta_1} v_{J,n}^{\delta_2}\,,
\end{equation*}
where $D_2$ is defined by~\eqref{eq.D2} and where for some explicit numerical constants $c(\delta_1,\delta_2)$,
\[
C_{A,B}(\delta_1,\delta_2)=
\begin{cases}
4A &\text{ if }(\delta_1,\delta_2)=((0,-),(0,+))
\\
4B &\text{ if }(\delta_1,\delta_2)=((0,+),(0,-))
\\
c(\delta_1,\delta_2)\sqrt{AB} &\text{ if }(\delta_1,\delta_2)\in (\{-1,1\}\times\{+,-\})^2
\\
1 &\text{ if }(\delta_1,\delta_2)=((0,-),(0,-))
\\
0 &\text{ if }(\delta_1,\delta_2)=((0,+),(0,+)),.
\end{cases}
\] 
\end{lemma}

\begin{proof}
Taking the Fourier transform of $u_{I,m}v_{J,n}$ in $y$, we transform the product into a convolution product
\[
\mathcal{F}_{y \to \eta}(u_{I,m}v_{J,n})(x,\eta)
	= \int f_{I,m}(\eta_1)h_m(|\eta_1|^{\frac{1}{2}}x)g_{J,n}(\eta-\eta_1)h_n(|\eta-\eta_1|^{\frac{1}{2}}x)\mathrm{d}\eta_1.\]
In Fourier variable, the Grushin operator acts as $\partial_{xx}-x^2|\eta|^2$, therefore we have
\[
\mathcal{F}(\Delta_G(u_{I,m}v_{J,n}))(x,\eta)
	= \int f_{I,m}(\eta_1)g_{J,n}(\eta-\eta_1)  (\partial_{xx}-x^2|\eta|^2)\left(h_m(|\eta_1|^{\frac{1}{2}}x)h_n(|\eta-\eta_1|^{\frac{1}{2}}x)\right)\mathrm{d}\eta_1.\]

We now use the fact that the Hermite functions are eigenvectors of the Harmonic oscillator:
\[
\frac{\mathrm{d}^2}{\mathrm{d} x^2}h_m(x)=-(2m+1)h_m+x^2h_m
\]
to deduce the formula
\begin{align*}
(\partial_{xx}&-x^2|\eta|^2)(h_m(|\eta_1|^{\frac{1}{2}}x)h_n (|\eta-\eta_1|^{\frac{1}{2}}x))\\
	=&-((2m+1)|\eta_1|+(2n+1)|\eta-\eta_1|)h_m(|\eta_1|^{\frac{1}{2}}x)h_n(|\eta-\eta_1|^{\frac{1}{2}}x)\\
	&+x^2(|\eta_1|^2+|\eta-\eta_1|^2-|\eta|^2)h_m(|\eta_1|^{\frac{1}{2}}x)h_n(|\eta-\eta_1|^{\frac{1}{2}}x)
	+2\partial_x(h_m(|\eta_1|^{\frac{1}{2}}x))\partial_x(h_n (|\eta-\eta_1|^{\frac{1}{2}}x))\,.
\end{align*}
Note that if $\eta_1$ and $\eta-\eta_1$ have the same sign, then $|\eta_1|^2+|\eta-\eta_1|^2-|\eta|^2=-2\eta_1(\eta-\eta_1)$, and if $\eta_1$ and $\eta-\eta_1$ are of opposite sign, $|\eta_1|^2+|\eta-\eta_1|^2-|\eta|^2=2|\eta_1||\eta-\eta_1|=-2\eta_1(\eta-\eta_1)$.

We now use the identities
\begin{equation*}
xh_m(x)=\sqrt{\frac{m}{2}}h_{m-1}+\sqrt{\frac{m+1}{2}}h_{m+1}\,,
\end{equation*}
\begin{equation*}
\frac{\mathrm{d}}{\mathrm{d} x}h_m(x)=\sqrt{\frac{m}{2}}h_{m-1}-\sqrt{\frac{m+1}{2}}h_{m+1}\,,
\end{equation*}
to remove the weight $x^2$ and simplify the products $\partial_x(h_m(|\eta_1|^{\frac{1}{2}}x))\partial_x(h_n (|\eta-\eta_1|^{\frac{1}{2}}x))$. We deduce that
\begin{align*}
(\partial_{xx}-x^2|\eta|^2)(h_m(|\eta_1|^{\frac{1}{2}}x)h_n &(|\eta-\eta_1|^{\frac{1}{2}}x))\\
	=&-((2m+1)|\eta_1|+(2n+1)|\eta-\eta_1|)h_m(|\eta_1|^{\frac{1}{2}}x)h_n(|\eta-\eta_1|^{\frac{1}{2}}x)\\
	&- 2\frac{\eta_1(\eta-\eta_1)}{\sqrt{|\eta_1||\eta-\eta_1|}}\left(\sqrt{\frac{m}{2}}h_{m-1}(|\eta_1|^{\frac{1}{2}}x)+\sqrt{\frac{m+1}{2}}h_{m+1}(|\eta_1|^{\frac{1}{2}}x)\right)\\
	&\hspace{42pt}\times\left(\sqrt{\frac{n}{2}}h_{n-1}(|\eta-\eta_1|^{\frac{1}{2}}x)+\sqrt{\frac{n+1}{2}}h_{n+1}(|\eta-\eta_1|^{\frac{1}{2}}x)\right)\\
	&+2\sqrt{|\eta_1||\eta-\eta_1|}\left(\sqrt{\frac{m}{2}}h_{m-1}(|\eta_1|^{\frac{1}{2}}x)-\sqrt{\frac{m+1}{2}}h_{m+1}(|\eta_1|^{\frac{1}{2}}x)\right)\\
	&\hspace{42pt}\times\left(\sqrt{\frac{n}{2}}h_{n-1}(|\eta-\eta_1|^{\frac{1}{2}}x)-\sqrt{\frac{n+1}{2}}h_{n+1}(|\eta-\eta_1|^{\frac{1}{2}}x)\right)\,.
\end{align*}
Now, we observe that by definition of the support of $\mathcal{F}_{y\to\eta}(u_{I,m})(x,\eta)$ and $\mathcal{F}_{y\to\eta}(v_{J,n})(x,\eta)$, we have $1+(2m+1)|\eta_1|\in[A,4A]$ and $1+(2n+1)|\eta-\eta_1|\in[B,4B]$. Therefore, for some numerical constants $c(\delta_1,\delta_2)$, using the notations~\eqref{def:F_delta} for $F^{\delta}_{I,m}$,  one can write
\begin{align*}
(1-\partial_{xx}+&x^2|\eta|^2)(h_m(|\eta_1|^{\frac{1}{2}}x)h_n(|\eta-\eta_1|^{\frac{1}{2}}x))\\
	=&\left(1+4A F_{I,m}^{(0,+)}(\eta_1)+4BF_{J,n}^{(0,+)}(\eta-\eta_1)\right)h_m(|\eta_1|^{\frac{1}{2}}x)h_n(|\eta-\eta_1|^{\frac{1}{2}}x)\\
	&+\sum_{(\delta_1,\delta_2)\in (\{-1,1\}\times\{+,-\})^2}c(\delta_1,\delta_2)\sqrt{AB} F_{I,m}^{\delta_1}(\eta_1)F_{J,n}^{\delta_2}(\eta-\eta_1)	h_{m+\delta_1}(|\eta_1|^{\frac{1}{2}}x)h_{n+\delta_2}(|\eta-\eta_1|^{\frac{1}{2}}x)\,.
\end{align*}
We denote  $D_2=D_1^2$. Then we  define $C_{A,B}((0,-),(0,-))=1$, $C_{A,B}((0,+),(0,+))=0$, $C_{A,B}((0,+),(0,-))=2A$, $C_{A,B}((0,-),(0,+))=2B$ and $C_{A,B}(\delta_1,\delta_2)=c(\delta_1,\delta_2)\sqrt{AB}$ if $(\delta_1,\delta_2)\in(\{-1,1\}\times\{+,-\})^2$.
With the notation from Definition~\ref{def:delta_shifted}
(and a similar notation for the $v_{J,n}^{\delta_2}$), 
we conclude that
\begin{multline*}
\mathcal{F}((\operatorname{Id}-\Delta_G)(u_{I,m}v_{J,n}))(x,\eta)\\
	= \sum_{(\delta_1,\delta_2)\in D_2}C_{A,B}(\delta_1,\delta_2) \int f_{I,m}^{\delta_1}(\eta_1)g_{J,n}^{\delta_2}(\eta-\eta_1) h_{m+\delta_1}(|\eta_1|^{\frac{1}{2}}x)h_{n+\delta_2}(|\eta-\eta_1|^{\frac{1}{2}}x)\,\mathrm{d}\eta_1\,,
\end{multline*}
so that
\begin{equation*}
(\operatorname{Id}-\Delta_G)(u_{I,m}v_{J,n})\\
	= \sum_{(\delta_1,\delta_2)\in D_2}C_{A,B}(\delta_1,\delta_2) u_{I,m}^{\delta_1}v_{J,n}^{\delta_2}\,.\qedhere
\end{equation*}
\end{proof}

\begin{lemma}[Action of $\Delta_G$ for the product of three terms]\label{lem:shift_3terms}There exists a finite set 
\begin{equation}\label{def:D3}
D_3\subset D_1\times D_1\times D_1
\end{equation}
(where $D_1$ is defined in~\eqref{eq.D1}) such that the following holds. Let $u^{(1)},u^{(2)},u^{(3)}\in L^2_G$, $A_1,A_2,A_3\in2^{\mathbb{N}}$, $I_1,I_2,I_3\in 2^{\mathbb{Z}}$ and $m_1,m_2,m_3\in\mathbb{N}$ such that $(m_i+1)I_i\sim A_i$ for $i=1,2,3$. Then
\begin{equation*}
(\operatorname{Id}-\Delta_G)(u^{(1)}_{I_1,m_1}u^{(2)}_{I_2,m_2} u^{(3)}_{I_3,m_3}) = \sum_{\delta=(\delta_1,\delta_2,\delta_3)\in D_3}C_{A_1,A_2,A_3}(\delta)(u^{(1)}_{I_1,m_1})^{\delta_1}(u^{(2)}_{I_2,m_2})^{\delta_2}(u^{(3)}_{I_3,m_3})^{\delta_3}\,,
\end{equation*}
where the shifted functions $(u^{(i)}_{I_i,m_i})^{\delta_i}$ are defined in Definition~\ref{def:delta_shifted} and for all $\delta\in D_3$,
\[
|C_{A_1,A_2,A_3}(\delta)|\lesssim \max\{A_1,A_2,A_3\}\,.
\]
\end{lemma}

\begin{proof}
Taking the Fourier transform in $y$, we transform the product $u^{(1)}_{I_1,m_1}u^{(2)}_{I_2,m_2} u^{(3)}_{I_3,m_3}$ into a convolution product
\begin{multline*}
\mathcal{F}_{y \to \eta}(u^{(1)}_{I_1,m_1}u^{(2)}_{I_2,m_2} u^{(3)}_{I_3,m_3})(x,\eta)
	= \int f^{(1)}_{I_1,m_1}(\eta_1)f^{(2)}_{I_2,m_2}(\eta_2)f^{(3)}_{I_3,m_3}(\eta-\eta_1-\eta_2)\\
	h_{m_1}(|\eta_1|^{\frac{1}{2}}x)h_{m_2}(|\eta_2|^{\frac{1}{2}}x)h_{m_3}(|\eta-\eta_1-\eta_2|^{\frac{1}{2}}x)\mathrm{d}\eta_1\mathrm{d}\eta_2\,.
\end{multline*}
Then, we write $h_{m_1,m_2,m_3}(\eta_1,\eta_2,\eta_3):=h_{m_1}(|\eta_1|^{\frac{1}{2}}x)h_{m_2}(|\eta_2|^{\frac{1}{2}}x)h_{m_3}(|\eta_3|^{\frac{1}{2}}x)$ and expand
\begin{align*}
(\partial_{xx}-x^2|\eta_1+\eta_2+\eta_3|^2)&(h_{m_1,m_2,m_3}(\eta_1,\eta_2,\eta_3))\\
	=&-((2m_1+1)|\eta_1|+(2m_2+1)|\eta_2|+(2m_3+1)|\eta_3|)h_{m_1,m_2,m_3}(\eta_1,\eta_2,\eta_3)\\
	&+x^2(|\eta_1|^2+|\eta_2|^2+|\eta_3|^2-|\eta_1+\eta_2+\eta_3|^2)h_{m_1,m_2,m_3}(\eta_1,\eta_2,\eta_3)\\
	&+2(\mathbb{I}_{1,2}+\mathbb{I}_{1,3}+\mathbb{I}_{2,3})\,,
\end{align*}
where
\begin{multline*}
\mathbb{I}_{i,j}
	=\partial_x(h_{m_i}(|\eta_i|^{\frac{1}{2}}x))\partial_x(h_{m_j}(|\eta_j|^{\frac{1}{2}}x))\\
	=\sqrt{|\eta_i||\eta_j|}\left(\sqrt{\frac{m_i}{2}}h_{m_i-1}(|\eta_i|^{\frac{1}{2}}x)-\sqrt{\frac{m_i+1}{2}}h_{m_i+1}(|\eta_i|^{\frac{1}{2}}x)\right)\\
	\left(\sqrt{\frac{m_j}{2}}h_{m_j-1}(|\eta_j|^{\frac{1}{2}}x)-\sqrt{\frac{m_j+1}{2}}h_{m_j+1}(|\eta_j|^{\frac{1}{2}}x)\right)\,.
\end{multline*}
In particular, expanding $\mathbb{I}_{i,j}$, this term writes as a combination of terms
\begin{equation}\label{eq:shift3}
C_{A_1,A_2,A_3}(\delta) F_{I_1,m_1}^{\delta_1}(\eta_1)F_{I_2,m_2}^{\delta_2}(\eta_2)F_{I_3,m_3}^{\delta_3}(\eta_3)h_{m_1+\delta_1,m_2+\delta_2,m_3+\delta_3}(\eta_1,\eta_2,\eta_3)\,,
\end{equation}
where $|C_{A_1,A_2,A_3}(\delta)|\lesssim\max\{A_1,A_2,A_3\}$ for $\delta=(\delta_1,\delta_2,\delta_3)$ in some finite set $D_3\subset D_1^3$.

By comparing the signs of $\eta_1,\eta_2$ and $\eta_3$, one can see that $ |\eta_1|^2+|\eta_2|^2+|\eta_3|^2-|\eta_1+\eta_2+\eta_3|^2$ is a linear combination of terms $|\eta_i\eta_j|$ and $\eta_i\eta_j$ for $i\neq j$. 
We now use the identity
\[
xh_m(x)=\sqrt{\frac{m}{2}}h_{m-1}+\sqrt{\frac{m+1}{2}}h_{m+1}
\]
to remove the weight $x^2$ and write the term 
\[
x^2(|\eta_1|^2+|\eta_2|^2+|\eta_3|^2-|\eta_1+\eta_2+\eta_3|^2)h_{m_1,m_2,m_3}(\eta_1,\eta_2,\eta_3)
\]
 as a linear combination of terms as~\eqref{eq:shift3} above.
Finally, we 
 conclude that
\begin{multline*}
\mathcal{F}((\operatorname{Id}-\Delta_G)(u^{(1)}_{I_1,m_1}u^{(2)}_{I_2,m_2} u^{(3)}_{I_3,m_3}))(x,\eta)\\
	= \sum_{\delta\in D_3}C_{A_1,A_2,A_3}(\delta)\int (f^{(1)}_{I_1,m_1})^{\delta_1}(\eta_1)(f^{(2)}_{I_2,m_2})^{\delta_2}(\eta_2)(f^{(3)}_{I_3,m_3})^{\delta_3}(\eta-\eta_1-\eta_2)\\
	h_{m_1+\delta_1}(|\eta_1|^{\frac{1}{2}}x)h_{m_2+\delta_2}(|\eta_2|^{\frac{1}{2}}x)h_{m_3+\delta_3}(|\eta-\eta_1-\eta_2|^{\frac{1}{2}}x)\mathrm{d}\eta_1\mathrm{d}\eta_2\,
\end{multline*}
for some $|C_{A_1,A_2,A_3}(\delta)|\lesssim\max\{A_1,A_2,A_3\}$, leading to the lemma.
\end{proof}

\subsection{A bound on the Sobolev norm of a high-low product}

We are now able to estimate the $H^{\ell}$ norm of a product $u_Av_B$ by its $L^2$ norm as follows.

\begin{corollary}\label{cor:derivativeSplit} Let ${\ell}\in[ 0,2]$. Then there exists $C>0$ such that for all $u,v\in H^{\ell}_G$, for all $A,B\in2^{\mathbb{N}}$, there holds
\[
    \|u_{A}v_{B}\|_{H^{\ell}_G} \leq C \max\{A,B\}^{\frac{{\ell}}{2}}\sum_{(\delta_1,\delta_2)\in D_2}\|(u_{A})^{\delta_1}(v_{B})^{\delta_2}\|_{L^2_G}\,,
\]
where $D_2$ is defined by~\eqref{eq.D2}.
	
For general $\ell\geq 0$, the same result holds up to taking bigger finite sets $D_2$ 
 (depending on $\ell$) and applying the shift several times for each function.
\end{corollary} 

\begin{proof} 
 We proceed by interpolation.


    In the case ${\ell}=0$, there is nothing to do. In the case ${\ell}=2$, we use Lemma~\ref{lem:shift} to write that for any $(m+1)I\sim A$ and $(n+1)J\sim B$,
    \begin{equation*}
        (\operatorname{Id}-\Delta_G )(u_{I,m}v_{J,n})\\
	    = \sum_{(\delta_1,\delta_2)\in D_2}C_{A,B}(\delta_1,\delta_2)u^{\delta_1}_{I,m} v_{J,n}^{\delta_2}\,,
    \end{equation*}
    and thus by summation and taking the $L^2_G$ norm, we obtain
    \begin{equation*}
        \left\|(\operatorname{Id}-\Delta_G)(u_{A} v_{B})\right\|_{L^2_G}^2
        \lesssim \bigg\|\sum_{(\delta_1,\delta_2)\in D_2} C_{A,B}(\delta_1,\delta_2) \sum_{\substack{(I,m)\in 2^{\mathbb{Z}}\times\mathbb{N}\\(m+1)I\sim A}}\sum_{\substack{(J,n)\in 2^{\mathbb{Z}}\times\mathbb{N}\\  (n+1)J\sim B}}u ^{\delta_1}_{I,m}v_{J,n}^{\delta_2} \bigg\|_{L^2_G}^2\,.
    \end{equation*}
    We use the triangle inequality on the sum over $(\delta_1,\delta_2)$ and the fact that $|C_{A,B}(\delta_1,\delta_2)|\lesssim \max(A,B)$ to deduce that
    \begin{align*}
        \|(\operatorname{Id}-\Delta_G)(u_{A}v_{B})\|_{L^2_G}^2
	    &\lesssim A^{2}\sum_{(\delta_1,\delta_2)\in D_2}\bigg\|\sum_{\substack{(I,m)\in 2^{\mathbb{Z}}\times\mathbb{N} \\ (m+1)I\sim A}}\sum_{\substack{(J,n)\in 2^{\mathbb{Z}}\times\mathbb{N}\\ (n+1)J\sim B}}  u ^{\delta_1}_{I,m} v_{J,n}^{\delta_2} \bigg\|_{L^2_G}^2\\
	    &\lesssim A^2\sum_{(\delta_1,\delta_2)\in D_2} \|(u_{A})^{\delta_1} (v_{B})^{\delta_2} \|_{L^2_G}^2\,.
    \end{align*}

To get the result when ${\ell}\in(0,2)$, it only remains to interpolate thanks to the inequality
\[\|u_{A}v_{B}\|_{H^{\ell}_G} \leqslant \|u_{A}v_{B}\|_{L^2_G}^{1-\frac{{\ell}}{2}}\|u_{A}v_{B}\|_{H^2_G}^{\frac{{\ell}}{2}}\]
when $u_A, v_B\in H^2_G$, and conclude by density.


Similarly, for even integer $\ell$, we apply Lemma~\ref{lem:shift} successively $\frac{\ell}{2}$ times to get the estimate, and conclude by interpolation for exponents between $\ell$ and $\ell+2$.
\end{proof}



Using Lemma~\ref{lem:shift_3terms} instead of~\ref{lem:shift}, we get the following adaptation of Corollary~\ref{cor:derivativeSplit} for the product of three terms.

\begin{corollary}\label{cor:derivativeSplit_3terms} Let ${\ell}\in[0,2]$ and $u,v,w\in H^{\ell}_G$. Then for all $A,B,C\in2^{\mathbb{N}}$, there holds
\[
    \|u_Av_Bw_C\|_{H^{\ell}_G} \lesssim \max\{A,B,C\}^{\frac{{\ell}}{2}}\sum_{(\delta_1,\delta_2,\delta_3)\in D_3}\|(u_A)^{\delta_1}(v_B)^{\delta_2}(w_C)^{\delta_3}\|_{L^2_G}\,,
\]
where $D_3$ is defined in~\eqref{def:D3}. 
%
For general $\ell\geq 0$, the same result holds up to taking bigger finite sets $D_3$ 
(depending on $\ell$) and applying the shift several times for each function.
\end{corollary} 

\begin{proof}
By interpolation over $\ell$, it is sufficient to establish the following inequality for even integer $\ell$: for all $B,C\leq A$,
\[
\left\|(\operatorname{Id}-\Delta_G)^{\ell/2}(u_{A}v_{B}w_C)\right\|_{L^2_G}^2
	\lesssim \sum_{(\delta_1,\delta_2,\delta_3) \in D_3}  A^{\ell}\|(u_{A})^{\delta_1}(v_{B})^{\delta_2}(w_{C})^{\delta_3}\|_{L^2_G}^2\,,
\]
then apply this result to $u, P_{\leq A}v$ and $P_{\leq A}w$.

 Since in the case $\ell=0$ there is nothing to do, let us assume $\ell=2$. Using Lemma~\ref{lem:shift_3terms}, we can write
\begin{equation*}
    (\operatorname{Id}-\Delta_G )(u_{I_1,m_1}v_{I_2,m_2}w_{I_3,m_3})
	= \sum_{(\delta_1,\delta_2,\delta_3)\in D_3}C_{A,B,C}(\delta_1,\delta_2,\delta_3)u^{\delta_1}_{I_1,m_1}v_{I_2,m_2}^{\delta_2}w_{I_3,m_3}^{\delta_3}\,,
\end{equation*}
and therefore 
\begin{equation*}
    \left\|(\operatorname{Id}-\Delta_G)(u_{A} v_{B}w_C)\right\|_{L^2_G}^2
    \lesssim \Bigg\|\sum_{(\delta_1,\delta_2,\delta_3)\in D_3} C_{A,B,C}(\delta_1,\delta_2,\delta_3) \sum_{m_1,m_2,m_3\in\mathbb{N}}\sum_{\substack{(m_1+1)I_1\sim A\\(m_2+1)I_2\sim B\\(m_3+1)I_3\sim C}}u ^{\delta_1}_{I_1,m_1}v_{I_2,m_2}^{\delta_2}w_{I_3,m_3}^{\delta_3}\Bigg\|_{L^2_G}^2\,.
\end{equation*}
We use the triangle inequality on the sum over $(\delta_1,\delta_2,\delta_3)$ and the fact that $|C_{A,B,C}(\delta_1,\delta_2,\delta_3)|\lesssim \max(A,B,C) \leqslant A$ to deduce that
\begin{align*}
    \|(\operatorname{Id}-\Delta_G)(u_{A}v_{B}w_C)\|_{L^2_G}^2
    &\lesssim A^{2}\sum_{(\delta_1,\delta_2,\delta_3)\in D_3} \|(u_{A})^{\delta_1} (v_{B})^{\delta_2}(w_{C})^{\delta_3}\|_{L^2_G}^2\,.
\end{align*}
In the case of a general even integer $\ell$, successive applications of Lemma~\ref{lem:shift_3terms} lead to a similar result.

\end{proof}
\subsection{A refined Sobolev estimate of a product}

Finally, we estimate the $H^{\ell}$ norm of a product $u_Av$ by the $L^2$ norm of products $u_Av_B$ and the $H^{\ell}$ norm of $v$ as follows.

In rough terms, we should have the following. Let ${\ell}\in [0,2]$, $u,v \in L^2_G$ and $A\in2^{\mathbb{N}}$. Then for all $\varepsilon>0$, we have
\[
    \|uv\|_{H^{\ell}_G}^2 \lesssim \sum_{(\delta_1,\delta_2)\in D_2}\sum_{A,B:B\leqslant A} A^{\ell+\varepsilon}\|(u_{A})^{\delta_1}(v_{B})^{\delta_2}\|_{L^2_G}^2+\sum_{\delta\in D_1}\sum_{A}A^{\varepsilon}\|(u_{A})^{\delta_1} \|_{L^{\infty}_G}^2\|v\|_{H^{\ell}_G}^2\,,
\]
where we recall that $D_1$ and $D_2$ are defined by~\eqref{eq.D1} and~\eqref{eq.D2}.

However, in order to get nice mapping properties in $L^p$ spaces for $p\neq 2$ which are necessary during the course of this proof, see Appendix~\ref{appendix:C}, we introduce a cutoff function $\chi\in\mathcal{C}_c^{\infty}[0,1)$ such that $\chi\equiv 1$ on $[0,\frac12]$. Then, for $A\in2^{\mathbb{N}}$, we define the projection $P_{\leq A}$ as the Fourier multiplier 
\begin{equation*}
P_{\leq A}=\displaystyle\chi\left(\frac{\operatorname{Id}-\Delta_G}{A}\right)\,,
\end{equation*}
which is a smooth counterpart for the projection on the Sobolev modes $1+(2m+1)|\eta|\leq A$, acting on the decomposition~\eqref{eq:decompo_Grushin} as
\[
\mathcal{F}_{y\to\eta}(P_{\leq A}u)(x,\eta)
	=  \sum_{(I,m)\in 2^{\mathbb{Z}}\times\mathbb{N}}\chi\left(\frac{1+(2m+1)|\eta|}{A}\right)\widehat{u_{I,m}}(x,\eta)\,.
\]
Note that since on the support of $\widehat{u_{I,m}}$, we have $|\eta|\in[I,2I]$, the sum can be restricted to the indices $1+(2m+1)I\leq A$. The projection $P_{\leq A}$ commutes with the block decomposition~\eqref{def:dyadic} and the Grushin operator: 
for $B\in 2^{\mathbb{N}}$, we have
\[
P_{\leq A}(v_B)=(P_{\leq A}v)_B
\]
and for $k\in\mathbb{R}$, we have
\[
P_{\leq A}((-\Delta_G)^{ k/2}v)=(-\Delta_G)^{ k/2}(P_{\leq A}v).
\]
Moreover, for all $u\in L^2_G$, we have $\|P_{\leq A}u\|_{L^2_G}\leq \|u\|_{L^2_G}$. We also denote $P_{>A}=\operatorname{Id}-P_{\leq A}$.




\begin{lemma}[Sobolev norm for the product of two terms]\label{lem.derivativeSplit-2}
Let ${\ell}\in[0,2]$, $u,v \in L^2_G$ and $A\in2^{\mathbb{N}}$. Then for all $\varepsilon>0$, we have
\[
    \|u_Av\|_{H^{\ell}_G}^2 \lesssim \sum_{(\delta_1,\delta_2)\in D_2}\sum_{B:B\leqslant A} A^{\ell}B^{\varepsilon}\|(u_{A})^{\delta_1}(P_{\leq A}v_{B})^{\delta_2}\|_{L^2_G}^2+\sum_{\delta\in D_1}\|(u_{A})^{\delta_1} \|_{L^{\infty}_G}^2\|v\|_{H^{\ell}_G}^2\,,
\]
where we recall that $D_1$ and $D_2$ are defined by~\eqref{eq.D1} and~\eqref{eq.D2}.


 As a consequence, we have
\[
    \|uv\|_{H^{\ell}_G}^2 \lesssim \sum_{(\delta_1,\delta_2)\in D_2}\sum_{A,B:B\leqslant A} A^{\ell+\varepsilon}\|(u_{A})^{\delta_1}(P_{\leq A}v_{B})^{\delta_2}\|_{L^2_G}^2+\sum_{\delta\in D_1}\sum_{A}A^{\varepsilon}\|(u_{A})^{\delta_1} \|_{L^{\infty}_G}^2\|v\|_{H^{\ell}_G}^2\,.
\]

For general $\ell\geq 0$, the same result holds up to taking bigger finite sets $D_1,D_2$ 
 (depending on $\ell$) and applying the shift several times for each function.
 
\end{lemma}



\begin{proof}
The consequence is a simple application of the Cauchy-Schwarz' inequality.
 
We decompose $v=P_{\leq A}v+P_{> A}v$ and $v=\displaystyle\sum_{B \in 2^{\mathbb{N}}}v_B$. Since $(P_{\leq A}v_B)=v_B$ for $4B\leq A$ and $(P_{\leq A}v_B)=0$ for $B>A$, we have 
\begin{equation*}
\|u_Av\|_{H^{\ell}_G}^2
	\lesssim \bigg\|\sum_{B\leqslant A} u_{A}(P_{\leq A}v_{B})\bigg\|_{H^{\ell}_G}^2
	+\bigg\|\sum_{4B>A}u_{A}(P_{>A}v_{B})\bigg\|_{H^{\ell}_G}^2\,.
\end{equation*}
We treat the two parts separately, and in each case we proceed by interpolation between successive even integers $\ell$.

\noindent{\it Step 1: upper bound for $\displaystyle\big\|\sum_{B\leqslant A} u_{A}(P_{\leq A}v_{B})\big\|_{H^{\ell}_G}$.} Concerning the first term, thanks to the Cauchy-Schwarz inequality, we have
    \begin{equation*}
        \bigg\|\sum_{B\leqslant A}(\operatorname{Id}-\Delta_G)^{{\ell}/2} (u_{A}  (P_{\leq A}v_{B}))\bigg\|_{L^2_G}^2
         \lesssim \left(\sum_{B\leqslant A}B^{-\varepsilon}\right)
        \left(\sum_{B\leqslant A}B^{\varepsilon}\left\|(\operatorname{Id}-\Delta_G)^{{\ell}/2}(u_{A}(P_{\leq A}v_{B}))\right\|_{L^2_G}^2\right)\,,
    \end{equation*}
    and as $\sum_{\substack{B \in 2^{\mathbb{N}}}}B^{-\varepsilon}\lesssim 1$, we infer that
    \begin{align*}
        \bigg\|\sum_{B\leqslant A}(\operatorname{Id}-\Delta_G)^{{\ell}/2} (u_{A}  (P_{\leq A}v_{B}))\bigg\|_{L^2_G}^2
  	    & \lesssim \sum_{B \leqslant A}B^{\varepsilon}\|(\operatorname{Id}-\Delta_G)^{{\ell}/2}(u_{A}  (P_{\leq A}v_{B}))\|_{L^2_G}^2\,.
    \end{align*}   
But using Corollary~\ref{cor:derivativeSplit} to $(P_{\leq A}v)$ instead of $v$, we have
  \begin{equation*}
        \left\|(\operatorname{Id}-\Delta_G)^{{\ell}/2}(u_{A}  v_{B})\right\|_{L^2_G}^2
        \lesssim A^{\ell}\sum_{(\delta_1,\delta_2)\in D_2} \|(u_{A})^{\delta_1} (P_{\leq A}v_{B})^{\delta_2} \|_{L^2_G}^2\,.
    \end{equation*}

\noindent{\it Step 2: upper bound for $\displaystyle\|\sum_{4B> A} u_{A}(P_{>A}v_{B})\|_{H^{\ell}_G}$.} For the second term, we establish the following estimate by interpolation:
    \begin{equation}\label{eq:term_stein}
        \bigg\|\sum_{4B>A}(\operatorname{Id}-\Delta_G)^{{\ell}/2}(u_{A}(P_{>A}v_{B}))\bigg\|_{L^2_G}^2
	    \lesssim \sum_{\delta\in D_1}\|(u_{A})^{\delta_1} \|_{L^{\infty}_G}^2\|v\|_{H^{\ell}_G}^2\,.
    \end{equation}
    In the case ${\ell}=0$, we only need to note that from the orthogonality of the different modes $v_B$, we have
    \begin{equation*}
        \left\|u_A\sum_{4B>A}(P_{>A}v_{B})\right\|_{L^2_G}^2
	    \lesssim  \|u_A\|_{L^{\infty}_G}^2\left\|\sum_{4B>A}(P_{>A}v_{B})\right\|_{L^2_G}^2
	    \lesssim  \|u _A\|_{L^{\infty}_G}^2\|v\|_{L^2_G}^2\,.
    \end{equation*}
    In the case ${\ell}=2$, we use Lemma~\ref{lem:shift} to get that 
for any $m, n \in \mathbb{N}$ and $I,J\in 2^{\mathbb{Z}}$ such that $(m+1)I\sim A$ and $(n+1)J\sim B$, 
    \begin{equation*}
        (\operatorname{Id}-\Delta_G)(u _{I,m}(P_{> A}v)_{J,n})\\
	    = \sum_{(\delta_1,\delta_2)\in D_2}C_{A,B}(\delta_1,\delta_2)u _{I,m}^{\delta_1} (P_{> A}v)_{J,n}^{\delta_2}\,,
    \end{equation*}
    and therefore obtain
    \begin{equation*}
        \left\|\sum_{4B> A}(\operatorname{Id}-\Delta_G)u_{A}  (P_{>A}v_{B})\right\|_{L^2_G}^2
	    \lesssim \Bigg\|\sum_{(\delta_1,\delta_2)\in D_2} \sum_{4B> A}  C_{A,B}(\delta_1,\delta_2) \sum_{m, n \in \mathbb{N}}\sum_{\substack{(m+1)I\sim A\\ (n+1)J\sim B}} u_{I,m}^{\delta_1} (P_{> A}v)_{J,n}^{\delta_2}\Bigg\|_{L^2_G}^2\,.
    \end{equation*}
    We  use the triangle inequality on the sum over $(\delta_1,\delta_2)$ to get that for fixed~$A$,
    \begin{align*}
        \left\|\sum_{4B> A}(\operatorname{Id}-\Delta_G)u_{A}  (P_{>A}v_{B})\right\|_{L^2_G}^2
	    &\lesssim \sum_{(\delta_1,\delta_2)\in D_2}\Bigg\|\sum_{4B>A}C_{A,B}(\delta_1,\delta_2) \sum_{m, n \in \mathbb{N}}\sum_{\substack{(m+1)I\sim A \\ (n+1)J\sim B}} u _{I,m}^{\delta_1} (P_{> A}v)_{J,n}^{\delta_2} \Bigg\|_{L^2_G}^2\\
	    &=\sum_{(\delta_1,\delta_2)\in D_2}\left\|\sum_{4B>A}C_{A,B}(\delta_1,\delta_2)u _{A}^{\delta_1} (P_{>A}v_{B})^{\delta_2} \right\|_{L^2_G}^2\\
	    &\lesssim \sum_{(\delta_1,\delta_2)\in D_2}\|u _{A}^{\delta_1}\|_{L^{\infty}_G}^2\left\|\sum_{4B> A}C_{A,B}(\delta_1,\delta_2) (P_{>A}v_{B})^{\delta_2} \right\|_{L^2_G}^2.
    \end{align*}
    Now, by orthogonality and the fact that $|C_{A,B}(\delta_1,\delta_2)|\lesssim B$ for $4B>A$, we have
    \begin{multline*}
        \left\|\sum_{4B>A}C_{A,B}(\delta_1,\delta_2) (P_{>A}v_{B})^{\delta_2} \right\|_{L^2_G}^2
        =\sum_{4B>A}|C_{A,B}(\delta_1,\delta_2)|^2\| (P_{>A}v_{B})^{\delta_2}\|_{L^2_G}^2\\
        \lesssim \sum_{4B>A} B^2\|(P_{>A}v_{B})\|_{L^2_G}^2
         \lesssim \|v\|^2_{H^2_G}\,,
    \end{multline*}
    which concludes the proof in the case ${\ell}=2$. 
    
Estimate~\eqref{eq:term_stein} is now a consequence of an interpolation result, stated and proven in Lemma~\ref{lem.interpol}, based on the Stein's interpolation theorem. Instead of using this lemma, one could invoke the interpolation result between Sobolev spaces $[H^0,H^{\ell}]_{\theta}= H^{\theta\ell}$ for $\theta\in[0,1]$, that holds for Sobolev spaces on $\mathbb{R}^d$ (see for instance~\cite{Bergh1976}) and could be extended to the setting of the Grushin operator. 
For general $\ell$, successive uses of Lemma~\ref{lem:shift} lead to the result.
\end{proof}
 
In what follows, we will actually use the trilinear version of Lemma~\ref{lem.derivativeSplit-2}.
 
\begin{corollary}[Sobolev norm for the product of three terms]\label{coro:derivativeSplit-2.3termes} Let $\ell \in [0,2]$ and $u, v, w \in L^2_G$. Then for all $\varepsilon >0$ there holds 
\begin{equation*}
    \|uvw\|_{H^{\ell}_G}^2 \lesssim \sum_{\substack{(\delta_1,\delta_2,\delta_3) \in D_3\\A,B,C: B,C\leqslant A}} A^{\ell+\varepsilon} \|(u_{A})^{\delta_1}(P_{\leq A}v_B)^{\delta_2}(P_{\leq A}w_C)^{\delta_3}\|_{L^2_G}^2 + \sum_{\delta_1\in D_1}\sum_{A}A^{\varepsilon} \|(u_{A})^{\delta_1}\|_{L^{\infty}_G}^2\|v\|^2_{H^{\ell}_G}\|w\|_{H^{\ell}_G}^2\,,
\end{equation*}
where $D_1$ is defined in~\eqref{eq.D1} and $D_3$ is defined in~\eqref{def:D3}.
For $\ell>2$, a similar result holds up to taking bigger finite sets $D_3$ 
 and taking successive shifted functions.
\end{corollary} 
 
\begin{proof} 
We mimick the proof of Lemma~\ref{lem.derivativeSplit-2}. It is enough to establish that for every $A$, we have
\[
    \|u_{A}vw\|_{H^{\ell}_G}^2 \lesssim \sum_{\substack{(\delta_1,\delta_2,\delta_3) \in D_3\\B,C: B,C\leqslant A}} A^{\ell} (BC)^{\varepsilon} \|(u_{A})^{\delta_1}(P_{\leq A}v_B)^{\delta_2}(P_{\leq A}w_C)^{\delta_3}\|_{L^2_G}^2 + \sum_{\delta_1\in D_1} \|(u_{A})^{\delta_1}\|_{L^{\infty}_G}^2\|v\|^2_{H^{\ell}_G}\|w\|_{H^{\ell}_G}^2\,.
\]
We start by writing that 
\[
    \|u_{A}vw\|^2_{H^{\ell}_G} \lesssim \bigg\| \sum_{B,C\leqslant A} u_{A}(P_{\leq A}v_B)(P_{\leq A}w_C)\bigg\|_{H^{\ell}_G}^2 + \bigg\|\sum_{B,C}u_{A}\left(v_Bw_C-(P_{\leq A}v_B)(P_{\leq A}w_C)\right)\bigg\|^2_{H^{\ell}_G}\,. 
\]

\noindent\textit{Step 1: upper bound for} $\displaystyle\|\sum_{B,C\leqslant A} u_{A}(P_{\leq A}v_B)(P_{\leq A}w_C)\|_{H^{\ell}_G}$.
By the Cauchy-Schwarz inequality, we have
\begin{multline*}
    \Bigg\|\sum_{B,C\leqslant A}(\operatorname{Id}-\Delta_G)^{\ell/2} (u_{A}  (P_{\leq A}v_B)(P_{\leq A}w_C))\Bigg\|_{L^2_G}^2\\
    \lesssim \Bigg(\sum_{B,C\leqslant A}(BC)^{-\varepsilon}\Bigg)
    \Bigg(\sum_{B,C\leqslant A}(BC)^{\varepsilon}\left\|(\operatorname{Id}-\Delta_G)^{\ell/2}(u_{A}(P_{\leq A}v_B)(P_{\leq A}w_C))\right\|_{L^2_G}^2\Bigg) \\
     \lesssim \sum_{{B,C\leqslant A}}(BC)^{\varepsilon}\left\|(\operatorname{Id}-\Delta_G)^{\ell/2}(u_{A}(P_{\leq A}v_B)(P_{\leq A}w_C))\right\|_{L^2_G}^2\,.
\end{multline*}
We then conclude thanks to Corollary~\ref{cor:derivativeSplit_3terms} that
\[
    \Bigg\|\sum_{B,C\leqslant A}(\operatorname{Id}-\Delta_G)^{\ell/2} (u_{A}  (P_{\leq A}v_B)(P_{\leq A}w_C))\Bigg\|_{L^2_G}^2\\
    	\lesssim\sum_{\substack{(\delta_1,\delta_2,\delta_3) \in D_3\\B,C: B,C\leqslant A}} A^{\ell} (BC)^{\varepsilon} \|(u_{A})^{\delta_1}(P_{\leq A}v_B)^{\delta_2}(P_{\leq A}w_C)^{\delta_3}\|_{L^2_G}^2\,.
\]

\noindent\textit{Step 2: upper bound for} $\displaystyle\bigg\|\sum_{B,C}u_A\left(v_Bw_C-(P_{\leq A}v_B)(P_{\leq A}w_C)\right)\bigg\|_{H^{\ell}_G}$. We split 
\[
v_Bw_C-(P_{\leq A}v_B)(P_{\leq A}w_C)
	=(P_{>A}v_{B})w_C+(P_{\leq A}v_B)(P_{>A}w_{C})\,.
\]
For the first term, using the algebra property of $H^{\ell}_G$ for ${\ell}>\frac{3}{2}$ we write
\begin{align*}
    \bigg\|\sum_{B,C}u_{A}(P_{>A}v_{B})w_C\bigg\|_{H^{\ell}_G}^2
    &\lesssim     \bigg\|\sum_{B}u_{A}(P_{>A}v_{B})\bigg\|_{H^{\ell}_G}^2    \bigg\|\sum_{\substack{C:C\leq A}}w_C\bigg\|_{H^{\ell}_G}^2\\
	&\lesssim     \bigg\|\sum_{\substack{B:4B>A}}(\operatorname{Id}-\Delta_G)^{\ell/2}(u_{A}(P_{>A}v_{B}))\bigg\|_{L^2_G}^2 \|w\|_{H^{\ell}_G}^2\,,
\end{align*}
hence inequality~\eqref{eq:term_stein} in the Step~2 of Lemma~\ref{lem.derivativeSplit-2} applies and gives the expected bound. A similar treatment can be applied to the second term.
\end{proof}

\section{Deterministic bilinear estimates}\label{sec:deterministicBilinear}

In this section, we first establish deterministic bilinear and trilinear estimates for the product of two Hermite functions with rescaling, of the form $h_m(\alpha_1\cdot)h_n(\alpha_2\cdot)$ and $h_{m_1}(\alpha_1\cdot)h_{m_2}(\alpha_2\cdot)h_{m_3}(\alpha_3^\cdot)$. We then deduce deterministic bilinear estimates for the products $u_{I,m}v_{J,n}$.

\subsection{Bilinear estimates for rescaled Hermite functions}

We start by establishing  bilinear estimates for the rescaled Hermite functions, which will be at the core of our future bilinear and trilinear smoothing estimates. Let us recall that $\lambda_m=\sqrt{2m+1}$.

\begin{lemma}\label{lem.bilinearRescalingHermite} Let $\alpha>0$, $m,n\in \mathbb{N}$, and assume that 
\(\lambda_n \leqslant \frac{\alpha}{4}\lambda_m.\) 
Then there holds
\[\|h_mh_n(\alpha \cdot)\|_{L^2}^2 \lesssim \frac{1}{\alpha\lambda_m}\,.\]
\end{lemma}

\begin{proof}
We cut the Hermite functions between the space regions delimited by the pointwise estimates of  Corollary~\ref{coro.pointwiseBoundsHermite}. 

(1) In the region $R_1:=\left\{ |x| \leqslant \frac{1}{2}\lambda_m\right\}$, we use the bound $|h_m(x)| \lesssim \lambda_m^{-\frac{1}{2}}$. Therefore there holds
    \begin{equation*}\|h_mh_n(\alpha \cdot)\|_{L^2(R_1)}^2
    	\lesssim \lambda_m^{-1} \|h_n(\alpha \cdot)\|_{L^2(R_1)}^2
    	\lesssim \lambda_m^{-1}\alpha^{-1}\,.
    \end{equation*}
    
(2) In the region $R_2:=\left\{\frac{1}{2}\lambda_m \leqslant |x| \leqslant 2\lambda_m\right\}$, we know by assumption that
    \(
2\lambda_n\leqslant \frac{\alpha}{2}\lambda_m\leqslant \alpha|x|\,.
\)
We therefore use the bounds $|h_m(x)| \lesssim \left(\lambda_m^{\frac{2}{3}}+ |x^2-\lambda_m^2|\right)^{-\frac{1}{4}}$ and  $|h_n(\alpha x)| \lesssim e^{-c(\alpha x)^2}$ with $c=\frac{1}{8}$, and get
    \[\|h_mh_n(\alpha \cdot)\|_{L^2(R_2)}^2  \lesssim \int_{\frac{1}{2}\lambda_m}^{2\lambda_m}\frac{e^{-2c\alpha ^2x^2}}{\sqrt{\lambda_m^{\frac{2}{3}}+ |x^2-\lambda_m^2|}}\,\mathrm{d}x\,.\] 
    We fix the limits of the integral by setting $y= \frac{x}{\lambda_m}$, which leads to the bound:
    \begin{align*}
        \int_{\frac{1}{2}\lambda_m}^{2\lambda_m}\frac{e^{-2c\alpha ^2x^2}}{\sqrt{\lambda_m^{\frac{2}{3}}+ |x^2-\lambda_m^2|}}\,\mathrm{d}x 
        & = \int_{\frac{1}{2}}^2 \frac{e^{-2c(\alpha \lambda_m)^2y^2}}{\sqrt{\lambda_m^{-\frac{4}{3}}+|y^2-1|}}\,\mathrm{d}y \\
        & \lesssim \int_{\frac{1}{2}}^2 \frac{e^{-2c(\alpha \lambda_m)^2y^2}}{\sqrt{|y^2-1|}}\,\mathrm{d}y\\
        & \lesssim e^{-\frac{c}{2}(\alpha \lambda_m)^2} \int_{\frac{1}{2}}^2\frac{\mathrm{d}y}{\sqrt{|y^2-1|}}\,.
    \end{align*}
    Now we remark that $e^{-\frac{c}{2}(\alpha \lambda_m)^2}\lesssim\frac{1}{\alpha \lambda_m}$, giving the required estimate.

(3) In the last region $R_3:= \left\{|x|\geqslant 2\lambda_m\right\}$ we use the exponential bounds for the two terms $|h_m(x)|\leqslant e^{-cx^2}$ and $|h_n(\alpha x)|\leqslant e^{-c(\alpha x)^2}$ with $c=\frac{1}{8}$. This leads to 
    \begin{align*}
        \|h_mh_n(\alpha \cdot)\|_{L^2(R_3)}^2
        &\lesssim \int_{|x|\geqslant 2 \lambda_m} e^{-2cx^2(1+\alpha^2)}\mathrm{d}x\,. 
    \end{align*}
Then using that for $C=2c(1+\alpha^2)$ and $X=2\lambda_m$, we have
\begin{equation*}
\int_X^{\infty}e^{-Cx^2}\mathrm{d}x
	\leq \frac{1}{2CX}\int_X^{\infty}2Cxe^{-Cx^2}\mathrm{d}x
	\leq \frac{1}{2CX}\,,
\end{equation*}
this gives the result.

All the previous bounds put together imply the lemma.
\end{proof}

In what follows, we will use the bilinear estimate from Lemma~\ref{lem.bilinearRescalingHermite} under the following form. 

\begin{corollary}[Bilinear estimates for rescaled Hermite functions]\label{cor.product.rescaled.hermite}
Let $m,n\in\mathbb{N}$ and $\alpha_1,\alpha_2 >0$. Then
\[
\left\|h_m\left(\alpha_1\cdot\right)h_n\left(\alpha_2\cdot\right)\right\|_{L^2}^2
	\lesssim \min\left\{\frac{1}{\alpha_1^2(2n+1)},\frac{1}{\alpha_2^2(2m+1)}\right\}^{\frac{1}{2}}\,.
\]
\end{corollary}

\begin{proof}
In the first scenario, assume that $\alpha_1\sqrt{2n+1}\leqslant \frac{1}{4}\alpha_2\sqrt{2m+1}$.
We start with a change of variable: 
\begin{equation*}
\left\|h_m\left(\alpha_1\cdot\right)h_n\left(\alpha_2\cdot\right)\right\|_{L^2}^2
	= \frac{1}{\alpha_1} \left\|h_mh_n\left(\frac{\alpha_2}{\alpha_1}\cdot\right)\right\|_{L^2}^2\,.
\end{equation*}
Denote $\alpha:=\frac{\alpha_2}{\alpha_1}$. Since $\alpha_1\sqrt{2n+1}\leqslant \frac{1}{4}\alpha_2\sqrt{2m+1}$, we have
\(
\lambda_n
	\leq \frac{\alpha}{4}\lambda_m.
\)
Then Lemma~\ref{lem.bilinearRescalingHermite} implies 
\[
    \left\|h_mh_n\left(\frac{\alpha_2}{\alpha_1}\cdot\right)\right\|_{L^2}^2
	\lesssim \frac{1}{\alpha\sqrt{2m+1}}
	= \sqrt{\frac{\alpha_1^2}{\alpha_2^2(2m+1)}}\,.
\]
Consequently,
\[
\|h_m(\alpha_1\cdot)h_n(\alpha_2\cdot)\|_{L^2}^2
	\lesssim\sqrt{\frac{1}{\alpha_2^2(2m+1)}}\,.\]
This is enough to conclude since by assumption, we have $\sqrt{\frac{1}{\alpha_2^2(2m+1)}}\leqslant \sqrt{\frac{1}{\alpha_1^2(2n+1)}}$.

In the second scenario, assume that $\alpha_2\sqrt{2m+1}\leqslant \frac{1}{4}\alpha_1\sqrt{2n+1}$. We exchange the roles of $m$ and $n$ and the roles of $\alpha_1$ and $\alpha_2$ to get
\[
\|h_m(\alpha_1\cdot)h_n(\alpha_2\cdot)\|_{L^2}^2
	\lesssim\sqrt{\frac{1}{\alpha_1^2(2n+1)}}\,.\]

In the last scenario, there exists $C_0>0$ such that $\frac{1}{4}\alpha_2\sqrt{2m+1}\leq \alpha_1\sqrt{2n+1}\leq 4\alpha_2\sqrt{2m+1}$. Then from Hölder's inequality, one obtains the bound
\[
\|h_m(\alpha_1\cdot)h_n(\alpha_2\cdot)\|_{L^2}^2
	\lesssim \|h_m(\alpha_1\cdot)\|_{L^4}^2\|h_n(\alpha_2\cdot)\|_{L^4}^2\,.
	\]
From the $L^4$ norm estimate $\|h_m\|_{L^4}\leqslant \frac{C}{(2m+1)^{\frac{1}{8}}}$, $m\in\mathbb{N}$ (see Lemma~\ref{lem.Hermite_Lp}), we deduce
\[
\|h_m(\alpha_1\cdot)h_n(\alpha_2\cdot)\|_{L^2}^2
	\lesssim \frac{1}{(\alpha_1^2(2m+1))^{\frac{1}{4}}(\alpha_2^2(2n+1))^{\frac{1}{4}}}\,.
\]
Since $\alpha_1\sqrt{2n+1}$ and $\alpha_2\sqrt{2m+1}$ differ from a factor at most $4$, in this situation, one can bound the right-hand side by either $\sqrt{\frac{1}{\alpha_2^2(2m+1)}}$ or $\sqrt{\frac{1}{\alpha_1^2(2n+1)}}$.
\end{proof}

\begin{corollary}[Trilinear estimates for rescaled Hermite functions]\label{coro.product.rescaled.hermite.trilinear} Let $m_1,m_2,m_3\in\mathbb{N}$ and $\alpha_{1}, \alpha_{2},\alpha_{3} >0$. Then 
\[
    \|h_{m_1}(\alpha_{1} \cdot)h_{m_2}(\alpha_{2} \cdot) h_{m_3}(\alpha_{3} \cdot)\|_{L^2}^2 \lesssim \frac{1}{\alpha_{1} \alpha_{2} \alpha_{3}} \min \left\{\frac{\alpha_{i} \alpha_{j}}{(2m_i+1)^{\frac{1}{6}}(2m_j+1)^{\frac{1}{2}}}\right\}_{i\neq j\in \{1,2,3\}}\,.
\]
As a consequence, if $A\in2^{\mathbb{N}}$, $I_1,I_2,I_3\in 2^{\mathbb{Z}}$ are such that $(m_1+1)I_1\sim A$, for $\alpha_1^2\in[I_1,2I_1]$, $\alpha_2^2\in[I_2,2I_2]$ and $\alpha_3^2\in[I_3,2I_3]$, we have
\[
   \alpha_1\alpha_2\alpha_3 \|h_{m_1}(\alpha_{1} \cdot)h_{m_2}(\alpha_{2} \cdot) h_{m_3}(\alpha_{3} \cdot)\|_{L^2}^2 \lesssim  C_{\{I_i,m_i\}}^2\,,
\]
where  one can choose
\[
C_{\{I_i,m_i\}}^2
	=\frac{\langle I_1\rangle(I_2I_3)^{\frac{1}{4}}}{A^{\frac{1}{2}}(2m_2+1)^{\frac{1}{12}}(2m_3+1)^{\frac{1}{12}}}\,.
\]
\end{corollary}
Note that the gain with power $\frac{1}{12}$ induced by the $L^{\infty}$ norm on the Hermite functions is not necessary in the sequel, but we keep this exponent for the sake of security.

\begin{proof} We start by breaking the symmetry of roles of the three terms
\[
    \|h_{m_1}(\alpha_1 \cdot)h_{m_2}(\alpha_2 \cdot) h_{m_3}(\alpha_3 \cdot)\|_{L^2}^2 \leqslant \|h_{m_1}(\alpha_1 \cdot)\|_{L^{\infty}}^2 \|h_{m_2}(\alpha_2 \cdot) h_{m_3}(\alpha_3 \cdot)\|_{L^2}^2\,.
\]
Then we use Lemma~\ref{lem.Hermite_Lp} and Corollary~\ref{cor.product.rescaled.hermite} to bound
\[
    \|h_{m_1}(\alpha_1 \cdot)h_{m_2}(\alpha_2 \cdot) h_{m_3}(\alpha_3 \cdot)\|_{L^2}^2 \lesssim \frac{1}{(2m_1+1)^{\frac{1}{6}}} \min\left\{\frac{1}{\alpha_2 (2m_3+1)^{\frac{1}{2}}},\frac{1}{\alpha_3 (2m_2+1)^{\frac{1}{2}}}\right\}\,.
\]
By symmetry between the roles of $m_1,m_2,m_3$, this implies
\[
    \|h_{m_1}(\alpha_1 \cdot)h_{m_2}(\alpha_2 \cdot) h_{m_3}(\alpha_3 \cdot)\|_{L^2}^2 \lesssim \frac{1}{\alpha_1 \alpha_2 \alpha_3} \min \left\{\frac{\alpha_i \alpha_j}{(2m_i+1)^{\frac{1}{6}}(2m_j+1)^{\frac{1}{2}}}\right\}_{i\neq j\in \{1,2,3\}}\,.
\]

In order to establish the second bound, for every $j=1,2,3$, we have $\alpha_j\leq (2I_j)^{\frac{1}{2}}$, therefore
\[
\min \left\{\frac{\alpha_i \alpha_j}{(2m_i+1)^{\frac{1}{6}}(2m_j+1)^{\frac{1}{2}}}\right\}_{i\neq j\in \{1,2,3\}}
	\lesssim \min \left\{\frac{ I_i^{\frac{1}{2}} I_j^{\frac{1}{2}}}{(2m_i+1)^{\frac{1}{6}}(2m_j+1)^{\frac{1}{2}}}\right\}_{i\neq j\in \{1,2,3\}}\,.
\]
Now, let us remark that
\[
    \min \left\{ \frac{ \alpha_i \alpha_j}{(2m_i+1)^{\frac{1}{6}}(2m_j+1)^{\frac{1}{2}}}\right\}_{i \neq j \in \{1,2,3\}} 
	\leqslant    \frac{I_2^{\frac{1}{2}} I_1^{\frac{1}{2}}}{(2m_2+1)^{\frac{1}{6}}(2m_1+1)^{\frac{1}{2}}}\,,
\]
and since the same bound holds when exchanging the roles of $(I_2,m_2)$ and $(I_3,m_3)$, we obtain by interpolation
\[
    \min \left\{ \frac{\alpha_i \alpha_j }{(2m_i+1)^{\frac{1}{6}}(2m_j+1)^{\frac{1}{2}}}\right\}_{i \neq j \in \{1,2,3\}} 
    \leqslant \frac{(I_2I_3)^{\frac{1}{4}}\langle I_1\rangle }{((2m_2+1)(2m_3+1))^{\frac{1}{12}}A^{\frac 12}}\,.\qedhere
\]
\end{proof}

\subsection{Bilinear and trilinear estimates for the product of unimodal blocks}

In this part, we provide a bilinear estimate for the ``building blocks'' $u_{I,m}v_{J,n}$, where $u_{I,m}$ and $v_{J,n}$ are frequency localized and unimodal.

\begin{proposition}[Bilinear block estimate]\label{prop.buildingBlock} Let $m, n \geqslant 0$ and $I, J\in 2^{\mathbb{Z}}$. Let $u_{I,m}$ (resp. $v_{J,n}$) in $L^2_G$ defined by~\eqref{def:BuildingBlock}
\[\mathcal{F}_{y\to\eta}( u_{I,m}) (x,\eta) = f_{I,m}(\eta)h_m(|\eta|^{\frac{1}{2}}x)
\] resp. 
\[\mathcal{F}_{y\to\eta}(  v_{J,n}) (x,\eta)=g_{J,n}(\eta)h_n(|\eta|^{\frac{1}{2}}x)\,,
\]
the function $f_{I,m}$ (resp. $g_{J,n}$) being supported on the set $\{|\eta|\in[I,2I]\}$ (resp. on the set $\{|\eta|\in[J,2J]\}$).
Then, one has
\begin{equation*}
    \|u_{I,m}v_{J,n}\|_{L^2_G}^2 \lesssim \min\left\{I,J\right\}\min\left\{\frac{ I }{2m+1},\frac{ J }{2n+1}\right\}^{\frac{1}{2}}\|u_{I,m}\|^2_{L^2_G}\|v_{J,n}\|^2_{L^2_G}\,.
\end{equation*}
As a consequence,
\[\|u_{I,m}v_{J,n}\|_{L^2_G}^2 \lesssim \min\left\{\frac{ J\langle I\rangle }{(1+(2m+1)I)^{\frac{1}{2}}},\frac{I \langle J\rangle}{(1+(2n+1)J)^{\frac{1}{2}}}\right\}\|u_{I,m}\|^2_{L^2_G}\|v_{J,n}\|^2_{L^2_G}\,.\]
\end{proposition}
\begin{proof}
From the Parseval formula, we have $\|u_{I,m}v_{J,n}\|_{L^2_G}^2=\|\hat u_{I,m} * \hat v_{J,n}\|_{L^2_{x,\eta}}^2$, where the convolution is the classical convolution product in the variable $\eta$. We expand the norm of this convolution product and obtain
\begin{align*}
    \|u_{I,m}v_{J,n}\|_{L^2_G}^2 = \int_{\mathbb{R}^3} f_{I,m}(\eta_1)f_{I,m}(\eta_2)g_{J,n}(\eta-\eta_1)g_{J,n}(\eta - \eta_2) \mathbb{I}_{m,n}(\eta_1,\eta_2,\eta-\eta_1,\eta-\eta_2)\,\mathrm{d}\eta_1\,\mathrm{d}\eta_2\,\mathrm{d}\eta\,,
\end{align*}
with 
\begin{equation*}
    \mathbb{I}_{m,n}(\eta_1,\eta_2,\eta'_1,\eta'_2) := \int_{\mathbb{R}} h_m(|\eta_1|^{\frac{1}{2}}x)h_m(|\eta_2|^{\frac{1}{2}}x)h_n(|\eta_1'|^{\frac{1}{2}}x)h_n(|\eta_2'|^{\frac{1}{2}}x)\,\mathrm{d}x\,.
\end{equation*}
By Cauchy-Schwarz' inequality, one has
\begin{equation*}
|\mathbb{I}_{m,n}(\eta_1,\eta_2,\eta'_1,\eta'_2)|
	\leqslant \|h_m(|\eta_1|^{\frac{1}{2}}\cdot)h_n(|\eta_1'|^{\frac{1}{2}}\cdot)\|_{L^2} \|h_m(|\eta_2|^{\frac{1}{2}}\cdot)h_n(|\eta_2'|^{\frac{1}{2}}\cdot)\|_{L^2}\,.
\end{equation*}
Using Corollary~\ref{cor.product.rescaled.hermite}, we deduce the estimate
\begin{equation*}
|\mathbb{I}_{m,n}(\eta_1,\eta_2,\eta'_1,\eta'_2)|
	\lesssim \min\left\{\frac{1}{|\eta_1|(2n+1)},\frac{1}{|\eta'_1|(2m+1)}\right\}^{\frac{1}{4}}
	\min\left\{\frac{1}{|\eta_2|(2n+1)},\frac{1}{|\eta'_2|(2m+1)}\right\}^{\frac{1}{4}}\,.
\end{equation*}
Going back to the blocks $u_{I,m}$ and $v_{J,n}$, this estimate implies
\begin{align*}
\|u_{I,m}v_{J,n}\|_{L^2_G}^2
    &	\lesssim \int_{\mathbb{R}} \left(\int_{\mathbb{R}} \min\left\{\frac{1}{|\eta_1|(2n+1)},\frac{1}{|\eta-\eta_1|(2m+1)}\right\}^{\frac{1}{4}}|f_{I,m}(\eta_1)g_{J,n}(\eta-\eta_1)|\,\mathrm{d}\eta_1\right)^2\,\mathrm{d}\eta\,.
\end{align*}
By extracting the minimum, we conclude
\[
\|u_{I,m}v_{J,n}\|_{L^2_G}^2
	\lesssim \min\left\{\frac{1}{\sqrt{2n+1}} \|\frac{|f_{I,m}(\cdot)|}{|\cdot|^{1/4}} * |g_{J,n}|\|^2_{L^2},
	\frac{1}{\sqrt{2m+1}} \||f_{I,m}| * \frac{|g_{J,n}(\cdot)|}{|\cdot|^{1/4}}\|^2_{L^2}\right\}\,.\]
	
Using the symmetry of the roles of $u$ and $v$, we estimate for instance $  \|\frac{|f_{I,m}(\cdot)|}{|\cdot|^{1/4}} * |g_{J,n}|\|^2_{L^2}$.
We apply Young's inequality
\[
   \ \|\frac{|f_{I,m}(\cdot)|}{|\cdot|^{1/4}} * |g_{J,n}|\|^2_{L^2} \leqslant  \|\frac{f_{I,m}(\cdot)}{|\cdot|^{1/4}} \|^2_{L^1}\|g_{J,n}\|^2_{L^2}\,.
\]
Thanks to Cauchy-Schwarz' inequality, since the interval $[I,2I]$ has length $I$, we get
\[
\|\frac{f_{I,m}(\cdot)}{|\cdot|^{1/4}} \|_{L^1}^2
	\lesssim I\|\frac{f_{I,m}(\cdot)}{|\cdot|^{1/4}} \|_{L^2}^2\,.\]
Moreover, we know that $\|g_{J,n}\|_{L^2}^2\lesssim \sqrt{J}\|v_{J,n}\|_{L^2}^2$, therefore
\[
   \ \|\frac{|f_{I,m}(\cdot)|}{|\cdot|^{1/4}} * |g_{J,n}|\|^2_{L^2} \leqslant I\sqrt{J}\|u_{I,m} \|^2_{L^2_G}\|v_{J,n}\|^2_{L^2_G}\,.
\]
But we also have from Young's inequality
\[
   \ \|\frac{|f_{I,m}(\cdot)|}{|\cdot|^{1/4}} * |g_{J,n}|\|^2_{L^2} \leqslant  \|\frac{f_{I,m}(\cdot)}{|\cdot|^{1/4}} \|^2_{L^2}\|g_{J,n}\|^2_{L^1}\,,
\]
where from Cauchy-Schwarz' inequality, $\|g_{J,n}\|^2_{L^1}\leqslant J\|g_{J,n}\|^2_{L^2}\lesssim J^{\frac{3}{2}}\|v_{J,n}\|^2_{L^2_G}$, so that actually
\[
   \ \|\frac{|f_{I,m}(\cdot)|}{|\cdot|^{1/4}} * |g_{J,n}|\|^2_{L^2} \leqslant  \min\left\{I,J\right\}\sqrt{J}\|u_{I,m} \|^2_{L^2_G}\|v_{J,n}\|^2_{L^2_G}\,.
\]
To conclude, using the symmetry of the roles of $u$ and $v$, we have proven that
\[
\|u_{I,m}v_{J,n}\|_{L^2_G}^2 \lesssim  \min\left\{I,J\right\}\min\left\{\frac{I}{2m+1},\frac{J}{2n+1}\right\}^{\frac{1}{2}}\|u_{I,m}\|_{L^2_G}^2\|v_{J,n}\|_{L^2_G}^2\,.\qedhere
\] 
\end{proof}

\section{Probabilistic bilinear and trilinear estimates}
\label{sec:probabilisticBilinear}

In this section we provide the proof of Theorem~\ref{th.bilinear_estimate}~(\textit{i}). 
The main idea of proof is to use the deterministic bilinear estimate for the products $u_{I,m}v_{J,n}$ given by Proposition~\ref{prop.buildingBlock} in order to prove a smoothing effect on quadratic and cubic expressions of a random function.

\subsection{Bilinear estimate for random interactions} 

To deal with purely random interactions of the form $(u_0^{\omega})^2$ and  $|u_0^{\omega}|^2$, we heavily rely on the frequency decoupling offered by the randomization. This is akin to the bilinear estimates obtained in~\cite{btz}. 

Let us recall some notation. We fix a sequence of independent identically distributed subgaussian random variables $(X_{I,m})_{(I,m)\in 2^{\mathbb{Z}}\times\mathbb{N}}$ and $u_0\in H^k_G$. We use decomposition~\eqref{eq:randomizationMap}
\[
    u_0^{\omega}= \sum_{(I,m)\in 2^{\mathbb{Z}}\times\mathbb{N}}X_{I,m}(\omega)u_{I,m}\,.
\]

For all $t\in\mathbb{R}$, we consider the time evolution under the linear flow~\eqref{eq:LS-G} with initial data $u_{I,m}$ denoted $z_{I,m}(t)=e^{it\Delta_G}u_{I,m}$, $(I,m)\in 2^{\mathbb{Z}}\times\mathbb{N}$. Then, we define the random counterpart  $z^{\omega}$ of $u$ 
as in~\eqref{eq:decomposition_Grushin_random}
\[
    z^{\omega}(t)=e^{it\Delta_G}u_0^{\omega} = \sum_{(I,m)\in 2^{\mathbb{Z}}\times\mathbb{N}}X_{I,m}(\omega)z_{I,m}(t)\,.
\]

\begin{proof}[Proof of~\eqref{eq:Bilin_zz} in Theorem~\ref{th.bilinear_estimate}~(\textit{i})]
We establish estimate~\eqref{eq:Bilin_zz}, that is, we bound the norms of the products $\|(z^{\omega})^2\|_{L^q_TH_G^{\ell}}$ and $\||z^{\omega}|^2\|_{L^q_TH_G^{\ell}}$ for $\ell=k+\frac{3}{2}$.
\medskip

\noindent\textit{Step 1: reduction of inequality~\eqref{eq:Bilin_zz} to deterministic estimates.} 
We start with $\|(z^{\omega})^2\|_{L^q_TH_G^{\ell}}$. Fix $t\in[0,T]$ and $(x,y)\in \mathbb{R}^2$. Applying Corollary~\ref{coro:WienerChaos}~(\textit{i}) with the norm  $L^q_TL^2_{x,y}$ and $\Psi_{(I,m),(J,n)}=(\operatorname{Id}-\Delta_G)^{\ell/2}(z_{I,m}z_{J,n})$, we obtain that outside a set of probability at most $e^{-cR^2}$ there holds
\[
    \|(z^{\omega})^2\|_{L^q_TH^{\ell}_G}^2 
    \leqslant R^4\sum_{(I,m)\in 2^{\mathbb{Z}}\times\mathbb{N}} \sum_{(J,n)\in 2^{\mathbb{Z}}\times\mathbb{N}}\|z_{I,m}z_{J,n}\|^2_{L^q_TH^{\ell}_G}\,.
\]
It remains to prove that  for every $u_0,v_0\in\mathcal{X}^k_1$, denoting $z_{I,m}(t)=e^{it\Delta_G}u_{I,m}$ and $\widetilde{z}_{J,n}(t)=e^{it\Delta_G}v_{J,n}$ for $t\in\mathbb{R}$, we have
\begin{equation}
    \label{eq:AlmostBilin_zz_time}
    \sum_{(I,m)\in 2^{\mathbb{Z}}\times\mathbb{N}} \sum_{(J,n)\in 2^{\mathbb{Z}}\times\mathbb{N}}\|z_{I,m}\tilde{z}_{J,n}\|^2_{L^q_TH^{\ell}_G} \lesssim T^{\frac{2}{q}}\|u_0\|_{\mathcal{X}^{k}_1}^2\|v_0\|_{\mathcal{X}_1^{k}}^2\,.
\end{equation}
Indeed we will conclude by taking $u_0=v_0$.

Similarly for $\||z^{\omega}|^2\|_{L^q_TH_G^{\ell}}$, applying Corollary~\ref{coro:WienerChaos}~(\textit{i}) with the norm $L^q_TL^2_{x,y}$ and $\Psi_{(I,m),(J,n)}=(\operatorname{Id}-\Delta_G)^{\ell/2}(z_{I,m}\overline{z_{J,n}})$, we obtain that outside a set of probability at most $e^{-cR^2}$ there holds
\[
    \||z^{\omega}|^2\|_{L^q_TH^{\ell}_G}^2
    \leqslant R^4\sum_{(I,m),(J,n)\in 2^{\mathbb{Z}}\times\mathbb{N}}\|z_{I,m}\overline{z_{J,n}}\|^2_{L^q_TH^{\ell}_G}
    +R^4 \left(\sum_{(I,m)\in 2^{\mathbb{Z}}\times\mathbb{N}}\||z_{I,m}|^2\|_{L^q_TH^{\ell}_G}\right)^2\,.
\]
The upper bound is handled using inequality~\eqref{eq:AlmostBilin_zz_time} and
establishing
\begin{equation}
    \label{eq:AlmostBilin_zz*_time}
    \sum_{(I,m)\in 2^{\mathbb{Z}}\times\mathbb{N}}\|z_{I,m}\tilde{z}_{I,m}\|_{L^q_TH^{\ell}_G} \lesssim T^{\frac{1}{q}}\|u_0\|_{\mathcal{X}_1^k}\|v_0\|_{\mathcal{X}_1^k}\,,
\end{equation}
which we apply to $u$ and $v=\overline{u}$.
\medskip

\noindent\textit{Step 2: Proof of~\eqref{eq:AlmostBilin_zz_time}.} 
We claim that~\eqref{eq:AlmostBilin_zz_time} is a consequence of the time independent inequality
\begin{equation}
    \label{eq:AlmostBilin_zz}
    \sum_{(I,m),(J,n)\in 2^{\mathbb{Z}}\times\mathbb{N}}\|u_{I,m}v_{J,n}\|^2_{H^{\ell}_G} \lesssim \|u_0\|_{\mathcal{X}_{1}^k}^2\|v_0\|_{\mathcal{X}_{1}^{k}}^2\,.
\end{equation}
Indeed we apply~\eqref{eq:AlmostBilin_zz} to $z_{I,m}(t)=e^{it\Delta_G}u_{I,m}$ and $\widetilde{z}_{J,n}(t)=e^{it\Delta_G}v_{J,n}$ for $t\in\mathbb{R}$ instead of $u_{I,m}$ and $v_{J,n}$, then we use $\mathcal{X}^k_1$ isometry property of the linear flow, where we recall that we have defined  the norm~\eqref{def:Xk_rho}
\[
    \|u_0\|_{\mathcal{X}^k_{1}}^2 = \sum_{(I,m)\in 2^{\mathbb{Z}}\times\mathbb{N}} (1+(2m+1)I)^{k}\left\langle I \right\rangle  \|u_{I,m}\|^2_{L^2}\,.
\]
Finally we integrate in time using the Hölder inequality, yielding to~\eqref{eq:AlmostBilin_zz_time}. 

In order to prove~\eqref{eq:AlmostBilin_zz}, fix $I,J\in 2^{\mathbb{Z}}$ and $m,n\in\mathbb{N}$. Let $A,B\in2^{\mathbb{N}}$ such that $(m+1)I\sim A$ and $(n+1)J\sim B$. We apply Corollary~\ref{cor:derivativeSplit} and get
\[
    \|u_{I,m}v_{J,n}\|^2_{H^{\ell}_G} \lesssim \max\{A,B\}^{\ell}\|u_{I,m}v_{J,n}\|^2_{L^2_G}\,.
\]
Then, from the unit block bound of Proposition~\ref{prop.buildingBlock}, it follows that 
\begin{equation}\label{eq:bilin_uv}
    \|u_{I,m}v_{J,n}\|_{H^{\ell}_G}^2 \lesssim \max\{A,B\}^{\ell}\min\left\{\frac{ J\langle I\rangle  }{A^{\frac{1}{2}}},\frac{ I\langle J\rangle  }{B^{\frac{1}{2}}}\right\}\|u_{I,m}\|^2_{L^2_G}\|v_{J,n}\|^2_{L^2_G} \,.
\end{equation}
Separating the cases $A\leqslant B$ and $A>B$, we obtain the bounds
\begin{multline*}
\sum_{\substack{(I,m),(J,n)\in 2^{\mathbb{Z}}\times\mathbb{N}}} \|u_{I,m}v_{J,n}\|^2_{H^{\ell}_G}
    \lesssim \sum_{\substack{(I,m),(J,n)\in 2^{\mathbb{Z}}\times\mathbb{N}\\A\leqslant B}} B^{\ell-\frac{1}{2}}  I\langle J\rangle \|u_{I,m}\|^2_{L^2_G}\|v_{J,n}\|^2_{L^2_G}\\
    + \sum_{\substack{(I,m),(J,n)\in 2^{\mathbb{Z}}\times\mathbb{N}\\A> B}} A^{\ell-\frac{1}{2}}  J\langle I\rangle \|u_{I,m}\|^2_{L^2_G}\|v_{J,n}\|^2_{L^2_G}\,,
\end{multline*}
implying~\eqref{eq:AlmostBilin_zz} since $\ell - \frac{1}{2} = k$. 
\medskip

\noindent\textit{Step 3: Proof of~\eqref{eq:AlmostBilin_zz*_time}.} 
We claim that~\eqref{eq:AlmostBilin_zz*_time} is a consequence of the time independent inequality
\begin{equation}
    \label{eq:AlmostBilin_zz*}
   \sum_{(I,m)\in 2^{\mathbb{Z}}\times\mathbb{N}}  \| u_{I,m}v_{I,m}\|_{H^{\ell}_G} \lesssim \|u_0\|_{\mathcal{X}^k_1}\|v_0\|_{\mathcal{X}^k_1}\,.
\end{equation}
Indeed we apply~\eqref{eq:AlmostBilin_zz*} to $z_{I,m}(t)=e^{it\Delta_G}u_{I,m}$ and $\widetilde{z}_{I,m}(t)=e^{it\Delta_G}v_{I,m}$ for $t\in\mathbb{R}$ instead of $u_{I,m}$, then we use $\mathcal{X}^{k}_{1}$ isometry property of the linear flow and integrate in time using the Hölder inequality. 
This yields~\eqref{eq:AlmostBilin_zz*_time}. 

In order to prove~\eqref{eq:AlmostBilin_zz*}, for $I\in 2^{\mathbb{Z}}$ and $m\in\mathbb{N}$, we denote $A\in2^{\mathbb{N}}$ such that $(m+1)I\sim A$. From~\eqref{eq:bilin_uv} applied to $I=J$, $m=n$ and $A=B$, 
we therefore write
\begin{align*}
    \sum_{\substack{(I,m)\in 2^{\mathbb{Z}}\times\mathbb{N}}}\|u_{I,m}v_{I,m}\|_{H^{\ell}_G}
    \lesssim \sum_{\substack{(I,m)\in 2^{\mathbb{Z}}\times\mathbb{N}}}A^{\frac{\ell}{2}-\frac{1}{4}}\langle I\rangle\|u_{I,m}\|_{L^2_G}\|v_{I,m}\|_{L^2_G}\,.
\end{align*}
An application of Cauchy-Schwarz' inequality implies~\eqref{eq:AlmostBilin_zz*} since $\ell - \frac{1}{2} = k$. 
\end{proof}

\subsection{Trilinear estimate for random interactions} 

\begin{proof}[Proof of~\eqref{eq:Bilin_zzz*} in Theorem~\ref{th.bilinear_estimate}~(\textit{i})]
We now estimate $\||z^{\omega}|^2z^{\omega}\|_{L^q_TH_G^{\ell}}$.\medskip

\noindent\textit{Step 1: reduction of inequality~\eqref{eq:Bilin_zzz*} to deterministic estimates.} Using Corollary~\ref{coro:WienerChaos}~(\textit{ii}), we know that outside a set of probability at most $e^{-cR^2}$, there holds 
\begin{multline}
\label{eq:probaDecouplingTrilinear}
    \||z^{\omega}|^2z^{\omega}\|_{L^q_TH^{\ell}_G}^2
    \leqslant R^6 \sum_{(I_1,m_1),(I_2,m_2),(I_3,m_3)\in 2^{\mathbb{Z}}\times\mathbb{N}}\|z_{I_1,m_1}z_{I_2,m_2}\overline{z_{I_3,m_3}}\|^2_{L^q_TH^{\ell}_G}\\
    +R^6\sum_{(I_2,m_2)\in 2^{\mathbb{Z}}\times\mathbb{N}} \left(\sum_{(I_1,m_1)\in 2^{\mathbb{Z}}\times\mathbb{N}}\||z_{I_1,m_1}|^2z_{I_2,m_2}\|_{L^q_TH^{\ell}_G}\right)^2\,.
\end{multline}

As in the two subsections above, using Hölder's inequality in time, inequality~\eqref{eq:Bilin_zzz*} is now a consequence of the time independent inequalities
\begin{equation}
    \label{eq:AlmostBilin_zzz*_1}
\sum_{(I_1,m_1),(I_2,m_2),(I_3,m_3)\in 2^{\mathbb{Z}}\times\mathbb{N}}\|u_{I_1,m_1}u_{I_2,m_2}\overline{u_{I_3,m_3}}\|^2_{H^{\ell}_G}
    \lesssim \|u_0\|^6_{\mathcal{X}_1^k}\,,
\end{equation}
\begin{equation}
    \label{eq:AlmostBilin_zzz*_2}
    \sum_{(I_2,m_2)\in 2^{\mathbb{Z}}\times\mathbb{N}} \left(\sum_{(I_1,m_1)\in 2^{\mathbb{Z}}\times\mathbb{N}}\||u_{I_1,m_1}|^2u_{I_2,m_2}\|_{H^{\ell}_G}\right)^2
    \lesssim \|u_0\|^6_{\mathcal{X}_1^k}\,.
\end{equation}
\medskip

\noindent\textit{Step 2: Proof of~\eqref{eq:AlmostBilin_zzz*_1}.} 
We rather prove inequality
\begin{equation}
    \label{eq:AlmostBilin_zzz*_0}
\sum_{(I_1,m_1),(I_2,m_2),(I_3,m_3)\in 2^{\mathbb{Z}}\times\mathbb{N}}    \|u^{(1)}_{I_1,m_1}u^{(2)}_{I_2,m_2}u^{(3)}_{I_3,m_3}\|_{H^{\ell}_G}^2 
    \lesssim \|u^{(1)}\|^2_{\mathcal{X}_1^k}\|u^{(2)}\|^2_{\mathcal{X}_1^k}\|u^{(2)}\|^2_{\mathcal{X}_1^k}\,,
\end{equation}
for fixed $u^{(1)},u^{(2)},u^{(3)}\in L^2_G$ decomposed as in~\eqref{eq:decompo_Grushin}, as this implies~\eqref{eq:AlmostBilin_zzz*_1} with $u^{(1)}=u^{(2)}=u_0$ and $u^{(3)}=\overline{u_0}$. Let $A_1,A_2,A_3$ be the dyadic integers such that $(m_1+1)I_1\sim A_1$, $(m_2+1)I_2\sim A_2$ and $(m_3+1)I_3\sim A_3$.

We apply Corollary~\ref{cor:derivativeSplit_3terms} and get
\[
    \|u^{(1)}_{I_1,m_1}u^{(2)}_{I_2,m_2}u^{(3)}_{I_3,m_3}\|_{H^{\ell}_G}^2 
    \lesssim \max\{A_1,A_2,A_3\}^{\ell} \|u^{(1)}_{I_1,m_1}u^{(2)}_{I_2,m_2}u^{(3^)}_{I_3,m_3}\|_{L^2_G}^2\,.
\]
Assuming up to permutation that $\max\{A_1,A_2,A_3\}=A_1$, we deduce
\[
    \|u^{(1)}_{I_1,m_1}u^{(2)}_{I_2,m_2}u^{(3)}_{I_3,m_3}\|_{H^{\ell}_G}^2 
    \lesssim A_1^{\ell}\|u^{(1)}_{I_1,m_1}u^{(2)}_{I_2,m_2}\|_{L^2_G}^2 \|u^{(3)}_{I_3,m_3}\|_{L^{\infty}_G}^2\,.
\]
Then, from the unit block bound of Proposition~\ref{prop.buildingBlock}, it follows that if $\max\{A_1,A_2,A_3\}=A_1$, then 
\begin{equation}\label{eq:AlmostBilin_zzz*_decomposed}
    \|u^{(1)}_{I_1,m_1}u^{(2)}_{I_2,m_2}u^{(3)}_{I_3,m_3}\|_{H^{\ell}_G}^2 
    \lesssim A_1^{\ell-\frac{1}{2}} \langle I_1\rangle\langle I_2\rangle \|u^{(1)}_{I_1,m_1}\|^2_{L^2_G}\|u^{(2)}_{I_2,m_2}\|^2_{L^2_G}\|u^{(3)}_{I_3,m_3}\|_{L^{\infty}_G}^2\,. 
\end{equation}
This implies that
\begin{multline*}
\sum_{\substack{(I_1,m_1),(I_2,m_2),(I_3,m_3)\in 2^{\mathbb{Z}}\times\mathbb{N}\\\max\{A_1,A_2,A_3\}=A_1}}    \|u^{(1)}_{I_1,m_1}u^{(2)}_{I_2,m_2}u^{(3)}_{I_3,m_3}\|_{H^{\ell}_G}^2 \\
    \lesssim \sum_{\substack{(I_1,m_1),(I_2,m_2),(I_3,m_3)\in 2^{\mathbb{Z}}\times\mathbb{N}\\\max\{A_1,A_2,A_3\}=A_1}}   A_1^{\ell}\frac{ \langle I_1\rangle\langle I_2\rangle }{A_1^{\frac{1}{2}}} \|u^{(1)}_{I_1,m_1}\|^2_{L^2_G}\|u^{(2)}_{I_2,m_2}\|^2_{L^2_G}\|u^{(3)}_{I_3,m_3}\|_{L^{\infty}_G}^2\,. 
\end{multline*}
By definition of $\mathcal{X}^k_\rho$, we get that for $k=\ell-\frac{1}{2}$,
\begin{multline*}
\sum_{\substack{(I_1,m_1),(I_2,m_2),(I_3,m_3)\in 2^{\mathbb{Z}}\times\mathbb{N}\\\max\{A_1,A_2,A_3\}=A_1}}    \|u^{(1)}_{I_1,m_1}u^{(2)}_{I_2,m_2}u^{(3)}_{I_3,m_3}\|_{H^{\ell}_G}^2 \\
    \lesssim \|u^{(1)}\|_{\mathcal{X}^k_1}^2\|u^{(2)}\|_{\mathcal{X}^0_1}^2\sum_{\substack{(I_3,m_3)\in \mathbb{N}^*\times \mathbb{N}}} \|u^{(3)}_{I_3,m_3}\|_{L^{\infty}_G}^2\,. 
\end{multline*}
It only remains to use estimate~\eqref{eq.largeDeviation_v0} from Lemma~\ref{cor:random_integrability_improvement} and the embedding $W^{\varepsilon,p}_G\hookrightarrow L^{\infty}_G$ for arbitrary small $\varepsilon>0$ and $p>\frac{3}{\varepsilon}$ to deduce that one has
\[
    \sum_{\substack{(I_3,m_3)\in \mathbb{N}^*\times \mathbb{N}}} \|u^{(3)}_{I_3,m_3}\|_{L^{\infty}_G}^2
	\lesssim \|u^{(3)}\|_{\mathcal{X}^{-\zeta(p)+\varepsilon}_{\zeta(p)+\frac{3}{2}-\frac{3}{p}}}^2\,.
\]
Remark that for large $p$ we have $\mathcal{X}_1^k \hookrightarrow \mathcal{X}^{-\zeta (p)+\varepsilon}_{\zeta(p)+\frac{3}{2}-\frac{3}{p}}$. This is indeed the case because $\varepsilon + \frac{1}{2}-\frac{3}{p}
\leqslant k$ for small $\varepsilon$ since $k> \frac{1}{2}$. 
We conclude that
\begin{equation*}
    \sum_{\substack{(I_1,m_1),(I_2,m_2),(I_3,m_3)\in 2^{\mathbb{Z}}\times\mathbb{N}\\\max\{A_1,A_2,A_3\}=A_1}}    \|u^{(1)}_{I_1,m_1}u^{(2)}_{I_2,m_2}u^{(3)}_{I_3,m_3}\|_{H^{\ell}_G}^2     \lesssim \|u^{(1)}\|_{\mathcal{X}^k_1}^2\|u^{(2)}\|_{\mathcal{X}^k_1}^2 \|u^{(3)}\|_{\mathcal{X}^k_1}^2\,. 
\end{equation*}
By symmetry, this inequality is also valid when $\max\{A_1,A_2,A_3\}=A_2$ or $A_3$. This implies~\eqref{eq:AlmostBilin_zzz*_0} and therefore~\eqref{eq:AlmostBilin_zzz*_1}. \medskip

\noindent\textit{Step 3: Proof of~\eqref{eq:AlmostBilin_zzz*_2}.} In order to prove~\eqref{eq:AlmostBilin_zzz*_2}, fix $I_1,I_2\in 2^{\mathbb{Z}}$ and $m_1,m_2\in\mathbb{N}$. 

If $\max\{A_1,A_2\}=A_1$, we use inequality~\eqref{eq:AlmostBilin_zzz*_decomposed} with $u^{(1)}=u^{(2)}=u_0$ and $u^{(3)}=\overline{u_0}$:
\begin{equation*}
    \|u_{I_1,m_1}u_{I_2,m_2}\overline{u}_{I_1,m_1}\|_{H^{\ell}_G}^2 
    \lesssim A_1^{\ell}\frac{ \langle I_1\rangle\langle I_2\rangle}{A_1^{\frac{1}{2}}}\|u_{I_1,m_1}\|^2_{L^2_G}\|u_{I_2,m_2}\|^2_{L^2_G}\|u_{I_1,m_1}\|_{L^{\infty}_G}^2\,. 
\end{equation*}
Otherwise, we have $\max\{A_1,A_2\}=A_2$. In this case, we still use inequality~\eqref{eq:AlmostBilin_zzz*_decomposed} with $u^{(1)}=u^{(2)}=u_0$ and $u^{(3)}=\overline{u_0}$:
\begin{equation*}
    \|u_{I_2,m_2}u_{I_1,m_1}\overline{u}_{I_1,m_1}\|_{H^{\ell}_G}^2 
    \lesssim A_2^{\ell}\frac{ \langle I_1\rangle\langle I_2\rangle}{A_2^{\frac{1}{2}}}\|u_{I_2,m_2}\|^2_{L^2_G}\|u_{I_1,m_1}\|^2_{L^2_G}\|u_{I_1,m_1}\|_{L^{\infty}_G}^2\,. 
\end{equation*}

We deduce by summation that
\begin{multline*}
\sum_{(I_2,m_2)\in 2^{\mathbb{Z}}\times\mathbb{N}} \left(\sum_{(I_1,m_1)\in 2^{\mathbb{Z}}\times\mathbb{N}}\||u_{I_1,m_1}|^2u_{I_2,m_2}\|_{H^{\ell}_G}\right)^2\\
	\lesssim \sum_{(I_2,m_2)\in 2^{\mathbb{Z}}\times\mathbb{N}} \left(\sum_{(I_1,m_1)\in 2^{\mathbb{Z}}\times\mathbb{N}}   A_1^{\frac{\ell}{2}-\frac{1}{4}} (\langle I_1\rangle\langle I_2\rangle)^{\frac 12}\|u_{I_1,m_1}\|_{L^2_G}\|u_{I_2,m_2}\|_{L^2_G} \|u_{I_1,m_1}\|_{L^{\infty}_G}\right)^2\\
	 +  \sum_{(I_2,m_2)\in 2^{\mathbb{Z}}\times\mathbb{N}} \left(\sum_{(I_1,m_1)\in 2^{\mathbb{Z}}\times\mathbb{N}}    \|u_{I_1,m_1}\|_{L^2_G} A_2^{\frac{\ell}{2}-\frac{1}{4}} (\langle I_1\rangle\langle I_2\rangle)^{\frac 12} \|u_{I_2,m_2}\|_{L^2_G}\|u_{I_1,m_1}\|_{L^{\infty}_G}\right)^2.
\end{multline*}
Now we apply Cauchy-Schwarz' inequality and get
\begin{multline*}
\sum_{(I_2,m_2)\in 2^{\mathbb{Z}}\times\mathbb{N}} \left(\sum_{(I_1,m_1)\in 2^{\mathbb{Z}}\times\mathbb{N}}\||u_{I_1,m_1}|^2u_{I_2,m_2}\|_{H^{\ell}_G}\right)^2\\
	\lesssim \sum_{(I_2,m_2)\in 2^{\mathbb{Z}}\times\mathbb{N}}
	\left(\sum_{(I_1,m_1)\in 2^{\mathbb{Z}}\times\mathbb{N}}  A_1^{\ell-\frac{1}{2}} \langle I_1\rangle \|u_{I_1,m_1}\|_{L^2_G}^2\right) \langle I_2\rangle\|u_{I_2,m_2}\|_{L^2_G}^2
	 \left(\sum_{(I_1,m_1)\in 2^{\mathbb{Z}}\times\mathbb{N}}     \|u_{I_1,m_1}\|_{L^{\infty}_G}^2\right) \\
	 +  \sum_{(I_2,m_2)\in 2^{\mathbb{Z}}\times\mathbb{N}} 
	 \left( \sum_{(I_1,m_1)\in 2^{\mathbb{Z}}\times\mathbb{N}}\langle I_1\rangle  \|u_{I_1,m_1}\|_{L^2_G}^2\right) A_2^{\ell-\frac{1}{2}} \langle I_2\rangle \|u_{I_2,m_2}\|_{L^2_G}^2
	 \left(\sum_{(I_1,m_1)\in 2^{\mathbb{Z}}\times\mathbb{N}}    \|u_{I_1,m_1}\|_{L^{\infty}_G}^2\right)\,.
\end{multline*}
In only remains to use the definition~\eqref{def:Xk_rho} of $\mathcal{X}^k_1$ and the estimate proven in the above step $\sum_{(I_1,m_1)\in 2^{\mathbb{Z}}\times\mathbb{N}}    \|u_{I_1,m_1}\|_{L^{\infty}_G}^2\lesssim\|u_0\|_{\mathcal{X}^k_1}^2$ to get~\eqref{eq:AlmostBilin_zzz*_2}. 
\end{proof}

\section{Deterministic-probabilistic trilinear estimate}\label{sec:proba-deterministic}

In this section, we establish Theorem~\ref{th.bilinear_estimate}~(\textit{ii}). Let $\varepsilon_0>0$ such that $u_0\in \mathcal{X}^k_{1+\varepsilon_0}$ and $v,w\in L^{\infty}_TH^{\ell}_G$.
We recall that $z$, $v$ and $w$ have a decomposition
\[
    z=\sum_{A\in 2^{\mathbb{N}}} z_A=\sum_{A\in 2^{\mathbb{N}}}\sum_{\substack{(I_1,m_1)\in 2^{\mathbb{Z}}\times\mathbb{N}\\(m_1+1)I_1\sim A}} z_{I_1,m_1}\,,
\]
\[
    v=\sum_{B\in 2^{\mathbb{N}}} v_B=\sum_{B\in 2^{\mathbb{N}}}\sum_{\substack{(I_2,m_2)\in 2^{\mathbb{Z}}\times\mathbb{N}\\(m_2+1)I_2\sim B}} v_{I_2,m_2} 
\]
and
\[
    w=\sum_{C\in 2^{\mathbb{N}}} w_C=\sum_{C\in 2^{\mathbb{N}}}\sum_{\substack{(I_3,m_3)\in 2^{\mathbb{Z}}\times\mathbb{N}\\(m_3+1)I_3\sim B}} w_{I_3,m_3}\,. 
\]

Let $\varepsilon >0$ to be chosen later. Using Corollary~\ref{coro:derivativeSplit-2.3termes}, we get that for all $t$, there holds
\begin{multline}\label{eq:7_debut}
    \|z^{\omega}vw(t)\|_{H^{\ell}_G}^2 \lesssim \sum_{\substack{(\delta_1,\delta_2,\delta_3) \in D_3\\A,B,C: B,C\leqslant A}} A^{\ell+\varepsilon} \|(z^{\omega}_{A})^{\delta_1}(P_{\leq A}v_B)^{\delta_2}(P_{\leq A}w_C)^{\delta_3}(t)\|_{L^2_G}^2\\
     + \sum_{\substack{\delta_1 \in D_1\\A\in 2^{\mathbb{N}}}}A^{\varepsilon}\|(z^{\omega}_{A})^{\delta_1}(t)\|_{L^{\infty}_G}^2\|v(t)\|^2_{H^{\ell}_G}\|w(t)\|^2_{H^{\ell}_G}\,.
\end{multline}
 Using estimate~\eqref{eq.largeDeviation.v2} for the second term in the right hand side, we infer that outside of a set of probability $e^{-cR^2}$, there holds  
\begin{equation*}
    \sum_{A\in 2^{\mathbb{N}}}A^{\varepsilon}\|z^{\omega}_{A}\|^2_{L^q_TL^{\infty}_G}\lesssim_{\varepsilon} R^2T^{\frac{2}{q}} \|u_0\|^2_{\mathcal{X}^k_{1}}\,. 
\end{equation*}
Moreover, since applying the shifts $\delta_1\in D_1$ to every mode $z_{A}$ give equivalent estimates for the $L^p$ norms, this leads similarly to
\begin{equation*}
    \left\|\left(\sum_{\substack{\delta_1 \in D_1\\A\in 2^{\mathbb{N}}}}A^{\varepsilon}\|(z^{\omega}_{A})^{\delta_1}\|_{L^{\infty}_G}^2\right)^{1/2}\|v\|_{H^{\ell}_G}\|w\|_{H^{\ell}_G}\right\|_{L^q_T}
    \lesssim RT^{\frac{1}{q}} \|u_0\|_{\mathcal{X}^k_{1}}\|v\|_{L^{\infty}_TH^{\ell}_G}\|w\|_{L^{\infty}_TH^{\ell}_G}\,. 
\end{equation*}

For the sake of simplicity, we note that the shifted function $u^{\delta}$ is nothing but a shift of the indices $(I,m)\in 2^{\mathbb{Z}}\times\mathbb{N}$ of order at most one and a multiplication of every mode $(I,m)$ by a function of modulus at most $1$ in Fourier variable. Similarly, the projection $P_{\leq A}$ is a multiplication of every mode $(I,m)$ by a function of modulus at most $1$ in Fourier variable. Therefore, we assume without loss of generality that in the right-hand side of~\eqref{eq:7_debut} there is no shift ($\delta_1=\delta_2=\delta_3=\varnothing$) and no projection $P_{\leq A}$, up to applying the proof to $(P_{\leq A}v_B)^{\delta_1}$ instead of $v_B$ and doing a similar transformation for $w$. In order to estimate $ \|z^{\omega}vw\|_{L^q_TH^{\ell}_G}$, we therefore write
\begin{align*}
\left\|\left(\sum_{\substack{A,B,C: B,C\leqslant A}} A^{\ell+\varepsilon} \|z^{\omega}_{A}v_Bw_C\|_{L^2_G}^2\right)^{1/2}\right\|_{L^q_T}
	&=\left\|\sum_{\substack{A,B,C: B,C\leqslant A}} A^{\ell+\varepsilon} \|z^{\omega}_{A}v_Bw_C\|_{L^2_G}^2\right\|_{L^{q/2}_T}^{1/2}\,,
\end{align*}
and our aim is to prove that for our choice of $q$ and $\varepsilon$, then with probability greater than $1-e^{-cR^2}$, we have
\begin{equation*}
\left\|\sum_{\substack{A,B,C: B,C\leqslant A}} A^{\ell+\varepsilon} \|z^{\omega}_{A}v_Bw_C\|_{L^2_G}^2\right\|_{L^{q/2}_T}
	\lesssim_{\varepsilon}R^2T^{\frac{2}{q}} \|u_0\|_{\mathcal{X}^k_{1+\varepsilon_0}}^2\|v\|_{L^{\infty}_TH^{\ell}_G}^2\|w\|_{L^{\infty}_TH^{\ell}_G}^2\,.
\end{equation*} 

 By homogeneity, it is enough to prove that  there exists $C_{\varepsilon}$ such that for every $R>0$, with probability greater than $1-e^{-cR^2}$, for every $v,w\in  L^{\infty}_TH^{\ell}_G$ satisfying $\|v\|_{L^{\infty}_TH^{\ell}_G}\leq 1$ and $\|w\|_{L^{\infty}_TH^{\ell}_G}\leq 1$,
\begin{equation}\label{eq:using_Lem47}
\left\|\sum_{\substack{A,B,C: B,C\leqslant A}} A^{\ell+\varepsilon} \|z^{\omega}_{A}v_Bw_C(t)\|_{L^2_G}^2\right\|_{L^{q/2}_T}
	\leq C_{\varepsilon}R^2T^{\frac{2}{q}} \|u_0\|_{\mathcal{X}^k_{1+\varepsilon_0}}^2\,.
\end{equation}

\subsection{Step 1: Decoupling \texorpdfstring{$z^{\omega}$}{z} from \texorpdfstring{$vw$}{vw}\texorpdfstring{ in~\eqref{eq:using_Lem47}}{}.}\label{subsec:preparatory}
We fix $t\in\mathbb{R}$ and $A,B,C\in 2^{\mathbb{N}}$ such that $B,C\leq A$. Then we use the Plancherel formula to get 
\begin{align*}
\|z^{\omega}_{A}v_Bw_C(t)\|_{L^2_G}^2
	&=\int\mathrm{d}y\int \mathrm{d}x(z_A^{\omega}\overline{z_A^{\omega}})(x,y)(v_B\overline{v_B}w_C\overline{w_C})(x,y)\\
	&=\int \mathrm{d}\eta\int\mathrm{d}x (\widehat{z_A^{\omega}}*\widehat{\overline{z_A^{\omega}}})(x,\eta)(\widehat{v_B}*\widehat{\overline{v_B}}*\widehat{w_C}*\widehat{\overline{w_C}})(x,\eta)\,.
\end{align*}
We use decomposition~\eqref{def:BuildingBlock}. For $(I,m)\in 2^{\mathbb{Z}}\times\mathbb{N}$, we denote
\[
\mathcal{F}_{y\to\eta}(z^{\omega}_{I,m})(t,x,\eta)=f_{I,m}^{\omega}(t,\eta)h_{m}(\sqrt{|\eta|}x)\,,
\]
where $f_{I,m}^{\omega}(t,\eta)=X_{I,m}(\omega)e^{-t(2m+1)|\eta|}f_{I,m}(\eta)$ and $f_{I,m}(\eta)=f_m(\eta){\bf 1}_{|\eta|\in [I,2I]}$. Similarly for $v$, we write
\[
\mathcal{F}_{y\to\eta}(v_{I,m})(t,x,\eta)=g_{I,m}^{\omega}(t,\eta)h_{m}(\sqrt{|\eta|}x)\,,
\]
where the dependence of $g_{I,m}^{\omega}$ along the variables $t$ and $\omega$ is not explicit. For $w$ we do the same by making use of functions $\widetilde{g}^{\omega}_{I,m}$. Then we expand everything:
\begin{multline*}
\|z^{\omega}_{A}v_Bw_C(t)\|_{L^2_G}^2
	=\int \mathrm{d}\eta \int\mathrm{d}x 	\int\mathrm{d}\eta_1 \int\mathrm{d}\eta_2 \int\mathrm{d}\eta_3 \int\mathrm{d}\eta_4\\ 
	\sum_{\substack{m_1,m'_1\in\mathbb{N}\\(m_1+1)I_1,(m'_1+1)I'_1\sim A}}f_{I_1,m_1}^{\omega}(t,\eta_1)\overline{f_{I'_1,m'_1}^{\omega}(t,\eta-\eta_1)}
	h_{m_1}(\sqrt{|\eta_1|}x)h_{m'_1}(\sqrt{|\eta-\eta_1|}x)\\
	 \sum_{\substack{m_2,m'_2\in\mathbb{N}\\(m_2+1)I_2,(m'_2+1)I'_2\sim B}}g_{I_2,m_2}^{\omega}(t,\eta_2)\overline{g_{I'_2,m'_2}^{\omega}(t,\eta_3)}
	 h_{m_2}(\sqrt{|\eta_2|}x)h_{m'_2}(\sqrt{|\eta_3|}x)\\
	\sum_{\substack{m_3,m'_3\in\mathbb{N}\\(m_3+1)I_3,(m'_3+1)I'_3\sim C}} \widetilde{g}_{I_3,m_3}^{\omega}(t,\eta_4)\overline{\widetilde{g}_{I'_3,m'_3}^{\omega}(t,\eta-\eta_2-\eta_3-\eta_4)}
	h_{m_3}(\sqrt{|\eta_4|}x)h_{m'_3}(\sqrt{|\eta-\eta_2-\eta_3-\eta_4|}x)\,.
\end{multline*}
Our aim is to apply the probabilistic decoupling to the series involving products $f_{I_1,m_1}^{\omega}\overline{f_{I'_1,m'_1}^{\omega}}$. However, since $v$ and $w$ may depend on $\omega$, we first isolate the terms involving $g^{\omega}_{I_2,m_2}$, $g^{\omega}_{I'_2,m'_2}$, $\widetilde{g}^{\omega}_{I_3,m_3}$ and $\widetilde{g}^{\omega}_{I'_3,m'_3}$.

The expanded formula is rather long, we reorganize it as follows. We define
\begin{multline}\label{def:I}
\mathbb{J}_{\{m_i,m'_i\}}(\eta,\eta_1,\eta_2,\eta_3,\eta_4)\coloneqq\int\mathrm{d}x h_{m_1}(\sqrt{|\eta_1|}x)h_{m'_1}(\sqrt{|\eta-\eta_1|}x)h_{m_2}(\sqrt{|\eta_2|}x)h_{m'_2}(\sqrt{|\eta_3|}x)\\
	h_{m_3}(\sqrt{|\eta_4|}x)h_{m'_3}(\sqrt{|\eta-\eta_2-\eta_3-\eta_4|}x)\,.
\end{multline}
We define the ``random'' part
\begin{multline}\label{def:J}
{\bf J}^{\omega}_{I_2,I'_2,I_3,I'_3,m_2,m'_2,m_3,m'_3}(t,\eta,\eta_2,\eta_3,\eta_4)
	\coloneqq  \int\mathrm{d}\eta_1 \sum_{\substack{(I_1,m_1),(I'_1,m'_1)\in 2^{\mathbb{Z}}\times\mathbb{N}\\(m_1+1)I_1,(m'_1+1)I'_1\sim A}}
	f_{I_1,m_1}^{\omega}(t,\eta_1)\overline{f_{I'_1,m'_1}^{\omega}(t,\eta-\eta_1)}\\
	\mathbb{J}_{\{m_i,m'_i\}}(\eta,\eta_1,\eta_2,\eta_3,\eta_4)
	|\eta_2\eta_3\eta_4(\eta-\eta_2-\eta_3-\eta_4)|^{1/4}\\
	{\bf 1 }_{|\eta_2|\in[I_2,2I_2]}{\bf 1 }_{|\eta_3|\in[I'_2,2I'_2]}{\bf 1 }_{|\eta_4|\in
[I_3,2I_3]}{\bf 1 }_{|\eta-\eta_2-\eta_3-\eta_4|\in
[I'_3,2I'_3]
}\,,
\end{multline}
whereas the ``deterministic'' part is written
\[
{\bf K}^{\omega}_{I_2,I'_2,I_3,I'_3,m_2,m'_2,m_3,m'_3}(t,\eta,\eta_2,\eta_3,\eta_4)
	=\frac{g_{I_2,m_2}^{\omega}(t,\eta_2)}{|\eta_2|^{1/4}}\frac{\overline{g_{I'_2,m'_2}^{\omega}(t,\eta_3)}}{|\eta_3|^{1/4}}\\
	\frac{\widetilde{g}_{I_3,m_3}^{\omega}(t,\eta_4)}{|\eta_4|^{1/4}}\frac{\overline{\widetilde{g}_{I'_3,m'_3}^{\omega}(t,\eta-\eta_2-\eta_3-\eta_4)}}{|\eta-\eta_2-\eta_3-\eta_4|^{1/4}}\,.
\]
Then the expanded formula becomes
\begin{multline*}
\|z^{\omega}_{A}v_Bw_C(t)\|_{L^2_G}^2
	=\int \mathrm{d}\eta \int\mathrm{d}\eta_2 \int\mathrm{d}\eta_3 \int\mathrm{d}\eta_4 \sum_{\substack{(I_2,m_2),(I'_2,m'_2)\in 2^{\mathbb{Z}}\times\mathbb{N}\\(m_2+1)I_2,(m'_2+1)I'_2\sim B}}\sum_{\substack{(I_3,m_3),(I'_3,m'_3)\in 2^{\mathbb{Z}}\times\mathbb{N} \\(m_3+1)I_3,(m'_3+1)I'_3\sim C}}\\
	{\bf J}^{\omega}_{A,I_2,I'_2,I_3,I'_3,m_2,m'_2,m_3,m'_3}(t,\eta,\eta_2,\eta_3,\eta_4)\\
	{\bf K}^{\omega}_{I_2,I'_2,I_3,I'_3,m_2,m'_2,m_3,m'_3}(t,\eta,\eta_2,\eta_3,\eta_4)\,.
\end{multline*}
Taking the absolute value, this leads by summation to a condensed expression
\[
\left\|\sum_{\substack{A,B,C: B,C\leqslant A}} A^{\ell+\varepsilon} \|z^{\omega}_{A}v_Bw_C(t)\|_{L^2_G}^2\right\|_{L^{q/2}_T}
	\leq \|A^{\ell+\varepsilon}{\bf J}^{\omega}{\bf K}^{\omega}\|_{L^p_{\psi}}\,,
\]
where, $\psi$ stands for a huge tuple 
\[\psi=(T,A,B,C,I_2,I'_2,I_3,I'_3,m_2,m'_2,m_3,m'_3,\eta,\eta_2,\eta_3,\eta_4)\,.
\]
Moreover,  we take $L^p_{\psi}=L^{p_1}_{\psi_1}\dots L^{p_d}_{\psi_d}$ norm with the exponents 
\[
L^p_{\psi}=L^{q/2}_TL^1_AL^1_{B,C,I_2,I'_2,I_3,I'_3,m_2,m'_2,m_3,m'_3,\eta,\eta_2,\eta_3,\eta_4}\,,
\]
where the $L^{p_j}$ norm implicitely denotes the $\ell^{p_j}$ norm when we consider sequence spaces. It is now possible to take successive Hölder inequalities, as one easily checks that for $\frac{1}{q_1}+\frac{1}{q'_1}=\frac{1}{p_1}$ and $\frac{1}{q_2}+\frac{1}{q'_2}=\frac{1}{p_2}$, we have
\[
\|fg\|_{L^{p_1}_xL^{p_2}_y}
	\leq \Big\|\|f\|_{L^{q_2}_y}\|g\|_{L^{q'_2}_y}\Big\|_{L^{p_1}_x}
	\leq \|f\|_{L^{q_1}_xL^{q_2}_y}\|g\|_{L^{q'_1}_xL^{q'_2}_y}\,.
\]
Our purpose is to isolate the $L^{\infty}_TH^{\ell}_G$ norm of $vw$ by using Hölder inequalities. For this, we note that it better to split between $(BC)^{\ell}{\bf K}^{\omega}$ and $A^{\ell+\varepsilon}{\bf J}^{\omega}(BC)^{-\ell}$. As we will see below, our actual splitting also involves some extra other terms 
 which aim at balancing the two terms and make every series converge.

We fix $p_1$ and $p_2$ to be very large exponents in $[2,\infty)$ to be chosen later, with respective conjugate exponents denoted $p'_1$ and $p'_2$: $\frac{1}{p_1}+\frac{1}{p'_1}=1$ and $\frac{1}{p_2}+\frac{1}{p'_2}=1$. For the term involving ${\bf J}^{\omega}$, we choose the exponents
\[
L^{p(J)}_{\psi}=L^{q/2}_TL^{1}_AL^{p_2}_{B,C}L^2_{I_2,I'_2,I_3,I'_3,m_2,m'_2,m_3,m'_3}L^{p_1}_{\eta,\eta_2,\eta_3,\eta_4}
\]
and for the term involving ${\bf K}^{\omega}$, we choose the exponents
\[
L^{p(K)}_{\psi}=L^{\infty}_TL^{\infty}_AL^{p_2'}_{B,C}L^{2}_{I_2,I'_2,I_3,I'_3,m_2,m'_2,m_3,m'_3}L^{p_1'}_{\eta,\eta_2,\eta_3,\eta_4}\,.
\]
The reason why we introduce the exponents $p_1$ and $p_2$ instead of taking them directly equal to $\infty$ is that we need finite exponents in order to apply the probabilistic decoupling to ${\bf J}^{\omega}$.

Our choice of exponents is compatible with the assumptions for applying successive Hölder inequalities, so that
\begin{multline}\label{eq:decoupling}
\left\|\sum_{\substack{A,B,C: B,C\leqslant A}} A^{\ell+\varepsilon} \|z^{\omega}_{A}v_Bw_C\|_{L^2_G}^2\right\|_{L^{q/2}_T}\\
	\leq \|A^{\ell+\varepsilon}(BC)^{-\ell+\varepsilon_2}(I_2I'_2I_3I'_3)^{1/2-\varepsilon_1}{\bf J}^{\omega}\|_{L^{p(J)}_{\psi}}
	\|(BC)^{\ell-\varepsilon_2}(I_2I'_2I_3I'_3)^{-1/2+\varepsilon_1}{\bf K}^{\omega}\|_{L^{p(K)}_{\psi}}\,.
\end{multline}
The exponents $\varepsilon_1,\varepsilon_2>0$ will be chosen later depending on the following step.
 
\subsection{Step 2: Evaluating the deterministic part \texorpdfstring{${\bf K}^{\omega}$ in~\eqref{eq:decoupling}}{K}.}\label{subsec:deterministic}

We first estimate the term
\[\|(BC)^{\ell-\varepsilon_2}(I_2I'_2I_3I'_3)^{-1/2+\varepsilon_1}{\bf K}^{\omega}\|_{L^{p(K)}_{\psi}}
\]
with
\[
L^{p(K)}_{\psi}=L^{\infty}_AL^{\infty}_TL^{p_2'}_{B,C}L^{2}_{I_2,I'_2,I_3,I'_3,m_2,m'_2,m_3,m'_3}L^{p_1'}_{\eta,\eta_2,\eta_3,\eta_4}\,.
\]
For the intuition, we should keep in mind that $p'_1,p'_2\approx 1$. By definition, of ${\bf K}^{\omega}$, we have
\begin{multline*}
|{\bf K}^{\omega}_{I_2,I'_2,I_3,I'_3,m_2,m'_2,m_3,m'_3}(t,\eta,\eta_2,\eta_3,\eta_4)|
	=\frac{|g^{\omega}_{I_2,m_2}(t,\eta_2)|}{|\eta_2|^{1/4}}\frac{|g^{\omega}_{I'_2,m'_2}(t,\eta_3)|}{|\eta_3|^{1/4}}\\
	\frac{|\widetilde{g}^{\omega}_{I_3,m_3}(t,\eta_4)|}{|\eta_4|^{1/4}}\frac{|\widetilde{g}^{\omega}_{I'_3,m'_3}(t,\eta-\eta_2-\eta_3-\eta_4)|}{|\eta-\eta_2-\eta_3-\eta_4|^{1/4}}\,,
\end{multline*}
therefore
\[
\|{\bf K}^{\omega}\|_{L^{p'_1}_{\eta,\eta_2,\eta_3,\eta_4}}
	=\left\|\left(\frac{|g^{\omega}_{I_2,m_2}|}{|\cdot|^{1/4}}\right)^{p'_1}* 
	\left(\frac{|g^{\omega}_{I'_2,m'_2}|}{|\cdot|^{1/4}}\right)^{p'_1}*
	 \left(	\frac{|\widetilde{g}^{\omega}_{I_3,m_3}|}{|\cdot|^{1/4}}\right)^{p'_1}*
	 \left(\frac{|\widetilde{g}^{\omega}_{I'_3,m'_3}|}{|\cdot|^{1/4}}\right)^{p'_1}(t,\eta)\right\|_{L^{1}_\eta}^{1/p'_1}\,.
\]
We apply Young's inequality 
\[
\|{\bf K}^{\omega}\|_{L^{p'_1}_{\eta,\eta_2,\eta_3,\eta_4}}
	\leq\left\|\left(\frac{|g^{\omega}_{I_2,m_2}|}{|\cdot|^{1/4}}\right)^{p'_1}* 
	\left(\frac{|g^{\omega}_{I'_2,m'_2}|}{|\cdot|^{1/4}}\right)^{p'_1}(t,\eta)\right\|_{L^1_{\eta}}^{1/p'_1}
	 \left\|\left(	\frac{|\widetilde{g}^{\omega}_{I_3,m_3}|}{|\cdot|^{1/4}}\right)^{p'_1}*
	 \left(\frac{|\widetilde{g}^{\omega}_{I'_3,m'_3}|}{|\cdot|^{1/4}}\right)^{p'_1}(t,\eta)\right\|_{L^1_\eta}^{1/p'_1}\,.
\]
Applying Young's inequality again, we get
\[
\left\|\left(\frac{|g^{\omega}_{I_2,m_2}|}{|\cdot|^{1/4}}\right)^{p'_1}* \left(\frac{|g^{\omega}_{I'_2,m'_2}|}{|\cdot|^{1/4}}\right)^{p'_1}(t,\eta)\right\|_{L^1_{\eta}}^{1/p'_1}
	\leq \left\|\left(\frac{|g^{\omega}_{I_2,m_2}|}{|\cdot|^{1/4}}\right)^{p'_1}(t,\eta)\right\|_{L^{1}_{\eta}}^{1/p'_1}
	 \left\|\left(\frac{|g^{\omega}_{I'_2,m'_2}|}{|\cdot|^{1/4}}\right)^{p'_1}(t,\eta)\right\|_{L^{1}_{\eta}}^{1/p'_1}\,.
\]
But since
\begin{equation*}
\left\|\left(\frac{|g^{\omega}_{I'_2,m'_2}|}{|\cdot|^{1/4}}\right)^{p'_1}(t,\eta)\right\|_{L^{1}_{\eta}}^{1/p'_1}
	=\left\|\frac{|g^{\omega}_{I'_2,m'_2}|}{|\cdot|^{1/4}}(t,\eta)\right\|_{L^{p'_1}_{\eta}}\,,
\end{equation*}
the Hölder inequality with the choice of exponents $\frac{1}{p'_1}=\frac{1}{2}+\frac{p_1-2}{2p_1}$ leads to
\begin{align*}
 \left\|\left(\frac{|g^{\omega}_{I_2,m_2}|}{|\cdot|^{1/4}}\right)^{p'_1}(t,\eta)\right\|_{L^{1}_{\eta}}^{1/p'_1} 
 	&\leq  \left\|\frac{|g^{\omega}_{I_2,m_2}|}{|\cdot|^{1/4}}(t,\eta)\right\|_{L^{2}_{\eta}} \left\|{\bf 1}_{|\eta|\in[I_2,2I_2]}\right\|_{L^{2p_1/(p_1-2)}_{\eta}}\\
 	&\leq \|v_{I_2,m_2}(t)\|_{L^2_G}(I_2)^{(p_1-2)/(2p_1)}\,.
\end{align*}
By symmetry of the roles of $I_2,m_2$ and $I'_2,m'_2$, we deduce
\[
\left\|\left(\frac{|g^{\omega}_{I_2,m_2}|}{|\cdot|^{1/4}}\right)^{p'_1}* \left(\frac{|g^{\omega}_{I'_2,m'_2}|}{|\cdot|^{1/4}}\right)^{p'_1}(t,\eta)\right\|_{L^1_{\eta}}^{1/p'_1} 
	\leq \|v_{I_2,m_2}(t)\|_{L^2_G}\|v_{I'_2,m'_2}(t)\|_{L^2_G}\left(I_2I'_2\right)^{(p_1-2)/(2p_1)}\,.
\]
By symmetry of the roles of $u$ and $v$, we conclude
\[
\|{\bf K}^{\omega}\|_{L^{p_1'}_{\eta,\eta_2,\eta_3,\eta_4}}
	\leq \|v_{I_2,m_2}(t)\|_{L^2_G}\|v_{I'_2,m'_2}(t)\|_{L^2_G}\|w_{I_3,m_3}(t)\|_{L^2_G}\|w_{I'_3,m'_3}(t)\|_{L^2_G}\left(I_2I'_2I_3I'_3\right)^{(p_1-2)/(2p_1)}\,.
\]
Note that $\frac{p_1-2}{2p_1}=\frac{1}{2}-\frac{1}{p_1}=\frac{1}{2}-\varepsilon_1$ where $\varepsilon_1=\frac{1}{p_1}$ is very small. This leads to 
\begin{multline*}
\|(I_2I'_2I_3I'_3)^{-1/2+\varepsilon_1}{\bf K}^{\omega}\|_{L^{2}_{I_2,I'_2,I_3,I'_3,m_2,m'_2,m_3,m'_3}L^{p_1'}_{\eta, \eta_2,\eta_3,\eta_4}}\\
	\leq \Bigg(\sum_{\substack{m_2,m'_2,m_3,m'_3,I_2,I'_2,I_3,I'_3\\(m_2+1)I_2,(m'_2+1)I'_2\sim B\\(m_3+1)I_3,(m'_3+1)I'_3\sim C}}
	\|v_{I_2,m_2}(t)\|_{L^2_G}^2\|v_{I'_2,m'_2}(t)\|_{L^2_G}^2
	\|w_{I_3,m_3}(t)\|_{L^2_G}^2\|w_{I'_3,m'_3}(t)\|_{L^2_G}^2\Bigg)^{1/2}\\
	=\|v_B(t)\|_{L^2_G}^2\|w_C(t)\|_{L^2_G}^2\,.
\end{multline*}
Finally we get
\begin{multline*}
\|(BC)^{\ell-\varepsilon_2}(I_2I'_2I_3I'_3)^{-1/2+\varepsilon_1}{\bf K}^{\omega}\|_{L^{p_2'}_{B,C}L^{2}_{I_2,I'_2,I_3,I'_3,m_2,m'_2,m_3,m'_3}L^{p_1'}_{\eta,\eta_2,\eta_3,\eta_4}}\\
	\leq \left(\sum_{B,C}(BC)^{(\ell-\varepsilon_2)p'_2}\|v_B(t)\|_{L^2_G}^{2p'_2}\|w_C(t)\|_{L^2_G}^{2p'_2}\right)^{1/p'_2}\,.
\end{multline*}
The exponent $p'_2$ is greater than $1$, which could be troublesome, but we recall that we have chosen $v$ and $w$ satisfying $\|v\|_{L^{\infty}_TH^{\ell}_G}\leq 1$ and $\|w\|_{L^{\infty}_TH^{\ell}_G}\leq 1$ in the assumption of the desired estimate~\eqref{eq:using_Lem47}. Therefore, for every $t,B,C$, we have $\|v_B(t)\|_{L^2_G}\leq 1$ and $\|w_C(t)\|_{L^2_G}\leq 1$, leading to
\begin{multline*}
\|(BC)^{\ell-\varepsilon_2}(I_2I'_2I_3I'_3)^{-1/2+\varepsilon_1}{\bf K}^{\omega}\|_{L^{p_2'}_{B,C}L^{2}_{I_2,I'_2,I_3,I'_3,m_2,m'_2,m_3,m'_3}L^{p_1'}_{\eta ,\eta_2,\eta_3,\eta_4}}\\
	\leq \left(\sum_{B,C}(BC)^{(\ell-\varepsilon_2)p'_2}\|v_B(t)\|_{L^2_G}^{2}\|w_C(t)\|_{L^2_G}^{2}\right)^{1/p'_2}\,.
\end{multline*}
We choose $\varepsilon_2$ such that $(\ell-\varepsilon_2)p'_2=\ell$, i.e.\@ $\varepsilon_2=\frac{\ell}{p_2}$. Then the upper bound is bounded by
\begin{align*}
\|(BC)^{\ell-\varepsilon_2}(I_2I'_2I_3I'_3)^{-1/2+\varepsilon_1} 	&{\bf K}^{\omega}\|_{L^{p_2'}_{B,C}L^{2}_{I_2,I'_2,I_3,I'_3,m_2,m'_2,m_3,m'_3}L^{p_1'}_{\eta ,\eta_2,\eta_3,\eta_4}}\\
	&\leq \left(\sum_{B,C}\|v_B(t)\|_{H^{\ell}_G}^{2}\|w_B(t)\|_{H^{\ell}_G}^{2}\right)^{1/p'_2}\\
	&\leq (\|v(t)\|_{H^{\ell}_G}\|w(t)\|_{H^{\ell}_G})^{2/p'_2}\\
	&\leq 1\,.
\end{align*}
Taking the $L^{\infty}_TL^{\infty}_A$ norm, we conclude that
\[
\|(BC)^{\ell-\varepsilon_2}(I_2I'_2I_3I'_3)^{-1/2+\varepsilon_1}{\bf K}^{\omega}\|_{L^{p(K)}_{\psi}}
	\leq 1\,.
\]

\subsection{Step 3: Evaluating the random part \texorpdfstring{${\bf J}^{\omega}$ in~\eqref{eq:decoupling}}{J}.}\label{subsec:probabilistic}

We now estimate the term 
\[\|A^{\ell+\varepsilon}(BC)^{-\ell+\varepsilon_2}(I_2I'_2I_3I'_3)^{1/2-\varepsilon_1}{\bf J}^{\omega}\|_{L^{p(J)}_{\psi}}\]
with $\varepsilon_1=\frac{1}{p_1}$, $\varepsilon_2=\frac{\ell}{p_2}$ and
\[
L^{p(J)}_{\psi}=L^{q/2}_TL^{1}_AL^{p_2}_{B,C}L^2_{I_2,I'_2,I_3,I'_3,m_2,m'_2,m_3,m'_3}L^{p_1}_{\eta,\eta_2,\eta_3,\eta_4}\,.
\]

\noindent\textit{Step 3.1: probabilistic decoupling.} We first apply the probabilistic decoupling: let us check that the assumptions in order to apply Corollary~\ref{coro:WienerChaos}~(\textit{iii}) are met. It is now more convenient to write
\[
A^{\ell+\varepsilon}(BC)^{-\ell+\varepsilon_2}(I_2I'_2I_3I'_3)^{1/2-\varepsilon_1}{\bf J}^{\omega}
	=\sum_{\substack{m_1,m'_1,I_1,I'_1\\(m_1+1)I_1,(m'_1+1)I'_1\sim A}}X_{I_1,m_1}(\omega)\overline{X_{I'_1,m'_1}(\omega)}\Psi_{(I_1,m_1),(I'_1,m'_1)}(\psi)\,,
\]
where from the definition~\eqref{def:J} of ${\bf J}^{\omega}$, we have
\begin{multline}\label{def:psi}
\Psi_{(I_1,m_1),(I'_1,m'_1)}(\psi)=A^{\ell+\varepsilon}(BC)^{-\ell+\varepsilon_2}(I_2I'_2I_3I'_3)^{1/2-\varepsilon_1}\int\mathrm{d}\eta_1 f_{I_1,m_1}(t,\eta_1)\overline{f_{I'_1,m'_1}}(t,\eta-\eta_1) \\
	\mathbb{J}_{\{m_i,m'_i\}}(\eta,\eta_1,\eta_2,\eta_3,\eta_4)
	|\eta_2\eta_3\eta_4(\eta-\eta_2-\eta_3-\eta_4)|^{1/4}\\
	{\bf 1 }_{|\eta_2|\in[I_2,2I_2]}{\bf 1 }_{|\eta_3|\in[I'_2,2I'_2]}{\bf 1 }_{|\eta_4|\in
[I_3,2I_3]}{\bf 1 }_{|\eta-\eta_2-\eta_3-\eta_4|\in
[I'_3,2I'_3]
}\,,
\end{multline}
and we recall that $\mathbb{J}$ is defined in~\eqref{def:I}.

Applying Corollary~\ref{coro:WienerChaos}~(\textit{iii}) with  
\[
\psi_-=(T,A,B,C,I_2,I'_2,I_3,I'_3,m_2,m'_2,m_3,m'_3),
\quad
\psi_+=(\eta,\eta_2,\eta_3,\eta_4)\,,
\]
\[
L^{p_-(J)}_{\psi_-}=L^{q/2}_TL^1_AL^{p_2}_{B,C}L^2_{I_2,I'_2,I_3,I'_3,m_2,m'_2,m_3,m'_3}\,,
\quad
L^{p_+(J)}_{\psi_+}=L^{p_1}_{\eta,\eta_2,\eta_3,\eta_4}\,,
\]
since
\begin{multline*}
\|A^{\ell+\varepsilon}(BC)^{-\ell+\varepsilon_2}		(I_2I'_2I_3I'_3)^{1/2-\varepsilon_1}{\bf J}^{\omega}\|_{L^{p(J)}_{\psi}}\\
	=\left\|\sum_{\substack{m_1,m'_1,I_1,I'_1\\(m_1+1)I_1,(m'_1+1)I'_1\sim A}}X_{I_1,m_1}(\omega)\overline{X_{I'_1,m'_1}(\omega)}\Psi_{(I_1,m_1),(I'_1,m'_1)}(\psi)\right\|_{L^{p(J)}_{\psi}}\,,
\end{multline*}
we get that outside of a set of probability $e^{-cR^2}$, there holds
\begin{multline}\label{eq:Jdecoupled}
\|A^{\ell+\varepsilon}(BC)^{-\ell+\varepsilon_2}		(I_2I'_2I_3I'_3)^{1/2-\varepsilon_1}{\bf J}^{\omega}\|_{L^{p(J)}_{\psi}}\\
	\leq R^2\left\|\left(\sum_{\substack{m_1,m'_1,I_1,I'_1\\(m_1+1)I_1,(m'_1+1)I'_1\sim A}}\|\Psi_{(I_1,m_1),(I'_1,m'_1)}(\psi)\|_{L^{p_+(J)}_{\psi_+}}^2\right)^{1/2}\right\|_{L^{p_-(J)}_{\psi_-}}\\
	+R^2\left\|\sum_{\substack{m_1,I_1\\(m_1+1)I_1\sim A}} \|\Psi_{(I_1,m_1),(I_1,m_1)}(\psi)\|_{L^{p_+(J)}_{\psi_+}}\right\|_{L^{p_-(J)}_{\psi_-}}\,.
\end{multline}

\noindent\textit{Step 3.2: deterministic trilinear estimates.} 
Now, we apply the deterministic trilinear estimates from Corollary~\ref{coro.product.rescaled.hermite.trilinear} in order to gain $\frac{1}{2}$ additional derivatives, \textit{i.e.}\@ $\frac{1}{2}$ powers of $A$. 

Let us fix $I_1,I'_1,m_1,m'_1$ and $\psi$. We estimate $\Psi_{(I_1,m_1),(I'_1,m'_1)}(\psi)$ defined in~\eqref{def:psi}. We start with $\mathbb{J}$ defined in~\eqref{def:I} as
\begin{multline*}
\mathbb{J}_{\{m_i,m'_i\}}(\eta,\eta_1,\eta_2,\eta_3,\eta_4)
	=\int\mathrm{d}x h_{m_1}(\sqrt{|\eta_1|}x)h_{m'_1}(\sqrt{|\eta-\eta_1|}x)h_{m_2}(\sqrt{|\eta_2|}x)h_{m'_2}(\sqrt{|\eta_3|}x)\\
	h_{m_3}(\sqrt{|\eta_4|}x)h_{m'_3}(\sqrt{|\eta-\eta_2-\eta_3-\eta_4|}x)\,.
\end{multline*}
We apply Cauchy-Schwarz' inequality
\begin{multline*}
|\mathbb{J}_{\{m_i,m'_i\}}(\eta,\eta_1,\eta_2,\eta_3,\eta_4)|
	\leq \| h_{m_1}(\sqrt{|\eta_1|}x)h_{m_2}(\sqrt{|\eta_2|}x)h_{m_3}(\sqrt{|\eta_4|}x)\|_{L^2_x}\\
	\| h_{m'_1}(\sqrt{|\eta-\eta_1|}x)h_{m'_2}(\sqrt{|\eta_3|}x)	h_{m'_3}(\sqrt{|\eta-\eta_2-\eta_3-\eta_4|}x)\|_{L^2_x}\,.
\end{multline*}
We apply Corollary~\ref{coro.product.rescaled.hermite.trilinear} for $\alpha_1^2=|\eta_1|\in[I_1,2I_1] $, $\alpha_2^2=|\eta_2|\in[I_2,2I_2]$ and $\alpha_3^2=|\eta_4|\in
[I_3,2I_3]$:
\[
   \alpha_1\alpha_2\alpha_3 \|h_{m_1}(\alpha_{1} \cdot)h_{m_2}(\alpha_{2} \cdot) h_{m_3}(\alpha_{3} \cdot)\|_{L^2}^2 \lesssim  C_{\{I_i,m_i\}}^2\,,
\]
where
\[
C_{\{I_i,m_i\}}^2
	=\frac{\langle I_1\rangle(I_2I_3)^{1/4}}{A^{1/2}((2m_2+1)(2m_3+1))^{1/12}}\,.
\]
We also apply Corollary~\ref{coro.product.rescaled.hermite.trilinear} for $\alpha_1^2=|\eta-\eta_1|\in[I'_1,2I'_1]$, $\alpha_2^2=|\eta_3|\in[I'_2,2I'_2]$ and $\alpha_3^2=|\eta-\eta_2-\eta_3-\eta_4|\in
[I'_3,2I'_3]
$. We conclude that
\[
|\eta_1(\eta-\eta_1)\eta_2\eta_3\eta_4(\eta-\eta_2-\eta_3-\eta_4)|^{1/4}|\mathbb{J}_{\{m_i,m'_i\}}(\eta,\eta_1,\eta_2,\eta_3,\eta_4)|\\
	\lesssim C_{\{I_i,m_i\}}C_{\{I'_i,m'_i\}}\,.
\]
Moreover, since
\[
\int \mathrm{d}\eta_1\left|\frac{f_{I_1,m_1}(t,\eta_1)}{|\eta_1|^{1/4}}\frac{\overline{f_{I'_1,m'_1}(t,\eta-\eta_1)}}{|\eta-\eta_1|^{1/4}}\right|
	=\left(\frac{|f_{I_1,m_1}|}{|\cdot|^{1/4}}\right)*\left(\frac{|f_{I'_1,m'_1}|}{|\cdot|^{1/4}}\right)(t,\eta)\,,
\]
we get the estimate
\begin{multline*}
|\Psi_{(I_1,m_1),(I'_1,m'_1)}(\psi)|
	\leq A^{\ell+\varepsilon}(BC)^{-\ell+\varepsilon_2}(I_2I'_2I_3I'_3)^{1/2-\varepsilon_1}\\
	C_{\{I_i,m_i\}}C_{\{I'_i,m'_i\}}\left(\frac{|f_{I_1,m_1}|}{|\cdot|^{1/4}}\right)*\left(\frac{|f_{I'_1,m'_1}|}{|\cdot|^{1/4}}\right)(t,\eta)\\
	{\bf 1 }_{|\eta_2|\in[I_2,2I_2]}{\bf 1 }_{|\eta_3|\in[I'_2,2I'_2]}{\bf 1 }_{|\eta_4|\in
[I_3,2I_3]}{\bf 1 }_{|\eta-\eta_2-\eta_3-\eta_4|\in
[I'_3,2I'_3]
}\,.
\end{multline*}

\noindent\textit{Step 3.3: Integration bounds.} 
 It is now time to evaluate $\|\Psi_{(I_1,m_1),(I'_1,m'_1)}(\psi)\|_{L^{p_+(J)}_{\psi_+}}$ with
\[
L^{p_+(J)}_{\psi_+}=L^{p_1}_{\eta,\eta_2,\eta_3,\eta_4}\,.
\]

First, as indicators functions are bounded by $1$, we have
\begin{multline*}
\left\|{\bf 1 }_{|\eta_2|\in[I_2,2I_2]}{\bf 1 }_{|\eta_3|\in[I'_2,2I'_2]}
	{\bf 1 }_{|\eta_4|\in
[I_3,2I_3]}{\bf 1 }_{|\eta-\eta_2-\eta_3-\eta_4|\in
[I'_3,2I'_3]
}\right\|_{L^{p_1}_{\eta_2,\eta_3,\eta_4}}\\
	\leq \left\|{\bf 1 }_{|\eta_2|\in[I_2,2I_2]}{\bf 1 }_{|\eta_3|\in[I'_2,2I'_2]}{\bf 1 }_{|\eta_4|\in
[I_3,2I_3]}\right\|_{L^{1}_{\eta_2,\eta_3,\eta_4}}^{1/p_1}
	\leq \left(I_2I'_2I_3\right)^{1/p_1}\,.
\end{multline*}
By symmetry of roles, recalling that $\varepsilon_1=\frac{1}{p_1}$, we write that
\begin{equation*}
\left\|{\bf 1 }_{|\eta_2|\in[I_2,2I_2]}{\bf 1 }_{|\eta_3|\in[I'_2,2I'_2]}
	{\bf 1 }_{|\eta_4|\in
[I_3,2I_3]}{\bf 1 }_{|\eta-\eta_2-\eta_3-\eta_4|\in
[I'_3,2I'_3]
}\right\|_{L^{p_1}_{\eta_2,\eta_3,\eta_4}}\\
	\leq \left(I_2I'_2I_3I'_3\right)^{3\varepsilon_1/4}\,.
\end{equation*}
Hence we have discarded the integration over $\eta_1,\eta_2,\eta_3$ and the term ${\bf 1 }_{|\eta_2|\in[I_2,2I_2]} {\bf 1 }_{|\eta_3|\in[I'_2,2I'_2]}\times$ ${\bf 1 }_{|\eta_4|\in [I_3,2I_3]}{\bf 1 }_{|\eta-\eta_2-\eta_3-\eta_4|\in
[I'_3,2I'_3]
}$:
\begin{multline*}
\|\Psi_{(I_1,m_1),(I'_1,m'_1)}(\psi)\|_{L^{p_1}_{\eta_2,\eta_3,\eta_4}}
	\leq A^{\ell+\varepsilon}(BC)^{-\ell+\varepsilon_2}(I_2I'_2I_3I'_3)^{1/2-\varepsilon_1/4}\\
	C_{\{I_i,m_i\}}C_{\{I'_i,m'_i\}}\left(\frac{|f_{I_1,m_1}|}{|\cdot|^{1/4}}\right)*\left(\frac{|f_{I'_1,m'_1}|}{|\cdot|^{1/4}}\right)(\eta)\,.
\end{multline*}

Then, we integrate in $\eta$. Fixing $I_1,I'_1,m_1,m'_1$, we have thanks to Young's  inequality with $1+\frac{1}{p_1}=2\frac{p_1+1}{2p_1}$,
\begin{equation*}
\left\|\left(\frac{|f_{I_1,m_1}|}{|\cdot|^{1/4}}\right)*\left(\frac{|f_{I'_1,m'_1}|}{|\cdot|^{1/4}}\right)(\eta)\right\|_{L^{p_1}_{\eta}}
	\leq \left\|\left(\frac{|f_{I_1,m_1}|}{|\cdot|^{1/4}}\right)\right\|_{L^{2p_1/(p_1+1)}_{\eta}}
		\left\|\left(\frac{|f_{I'_1,m'_1}|}{|\cdot|^{1/4}}\right)\right\|_{L^{2p_1/(p_1+1)}_{\eta}}
\end{equation*}
and from Hölder's inequality with $\frac{p_1+1}{2p_1}=\frac{1}{2}+\frac{1}{2p_1}$,
\begin{align*}
\left\|\left(\frac{|f_{I_1,m_1}|}{|\cdot|^{1/4}}\right)*\left(\frac{|f_{I'_1,m'_1}|}{|\cdot|^{1/4}}\right)(t,\eta)\right\|_{L^{p_1}_{\eta}}
	&\leq (I_1I'_1)^{1/(2p_1)}\left\|\left(\frac{|f_{I_1,m_1}|}{|\cdot|^{1/4}}\right)(t,\eta)\right\|_{L^{2}_{\eta}}
		\left\|\left(\frac{|f_{I'_1,m'_1}|}{|\cdot|^{1/4}}\right)(t,\eta)\right\|_{L^{2}_{\eta}}\\
	&= (I_1I'_1)^{1/(2p_1)}\|z_{I_1,m_1}(t)\|_{L^2_G}\|z_{I'_1,m'_1}(t)\|_{L^2_G}\\
	&= (I_1I'_1)^{1/(2p_1)}\|u_{I_1,m_1}\|_{L^2_G}\|u_{I'_1,m'_1}\|_{L^2_G}\,,
\end{align*}
where in the latter equality we have used the $L^2$ isometry property of the linear flow. We have proven that
\begin{multline*}
\|\Psi_{(I_1,m_1),(I'_1,m'_1)}(\psi)\|_{L^{p_+(J)}_{\psi_+}}
	\leq A^{\ell+\varepsilon}(BC)^{-\ell+\varepsilon_2}(I_2I'_2I_3I'_3)^{1/2-\varepsilon_1/4}\\
	C_{\{I_i,m_i\}}C_{\{I'_i,m'_i\}}(I_1I'_1)^{1/(2p_1)}\|u_{I_1,m_1}\|_{L^2_G}\|u_{I'_1,m'_1}\|_{L^2_G}\,.
\end{multline*}
Recalling that $\varepsilon_1=\frac{1}{p_1}$ and replacing $C_{\{I_i,m_i\}}^2
	=\frac{\langle I_1\rangle(I_2I_3)^{1/4}}{A^{1/2}((2m_2+1)(2m_3+1))^{1/12}}$
	 by its value, this leads to
\begin{multline*}
\|\Psi_{(I_1,m_1),(I'_1,m'_1)}(\psi)\|_{L^{p_+(J)}_{\psi_+}}
		\leq A^{\ell+\varepsilon}(BC)^{-\ell+\varepsilon_2}(I_2I'_2I_3I'_3)^{1/2-\varepsilon_1/4}\\
	\frac{(\langle I_1\rangle\langle I'_1\rangle)^{1/2}(I_2I'_2I_3I'_3)^{1/8}}{A^{1/2}((2m_2+1)(2m'_2+1)(2m_3+1)(2m'_3+1))^{1/24}}
	(I_1I'_1)^{\varepsilon_1/2}\|u_{I_1,m_1}\|_{L^2_G}\|u_{I'_1,m'_1}\|_{L^2_G}\,.
\end{multline*}
\medskip

\noindent\textit{Step 3.4: decoupled summation.} We now consider the sums over $I_1,I'_1,m_1,m'_1$ in~\eqref{eq:Jdecoupled}.

For the second term in the right-hand side of \eqref{eq:Jdecoupled}, we take $I'_1=I_1$, $m'_1=m_1$ and incorporate the term $\langle I_1\rangle^{1+\varepsilon_1}$, which is the only term still dependent of $I_1,m_1$ in the upper bound. By summation, this leads to
\begin{equation*}
\sum_{\substack{m_1,I_1\\(m_1+1)I_1\sim A}}\langle I_1\rangle^{1+\varepsilon_1}\|u_{I_1,m_1}\|_{L^2_G}^2
	= \|u_{A}\|_{\mathcal{X}^0_{1+\varepsilon_1}}^2\,.
\end{equation*}
For the first term in the right-hand side of~\eqref{eq:Jdecoupled}, we also have
\[
\left(\sum_{\substack{m_1,m'_1,I_1,I'_1\\(m_1+1)I_1,(m'_1+1)I'_1\sim A}}\langle I_1\rangle^{1+\varepsilon_1}\langle I'_1\rangle^{1+\varepsilon_1}\|u_{I_1,m_1}\|_{L^2_G}^2\|u_{I'_1,m'_1}\|_{L^2_G}^2\right)^{1/2}
	= \|u_{A}\|_{\mathcal{X}^0_{1+\varepsilon_1}}^2\,.
\]
Therefore~\eqref{eq:Jdecoupled} becomes
\begin{multline*}
\|A^{\ell+\varepsilon}(BC)^{-\ell+\varepsilon_2}		(I_2I'_2I_3I'_3)^{1/2-\varepsilon_1}{\bf J}^{\omega}\|_{L^{p(J)}_{\psi}}\\
	\lesssim R^2\left\|A^{\ell-\frac{1}{2}+\varepsilon}(BC)^{-\ell+\varepsilon_2}\frac{(I_2I'_2I_3I'_3)^{1/2+1/8-\varepsilon_1/4}}
	{((2m_2+1)(2m'_2+1)(2m_3+1)(2m'_3+1))^{1/24}}
	\|u_{A}\|_{\mathcal{X}^0_{1+\varepsilon_1}}^2\right\|_{L^{p_-(J)}_{\psi_-}}\,,
\end{multline*}
where
\[
L^{p_-(J)}_{\psi_-}=L^{q/2}_TL^1_AL^{p_2}_{B,C}L^2_{I_2,I'_2,I_3,I'_3,m_2,m'_2,m_3,m'_3}\,.
\]
\medskip

\noindent\textit{Step 3.5: Hölder bounds.} 
Let us now compute the norm $L^2_{I_2,m_2}$. Thanks to the condition $1+(2m_2+1)I_2\in[B, 2B]$ implying $m_2\leq \frac{B}{I_2}$ and $I_2\leq B$, we have
\begin{align*}
\sum_{\substack{m_2,I_2\\(m_2+1)I_2\sim B}}\frac{I_2^{1+1/4-\varepsilon_1/2}}{(m_2+1)^{1/12}}
	&\leq \sum_{I_2\in 2^{\mathbb{Z}}: I_2\leq B}I_2^{1+1/4-\varepsilon_1/2}\sum_{m_2\leq \frac{B}{I_2}}\frac{1}{(m_2+1)^{1/12}}\\
	&\leq \sum_{I_2\in 2^{\mathbb{Z}}: I_2\leq B}I_2^{1+1/4-\varepsilon_1/2}\left(\frac{B}{I_2}\right)^{1-1/12}\\
	&\lesssim B^{1-1/12}\sum_{I_2\in 2^{\mathbb{Z}}: I_2\leq B}I_2^{1/4+1/12-\varepsilon_1/2}\,.
\end{align*}
The series over $I_2\leq 1$ is convergent since for $\varepsilon_1>0$ sufficiently small, the exponent $1/4+1/12-\varepsilon_1/2$ is positive. Moreover, one can see that the gain of $1/12$ is actually useless here in the argument. For the series over $I_2\geq 1$, we get a bound $B^{1/4+1/12-\varepsilon_1/2}$, so that
\[
\sum_{\substack{m_2,I_2\\(m_2+1)I_2\sim B}}\frac{I_2^{1+1/4-\varepsilon_1/2}}{(m_2+1)^{1/12}}
	\lesssim B^{1+1/4-\varepsilon_1/2}
	\lesssim B^{5/4}\,.
\]
We do the same for all the other indices $I'_2,m'_2,I_3,m_3,I'_3,m'_3$ and get
\begin{equation*}
\|A^{\ell+\varepsilon}(BC)^{-\ell+\varepsilon_2}		(I_2I'_2I_3I'_3)^{1/2-\varepsilon_1}{\bf J}^{\omega}\|_{L^{p(J)}_{\psi}}\\
	\lesssim R^2\left\|A^{\ell-1/2+\varepsilon}(BC)^{-\ell+\varepsilon_2+5/4}
		\|u_{A}\|_{\mathcal{X}^0_{1+\varepsilon_1}}^2\right\|_{L^{q/2}_TL^1_AL^{p_2}_{B,C}}\,.
\end{equation*}

We now take the $L^{p_2}$ norm over $B,C$. Since $\ell>\frac{3}{2}$, the exponent in front of the term $(BC)$ is $-\ell+\varepsilon_2+\frac 54<-\frac{3}{2}+\varepsilon_2+\frac{5}{4}$. When $\varepsilon_2<\frac{1}{4}$ is small, we see that this exponent is negative, so that the $L^{p_2}$ norm over $B$ and $C$ is convergent and bounded by some constant $C_0$. This leads to
\begin{equation*}
\|A^{\ell+\varepsilon}(BC)^{-\ell+\varepsilon_2}		(I_2I'_2I_3I'_3)^{1/2-\varepsilon_1}{\bf J}^{\omega}\|_{L^{p(J)}_{\psi}}\\
	\lesssim R^2\left\|A^{\ell-\frac{1}{2}+\varepsilon}\|u_{A}\|_{\mathcal{X}^0_{1+\varepsilon_1}}^2\right\|_{L^{q/2}_TL^1_A}\,.
\end{equation*}
We finally conclude that
\begin{equation*}
\|A^{\ell+\varepsilon}(BC)^{-\ell+\varepsilon_2}		(I_2I'_2I_3I'_3)^{1/2-\varepsilon_1}{\bf J}^{\omega}\|_{L^{p(J)}_{\psi}}\\
	\lesssim R^2\|\|u_0\|_{\mathcal{X}^{\ell-\frac{1}{2}+\varepsilon}_{1+\varepsilon_1}}^2\|_{L^{q/2}_T}\,.
\end{equation*}
It only remains to take the $L^{q/2}_T$ norm, but one can see that our upper bound does not depend on $T$ anymore, so this only adds a $T^{2/q}$ factor:
\begin{equation*}
\|A^{\ell+\varepsilon}(BC)^{-\ell+\varepsilon_2}		(I_2I'_2I_3I'_3)^{1/2-\varepsilon_1}{\bf J}^{\omega}\|_{L^{p(J)}_{\psi}}\\
	\lesssim R^2T^{2/q}\|u_0\|_{\mathcal{X}^{\ell-\frac{1}{2}+\varepsilon}_{1+\varepsilon_1}}^2\,.
\end{equation*}
When $\ell+\varepsilon<k+\frac{1}{2}$ (this is equivalent to taking $\varepsilon$ small enough), we conclude that
\[
    \|z^{\omega}vw\|^2_{L^q_TH^{\ell}_G} \lesssim_{\varepsilon} R^2T^{\frac{2}{q}}\|u_0\|_{\mathcal{X}^k_{1+\varepsilon_1}}^2\,,
\]
where $\varepsilon_1=\frac{1}{p_1}$ can be chosen arbitrarily small and in particular no greater than $\varepsilon_0$. Using the homogeneity, we remove the assumption $\|v\|_{L^{\infty}_TH^{\ell}_G}\leq 1$, $\|w\|_{L^{\infty}_TH^{\ell}_G}\leq 1$
and deduce that up to removing an extra set of probability not larger than $e^{-cR^2}$, inequality~\eqref{eq:Bilin_zv} holds. This concludes the proof of Theorem~\ref{th.bilinear_estimate}~(\textit{ii}).

\section{Local well-posedness}\label{sec:localWellPosedness}

This section is devoted to the proof of Theorem~\ref{th.main} . 

We fix $k \in (1,\frac{3}{2})$ and $u_0 \in H^k_G$. We assume that there exists $\varepsilon_0>0$ such that $u_0\in \mathcal{X}^k_{1+\varepsilon}$. We denote by  $u_0^{\omega}$ its associate randomization, defined by~\eqref{eq:randomizationMap}, and recall that we write $z^{\omega}(t) := e^{it\Delta_G}u_0^{\omega}$ for the solution to the linear flow~\eqref{eq:LS-G} associated to the initial data $u_0^{\omega}$. 
We seek for a solution $u$ to~\eqref{eq:NLSG} of the form
\[
    u(t)=z^{\omega}(t)+v(t)\,,
\]
where $v(0)=0$ and $v(t) \in H^{\ell}_G$ with $\ell \in (\frac{3}{2},k+\frac{1}{2})$. We will prove local well-posedness for $v\in\mathcal{C}^0\left([0,T],H^{\ell}_G\right)$ solving
\begin{equation}
    \label{eq:NLSGv}
    \begin{cases}
    i\partial_tv-\Delta_G  v=|z^{\omega}+v|^2(z^{\omega}+v)
\\
v(0)=0\,,
\end{cases}
\end{equation} 
thanks to a fixed point argument. We consider the map $\Phi : v\in\mathcal{C}^0\left([0,T],H^{\ell}_G\right) \mapsto \Phi(v)\in\mathcal{C}^0\left([0,T],H^{\ell}_G\right)$ defined by
\begin{equation}
    \label{eq:mapContraction}
    \Phi(v) : t\in[0,T] \mapsto -i\int_0^te^{i(t-t')\Delta_G} \left(|z^{\omega}+v|^{2}(z^{\omega}+v)\right)\,\mathrm{d}t'\,. 
\end{equation}
Observe that by the Duhamel formula, $v$ solves~\eqref{eq:NLSGv} if and only if $\Phi(v)=v$. Note that $v$ may depend on $\omega$.

We introduce the following set of initial data: 
\[
E_{R,T}=\{\omega\in\Omega \mid~\eqref{eq.1}, \eqref{eq.3}, \eqref{eq.4} \text{ and }~\eqref{eq.5} \text{ hold}\}\,,
\]
where
\begin{equation}
    \label{eq.1}
    \|(z^{\omega})^2\|_{L^{1}_TH^{\ell}_G}+\||z^{\omega}|^2\|_{L^{1}_TH^{\ell}_G} \leqslant TR^2 \|u_0\|_{\mathcal{X}^{k}_1}^2\,,
\end{equation}
\begin{equation}
    \label{eq.3}
    \||z^{\omega}|^2z^{\omega}\|_{L^1_TH^{\ell}_G} \leqslant TR^3 \|u_0\|_{\mathcal{X}^{k}_1}^3\,,
\end{equation}
\begin{equation}
    \label{eq.4}
    \|z^{\omega}vw\|_{L^1_TH^{\ell}_G} \leqslant  TR\|u_0\|_{\mathcal{X}^k_{1}}\|v\|_{L^{\infty}_TH^{\ell}_G}\|w\|_{L^{\infty}_TH^{\ell}_G}\text{ for all } v,w\in L^{\infty}_T{H}^{\ell}_G\,,
\end{equation}
\begin{equation}
    \label{eq.5}
    \|z^{\omega}\|_{L^2_TL^{\infty}_G}^2 \leqslant TR^2\|u_0\|^2_{\mathcal{X}^k_{1}}\,.
\end{equation}

We have the following estimate of $E_{R,T}$. 

\begin{lemma}
Let $u_0\in H^k_G$. Then there exists a constant $c>0$ which depends on the basis function $u_0$ of the randomization, such that for all $R,T>0$, $\mathbb{P}(\Omega \setminus E_{R,T}) \leqslant e^{-cR^2}$. 
\end{lemma}
\begin{proof} Outside a set of probability at most $e^{-cR^2}$ the bounds~\eqref{eq.1} and~\eqref{eq.3} follow from Theorem~\ref{th.bilinear_estimate}~(\textit{i}). Similarly, and~\eqref{eq.4} follows from Theorem~\ref{th.bilinear_estimate}~(\textit{ii}) with $u_0\in \mathcal{X}^{k}_{1}\subset\mathcal{X}^{k-\varepsilon}_{1+\varepsilon}$ for $\varepsilon$ chosen small enough so that $\ell<k-\varepsilon+\frac{1}{2}$. Moreover,  \eqref{eq.5} follows from Proposition~\ref{prop.random.integrability.improvement}.
\end{proof}

The key estimate in proving Theorem~\ref{th.main} is the following.  

\begin{proposition}[A priori estimate]\label{prop.contraction} Let $u_0^{\omega} \in E_{R,T}$ and $\ell \in (\frac{3}{2},k+\frac{1}{2})$. Then for any $v \in \mathcal{C}^0([0,T],H^{\ell}_G)$ there holds
\begin{equation}
    \label{eq.contraction}
    \|\Phi(v)\|_{L^{\infty}_TH^{\ell}_G} \lesssim T\left(\|v\|_{L^{\infty}_TH^{\ell}_G}^3 + (R\|u_0\|_{\mathcal{X}^k_1})^3\right)\,.
\end{equation}
Similarly, for any $v_1, v_2 \in \mathcal{C}^0([0,T],H^{\ell}_G)$ there holds
\begin{equation}
    \label{eq.contractionBis}
    \|\Phi(v_2)-\Phi(v_1)\|_{L^{\infty}_TH^{\ell}_G} \lesssim T\|v_2-v_1\|_{L^{\infty}_TH^{\ell}_G} \left((R\|u_0\|_{\mathcal{X}^k_1})^2+\|v_1\|^2_{L^{\infty}_TH^{\ell}_G} + \|v_2\|^2_{L^{\infty}_TH^{\ell}_G}\right)\,.
\end{equation}
\end{proposition}

\begin{proof}
Let $\omega \in E_{R,T}$. Estimates~\eqref{eq.contraction} and~\eqref{eq.contractionBis} reduce to multilinear estimates, as a consequence of the triangle inequality in the Duhamel formula~\eqref{eq:mapContraction} and the fact that $e^{it\Delta_G}$ is an isometry in $H^{\ell}_G$. Indeed, for $t \in[0,T]$, there holds
\begin{align*}
    \|\Phi(v)(t)\|_{H^{\ell}_G} &\leqslant \int_0^t\|e^{i(t-t')\Delta_G}|z^{\omega}+v|^2(z^{\omega}+v)(t')\|_{H^{\ell}_G}\,\mathrm{d}t' \\
    & \lesssim \int_0^t\||z^{\omega}+v|^2(z^{\omega}+v)(t')\|_{H^{\ell}_G}\,\mathrm{d}t' \\
    & \lesssim \||z^{\omega}+v|^2(z^{\omega}+v)\|_{L^1_TH^{\ell}_G}\,. 
\end{align*}
We expand the cubic term and bound
\begin{multline}
    \label{eq:contraction1terms}
    \|\Phi (v)\|_{L^{\infty}_TH^{\ell}_G} \lesssim \||z^{\omega}|^2z^{\omega}\|_{L^{1}_TH^{\ell}_G}  + \|(z^{\omega})^2  \bar v\|_{L^{1}_TH^{\ell}_G} + \||z^{\omega}|^2 v\|_{L^{1}_TH^{\ell}_G}\\
    + \|z^{\omega}|v|^2\|_{L^{1}_TH^{\ell}_G}  + \| z^{\omega} \bar v^2\|_{L^{1}_TH^{\ell}_G}  + T\||v|^2v\|_{L^{\infty}_TH^{\ell}_G}\,.
\end{multline}
Similarly, let $v_1, v_2 \in H^{\ell}_G$, then we have
\begin{equation*}
    \|\Phi(v_2)-\Phi(v_1)\|_{L^{\infty}_TH^{\ell}_G}
    \lesssim \||z^{\omega}+v_2|^2(z^{\omega}+v_2)-|z^{\omega}+v_1|^2(z^{\omega}+v_1)\|_{L^1_TH^{\ell}_G}\,,
\end{equation*}
so that
\begin{multline}
    \label{eq:contraction2terms}
    \|\Phi (v_2)-\Phi(v_1)\|_{L^{\infty}_TH^{\ell}_G} \lesssim \|(z^{\omega})^2  (\overline{v_2-v_1})\|_{L^{1}_TH^{\ell}_G} + \||z^{\omega}|^2 (v_2-v_1)\|_{L^{1}_TH^{\ell}_G}\\
    + \|z^{\omega}(|v_2|^2-|v_1|^2)\|_{L^{1}_TH^{\ell}_G}
     + \| z^{\omega} (\bar v_2^2-\bar v_1^2)\|_{L^{1}_TH^{\ell}_G}
       + T\||v_2|^2v_2-|v_1|^2v_1\|_{L^{\infty}_TH^{\ell}_G}\,.
\end{multline}
We now provide upper bounds for all the terms in~\eqref{eq:contraction1terms} and~\eqref{eq:contraction2terms}.
\begin{itemize}
\item We begin with $(z^{\omega})^2 \bar v$, $(z^{\omega})^2  (\overline{v_2-v_1})$, $|z^{\omega}|^2 v$ and $|z^{\omega}|^2 (v_2-v_1)$. Let $w=v$ or $w=v_2-v_1$. By the product law in $H^k_G$ of Proposition~\ref{lem.productLaw}, we have 
\begin{align*}
    \|(z^{\omega})^2\bar{w}(t)\|_{H_G^{\ell}} 
    & \lesssim \left(\|(z^{\omega})^2(t)\|_{H_G^{\ell}} + \|z^{\omega}(t)\|^2_{L^{\infty}_G}\right) \|w(t)\|_{H^{\ell}_G}\,.
\end{align*}
Using~\eqref{eq.1} and~\eqref{eq.5} from the assumption $\omega \in E_{R,T}$, this gives:
\begin{equation*}
    \|(z^{\omega})^2\bar{w}\|_{L^1_TH_G^{\ell}}  
    \lesssim \left(\|(z^{\omega})^2\|_{L^1_TH_G^{\ell}} + \|z^{\omega}\|^2_{L^2_TL^{\infty}_G}\right)\|w\|_{L^{\infty}_TH_G^{\ell}} \\
    \lesssim TR^2\|u_0\|_{\mathcal{X}_1^k}^2 \|w\|_{L^{\infty}_TH^{\ell}_G}\,.
\end{equation*}
Similarly one has 
\begin{equation*}
    \||z^{\omega}|^2 w\|_{L^1_TH_G^{\ell}}\lesssim TR^2\|u_0\|_{\mathcal{X}_{1}^k}^2 \|w\|_{L^{\infty}_TH^{\ell}_G}\,.
\end{equation*}

We have both proven
\[
     \|(z^{\omega})^2\bar{v}\|_{L^1_TH_G^{\ell}} +\||z^{\omega}|^2 v\|_{L^1_TH_G^{\ell}}\lesssim TR^2\|u_0\|_{\mathcal{X}_{1}^k}^2 \|v\|_{L^{\infty}_TH^{\ell}_G}
\]
and
\[
     \|(z^{\omega})^2(\overline{v_2- v_1})\|_{L^1_TH_G^{\ell}} +\||z^{\omega}|^2 (v_2-v_1)\|_{L^1_TH_G^{\ell}}\lesssim TR^2\|u_0\|_{\mathcal{X}_{1}^k}^2 \|v_2-v_1\|_{L^{\infty}_TH^{\ell}_G}\,.
\]
\item Let us estimate $z^{\omega}|v|^2$, $z^{\omega}(|v_2|^2-|v_1|^2)$, $z^{\omega}\bar v^2$ and $z^{\omega} (\bar v_2^2-\bar v_1^2)$. Using~\eqref{eq.4} from the assumption that $\omega \in E_{R,T}$, we infer:
\begin{equation*}
    \|z^{\omega}|v|^2\|_{L^1_TH^{\ell}_G} +\|z^{\omega}\bar v^2\|_{L^1_TH^{\ell}_G} 
     \lesssim TR\|u_0\|_{\mathcal{X}_{1}^k} \|v\|_{L^{\infty}_TH^{\ell}_G}^2
\end{equation*}
and
\begin{align*}
    \|z^{\omega}(|v_2|^2-|v_1|^2)\|_{L^1_TH^{\ell}_G} 
    &+\|z^{\omega}(\bar v_2^2-\bar v_1^2)\|_{L^1_TH^{\ell}_G} \\
     &\lesssim TR\|u_0\|_{\mathcal{X}_{1}^k} (\|v_2+v_1\|_{L^{\infty}_TH^{\ell}_G}^2+\|v_2-v_1\|_{L^{\infty}_TH^{\ell}_G}^2)\|v_2-v_1\|_{L^{\infty}_TH^{\ell}_G}^2\\
     &\lesssim TR\|u_0\|_{\mathcal{X}_{1}^k} (\|v_2\|_{L^{\infty}_TH^{\ell}_G}^2+\|v_1\|_{L^{\infty}_TH^{\ell}_G}^2)\|v_2-v_1\|_{L^{\infty}_TH^{\ell}_G}^2\,.
\end{align*}
\item Observe that thanks to the algebra property of $H^{\ell}_G$ (since $\ell >\frac{3}{2}$) of Lemma~\ref{lem.productLaw}, we have
\begin{equation*}
    \||v|^2v\|_{L^{\infty}_TH^{\ell}_G} \lesssim \|v\|_{L^{\infty}_TH^{\ell}_G}^3
\end{equation*}
and
\[
    \||v_2|^2v_2-|v_1|^2v_1\|_{L^{\infty}_TH^{\ell}_G} \lesssim (\|v_2\|_{L^{\infty}_TH^{\ell}_G}+\|v_1\|_{L^{\infty}_TH^{\ell}_G})^2\|v_2-v_1\|_{L^{\infty}_TH^{\ell}_G}\,.
\]
\end{itemize}
All these bounds combined together with assumption~\eqref{eq.3} in estimates~\eqref{eq:contraction1terms} and~\eqref{eq:contraction2terms} imply~\eqref{eq.contraction} and~\eqref{eq.contractionBis}.
\end{proof}


\begin{proof}[Proof of Theorem~\ref{th.main}] Let $\omega \in E_{R,T}$. Thanks Proposition~\ref{prop.contraction}, we know that there exists $C>0$ such that the map $\Phi : \mathcal{C}^0([0,T],H_G^{\ell}) \to \mathcal{C}^0([0,T],H_G^{\ell})$ is bounded Lipschitz on finite balls: if $\|v\|_{L^{\infty}_TH^{\ell}_G} \leqslant R\|u_0\|_{\mathcal{X}^k_1}$, we have
\[
\|\Phi (v)\|_{L^{\infty}_TH^{\ell}_G} \leqslant CT(R\|u_0\|_{\mathcal{X}^k_1})^3\,,
\]
and for $\|v_1\|_{L^{\infty}_TH^{\ell}_G} \leqslant R$ and $\|v_2\|_{L^{\infty}_TH^{\ell}_G} \leqslant R$, we have
\[
\|\Phi (v_2)-\Phi(v_1)\|_{L^{\infty}_TH^{\ell}_G} \leqslant CT(R\|u_0\|_{\mathcal{X}^k_1})^2\|v_2-v_1\|_{L^{\infty}_TH^{\ell}_G}\,.
\]
Thus, taking $T = \frac{1}{2C(R\|u_0\|_{\mathcal{X}^k_1})^2}$, we see that $\Phi$ stabilizes the ball $B(0,R\|u_0\|_{\mathcal{X}^k_1})$ in $\mathcal{C}^0([0,T],H_G^{\ell})$, moreover, $\Phi$ is a contraction on the ball $B(0,R\|u_0\|_{\mathcal{X}^k_1})$. The existence and uniqueness of $v$ solving~\eqref{eq:NLSGv} then follows from standard contraction mapping arguments. 

We have obtained that for any $\omega \in E_{R,T}$ (and $T=\frac{1}{2C(R\|u_0\|_{\mathcal{X}^k_1})^2}$), there exists a unique solution to~\eqref{eq:NLSG} in the space 
\[
    e^{it\Delta_G}u_0^{\omega}+\mathcal{C}^0([0,T],H^{\ell}_G) \subset \mathcal{C}^0([0,T],H^k_G)\,. 
\]
Then the set
\[
    E := \bigcup_{k \geqslant 1} \bigcap _{n \geqslant k}E_{n,\frac{1}{2n^2}}
\]
satisfies the requirements of Theorem~\ref{th.main}~(\textit{i}). Indeed, it remains to see that $\mathbb{P} (\Omega \setminus E)=0$. Since the sequence of sets $\bigcup_{n \geqslant k} E_{n,\frac{1}{2n^2}}$ is non-increasing, we have
\begin{equation*}
    \mathbb{P}(\Omega \setminus E)
     \leqslant \limsup_{k \to \infty} \mathbb{P}\left(\bigcup_{n \geqslant k} E_{n,\frac{1}{2n^2}}\right)
    \leqslant \limsup_{k \to \infty} \sum_{n \geqslant k} e^{-cn^2}\,,
\end{equation*}
which is $0$ since $\sum e^{-cn^2}$ converges. 
\end{proof}

\appendix 
\section{Appendices}

\subsection{Pointwise estimates on Hermite functions}\label{appendix:A}


The purpose of this appendix is to explain how to prove the estimates of Corollary~\ref{coro.pointwiseBoundsHermite} as a consequence of the pointwise estimates on the Hermite functions from Theorem~\ref{th.hermiteBounds}. 

\begin{proof}[Proof of Corollary~\ref{coro.pointwiseBoundsHermite}] We study the bounds on distinct regions of space. Let us fix $m \in\mathbb{N}$.  

\noindent(1) For $|x|\leqslant \frac{1}{2} \lambda_m$, there holds $|x^2-\lambda_m^2| \geqslant \frac{3}{4}\lambda_m^2$ thus
    \[|h_m(x)|\lesssim |x^2-\lambda_m^2|^{-\frac{1}{4}} \lesssim \lambda_m^{-\frac{1}{2}}\,.\]

\noindent(2) For $\frac{1}{2}\lambda_m\leqslant |x| \leqslant \lambda_m - \lambda_m^{-\frac{1}{3}}$, we have $x^2\leqslant \lambda_m^2-2\lambda_m^{\frac{2}{3}}+\lambda_m^{-\frac{2}{3}} \leqslant \lambda_m^2-\lambda_m^{\frac{2}{3}}$ thus $|\lambda_m^2-x^2| =(\lambda_m^2-x^2)\geqslant \lambda_m^{\frac{2}{3}}$, which implies that: 
    \[\lambda_m^{\frac{2}{3}}+ |\lambda_m^2-x^2| \leqslant 2 |\lambda_m^2-x^2|\,.\]
    Finally, we get
    \[|h_m(x)|\lesssim |\lambda_m^2-x^2|^{-\frac{1}{4}} \lesssim \left(\lambda_m^{\frac{2}{3}} + |\lambda_m^2-x^2|\right)^{-\frac{1}{4}}\,.\]

\noindent(3) For $| |x|-\lambda_m| \leqslant \lambda_m^{-\frac{1}{3}}$, we have $|x^2-\lambda_m^2| \leqslant \left||x|-\lambda_m\right|\cdot\left||x|+\lambda_m\right|\lesssim \lambda_m^{\frac{2}{3}}$, so that 
    \[|h_m(x)|\lesssim \lambda_m^{-\frac{1}{6}}=(\lambda_m^{\frac{2}{3}})^{-\frac{1}{4}} \lesssim \left(\lambda_m^{\frac{2}{3}} + |\lambda_m^2-x^2|\right)^{-\frac{1}{4}}\,.\] 

\noindent(4) For $\lambda_m+\lambda_m^{-\frac{1}{3}} \leqslant |x|\leqslant 2\lambda_m$ there holds $|x^2-\lambda_m^2|\gtrsim \lambda_m^{\frac{2}{3}}$, thus the crude bound $e^{-s_m(x)}\leqslant 1$ gives
    \[|h_m(x)|\leqslant \frac{e^{-s_m(x)}}{|x^2-\lambda_m^2|^{\frac{1}{4}}} \lesssim \left(\lambda_m^{\frac{2}{3}} + |x^2-\lambda_m^2|\right)^{-\frac{1}{4}}\,.\]

\noindent(5) Let $|x| \geqslant 2 \lambda_m$.  Observe that by change of variable $t=\lambda_m y$ we have
    \begin{equation*}
        s_m(x) 
        = \lambda_m^2\int_1^{\frac{x}{\lambda_m}} \sqrt{y^2-1}\,\mathrm{d}y
         \geqslant \lambda_m^2 \int_1^{\frac{x}{\lambda_m}} (y-1)\,\mathrm{d}y 
         = \frac{(x-\lambda_m)^2}{2}\,,
    \end{equation*}
    where we used that for $y\geqslant 1$, $\sqrt{y^2-1}\geqslant y-1$. Then, observe that $x-\lambda_m \geqslant \frac{x}{2}$ by definition of $x$. This implies $s_m(x) \geqslant \frac{x^2}{8}$ and finally, since $|x^2-\lambda_m^2|\geqslant \lambda_m^2 \geqslant 1$, we conclude:
    \[|h_m(x)| \lesssim \frac{e^{-s_m(x)}}{|x^2-\lambda_m^2|^{\frac{1}{4}}} \lesssim e^{-\frac{1}{8}x^2}\,.\qedhere\]
\end{proof}

\subsection{Algebra property, product laws and local Cauchy theory}\label{appendix:B}

\subsubsection{Proof of the functional inequalities}

In the Grushin case, the proof of Proposition~\ref{lem.productLaw} is a consequence of the following results. In the context of the Heisenberg sub-Laplacian, the proof of Proposition~\ref{lem.productLaw} about the algebra property of the Sobolev spaces $H^k$ can be found in~\cite{bahouriGallagher} and relies on representation theoretic formul\ae.


\begin{lemma}[See the proof of Lemma~3.6 in~\cite{BahouriFermanianGallagher2016}]
Let $H=\partial_{xx}+x^2$ denote the Harmonic oscillator. For all $k\geq 0$, there exists $C(k)>0$ such that for all $m\in\mathbb{N}$,
\[
\frac{1}{C(k) }\|H^{k/2}h_m\|_{L^2_x}
	\leq \|\partial_x^{k}h_m\|_{L^2_x}+\|x^{k}h_m\|_{L^2_x}
	\leq C(k)\|H^{k/2}h_m\|_{L^2_x}\,.
\]
\end{lemma}

\begin{corollary}\label{cor.equivNorms}
For all $k\geq 0$, there exists $C(k)>0$ such that for all $u\in H^k_G$, there holds
\begin{equation*}
\frac{1}{C(k)}\|(\operatorname{Id}-\Delta_G)^{k/2}u\|_{L^2_G}
	\leq \|\langle \partial _x\rangle ^ku\|_{L^2_G} + \|\langle x\partial_y\rangle ^k u\|_{L^2_G}
	\leq C(k)\|(\operatorname{Id}-\Delta_G)^{k/2}u\|_{L^2_G}\,.
\end{equation*}
\end{corollary}

\begin{proof}
We decompose $u$ as
\[
\mathcal{F}_{y\to\eta}u(x,\eta)=\sum_m f_m(\eta)h_m(|\eta|^{\frac{1}{2}}x)\,.
\]
Then
\begin{align*}
\|(\operatorname{Id}-\Delta_G)^{k/2}u\|_{L^2_x}^2
	&=\sum_m\int |f_m(\eta)|^2\mathrm{d}\eta \int (1+(2m+1)|\eta|)^k h_m(|\eta|^{\frac{1}{2}}x)^2\mathrm{d}x.
\end{align*}
Hence we see that $\|(\operatorname{Id}-\Delta_G)^{k/2}u\|_{L^2_G}\sim_k \|u\|_{L^2_x}+\|(-\Delta_G)^{k/2}u\|_{L^2_G}$,
and
\begin{align*}
\|(-\Delta_G)^{k/2}u\|_{L^2_G}^2
	&\sim_k	\sum_m\int |f_m(\eta)|^2\mathrm{d}\eta \int((2m+1)|\eta|)^k h_m(|\eta|^{\frac{1}{2}}x)^2\mathrm{d}x\\
	&=\sum_m\int |f_m(\eta)|^2\mathrm{d}\eta \int |\eta|^k(H^{k/2} h_m)(|\eta|^{\frac{1}{2}}x)^2\mathrm{d}x.
\end{align*}
Now we use a change of variables to get
\begin{align*}
\|(-\Delta_G)^{k/2}u\|_{L^2_G}^2
	&\sim_k\sum_m\int |f_m(\eta)|^2\mathrm{d}\eta  |\eta|^{k-1/2} \int (H^{k} h_m)(x)^2\mathrm{d}x.
\end{align*}
Then we use that $\|H^{k/2}h_m\|_{L^2_x}\sim_k \|\partial_x^{k}h_m\|_{L^2_x}+\|x^{k}h_m\|_{L^2_x}$ and get 
\begin{align*}
\|(-\Delta_G)^{k/2}u\|_{L^2_G}^2
	&\sim_k\sum_m\int |f_m(\eta)|^2\mathrm{d}\eta |\eta|^{k-1/2} \int ((\partial_x^{k}+x^{k}) h_m)(x)^2\mathrm{d}x\\
	&= \sum_m\int |f_m(\eta)|^2\mathrm{d}\eta  \int|\eta|^{k} ((\partial_x^{k}+x^{k}) h_m)(|\eta|^{\frac{1}{2}}x)^2\mathrm{d}x\\
	&=	\sum_m\int |f_m(\eta)|^2\mathrm{d}\eta \int (\partial_x^k +(|\eta|x)^k)h_m(|\eta|^{\frac{1}{2}}x)^2\mathrm{d}x\\
	&=\|\partial_x^ku\|_{L^2_G}^2+\|(x\partial_y)^ku\|_{L^2_G}^2.\qedhere
\end{align*}
\end{proof}


We are now ready to give a proof of the product laws. 

\begin{proof}[Proof of Proposition~\ref{lem.productLaw}] 
\noindent (\textit{i}) We start by using the above Corollary~\ref{cor.equivNorms} and get
\begin{align*}
    \|uv\|_{H^k_G}
    &\lesssim \|\langle\partial_x\rangle ^k(uv)\|_{L^2_{x,y}} + \|\langle x \partial_y \rangle ^k (uv)\|_{L^2_{x,y}}\,. 
\end{align*}
Now, the classical product rule in Sobolev spaces applied in $x$ (resp. in $y$) implies 
\[\|\langle\partial_x\rangle ^k(uv)\|_{L^2_x} \lesssim \|\langle \partial _x \rangle ^ku\|_{L^2_x}\|v\|_{L^{\infty}_x} + \|u\|_{L^{\infty}_x}\|\langle \partial_x\rangle ^k v\|_{L^2_x}\,,\]
resp. 
\[\|\langle x\partial_y\rangle ^k(uv)\|_{L^2_y} \lesssim \|\langle x\partial _y \rangle ^ku\|_{L^2_y}\|v\|_{L^{\infty}_y} + \|u\|_{L^{\infty}_y}\|\langle x\partial_y\rangle ^k v\|_{L^2_y}\,,\]
which combined with Hölder estimates in $y$ (resp. $x$) yields
\begin{align*}
  \|uv\|_{H^k_G}
  \lesssim &   \|\langle \partial _x \rangle ^ku\|_{L^2_{x,y}}\|v\|_{L^{\infty}_{x,y}} + \|u\|_{L^{\infty}_{x,y}}\|\langle \partial_x\rangle ^k v\|_{L^2_{x,y}} \\
  & + \|\langle x\partial _y \rangle ^ku\|_{L^2_{x,y}}\|v\|_{L^{\infty}_{x,y}} + \|u\|_{L^{\infty}_{x,y}}\|\langle x\partial_y\rangle ^k v\|_{L^2_{x,y}}\,.
\end{align*}
Proposition~\ref{lem.productLaw}~(\textit{i}) follows by an other application of Corollary~\ref{cor.equivNorms}. \medskip

\noindent (\textit{ii}) is a direct consequence of (\textit{i}) and the Sobolev embedding $H^k_G\hookrightarrow L^{\infty}_G$ when $k>\frac{3}{2}$. \medskip 

\noindent (\textit{iii}) When $p$ is an integer, the result follows from ~(\textit{ii}) by iteration.
\end{proof}

\begin{lemma}[Limit Sobolev embedding] The following statements hold.  
\begin{enumerate}[label=(\roman*)]
    \item There exists $C>0$ such that for every $p>2$ and $u\in H^{\frac 32}_G$, there holds
\begin{equation}
    \label{eq:trudinger1}
    \|u\|_{L^p_G} \leq C \sqrt{p}\|u\|_{H_G^{\frac{3}{2}}}\,.
\end{equation}
\item (Brezis-Gallouët) For any $k>\frac{3}{2}$, there exists $C_k>0$ such that there holds 
\begin{equation}
    \label{eq:trudinger2}
    \|u\|_{L^{\infty}_G}\leq C_k \|u\|_{H^{\frac{3}{2}}_G} \log^{\frac{1}{2}}\left(1 + \frac{\|u\|_{H^{k}_G}}{\|u\|_{H^{\frac{3}{2}}_G}}\right)\,.  
\end{equation}
\end{enumerate}
\end{lemma}

\begin{proof} (\textit{i}) Let $p>2$ and $u \in L^p_G$. By the triangle inequality we have
\[
    \|u\|_{L^p_G} \leqslant \sum_{A \in 2^{\mathbb{N}}} \|u_A\|_{L^p_G}\,.
\]
Then, the Sobolev embedding yields
\[
    \|u\|_{L^p_G} \lesssim \sum_{A \in 2^{\mathbb{N}}} \|u_A\|_{H_G^{3\left(\frac{1}{2}-\frac{1}{p}\right)}}
    \lesssim \sum_{A \in 2^{\mathbb{N}}} A^{-\frac{3}{2p}}\|u_A\|_{H_G^{\frac{3}{2}}}\,.
\]
An application of the Cauchy-Schwarz inequality and orthogonality provide us with
\[
    \|u\|_{L^p_G} \lesssim \left(\sum_{A \in 2^{\mathbb{N}}} A^{-\frac{3}{p}}\right)^{\frac{1}{2}}\|u\|_{H^{\frac{3}{2}}_G}\,,
\]
which gives the conclusion since $\sum_{A \in 2^{\mathbb{N}}} A^{-\frac{3}{p}} \sim Cp$ as $p$ goes to infinity.

(\textit{ii}) Let $A_0\geqslant 1$ be a dyadic integer. Start by writing
\(
    u= \sum_{A\leqslant A_0}u_A + \sum_{A>A_0}u_A\,,
\)
and by Cauchy-Schwarz
\[
    \|u\|_{L^{\infty}} \leqslant \log^{\frac{1}{2}}(A_0) \left(\sum_{A \leqslant A_0}\|u_A\|^2_{L^{\infty}}\right)^{\frac{1}{2}} + \sum_{A >A_0} \|u_A\|_{L^{\infty}}^2\,.
\]

Observe that for $p\geqslant 2$, $\|u_A\|_{L^p_G} \lesssim A^{\frac{3}{2}\left(\frac{1}{2}-\frac{1}{p}\right)}\|u_A\|_{L^2_G}$, which gives when letting $p \to \infty$
\begin{equation}
    \label{eq:trudBound1}
    \|u_A\|_{L^{\infty}_G} \lesssim \|u_A\|_{H^{\frac{3}{2}}_G}\,,
\end{equation}
and this latter inequality also implies
\begin{equation}
    \label{eq:trudBound2}
    \|u_A\|_{L^{\infty}_G} \lesssim A^{-\frac{1}{2}(k-\frac{3}{2})}\|u_A\|_{H^k_G}\,.
\end{equation}

The bound~\eqref{eq:trudBound1} when $A\leqslant A_0$ and~\eqref{eq:trudBound2} when $A>A_0$ imply
\[
    \|u\|_{L^{\infty}_G} \lesssim \|u\|_{H^{\frac{3}{2}}_G}\log^{\frac{1}{2}} A_0 + A_0^{-\frac{1}{2}(k-\frac{3}{2})}\|u\|_{H^k_G}\,,
\]
which gives the result after optimisation in $A_0$.
\end{proof}

\subsubsection{Deterministic local Cauchy theory}

We finish this appendix with a summary of well-posedness results for
~\eqref{eq:NLSH} (resp.~\eqref{eq:NLSG}) for $k>2$ (resp. $k>\frac{3}{2}$).
To the best of our knowledge, the best well-posedness result for~\eqref{eq:NLSH} is the following. 

\begin{proposition}[Well-posedness for~\eqref{eq:NLSH}, see~\cite{bahouriGallagher}] For $k>2$, the Cauchy problem for~\eqref{eq:NLSH} is locally well-posed in $\mathcal{C}^0([0,T^*),H^k(\mathbb{H}^1))$ and $T^*=T^*(\|u_0\|_{H^k(\mathbb{H}^1)}) \gtrsim \|u_0\|_{H^k(\mathbb{H}^1)}^{-2}$.
\end{proposition}

Similarly, we recall the best local theory for~\eqref{eq:NLSG}. 

\begin{proposition}[Well-posedness theory for~\eqref{eq:NLSG}]\label{prop:localTheory} The following well-posedness statements hold:
\begin{enumerate}[label=(\roman*)] 
\item For $k>\frac{3}{2}$, the Cauchy problem for~\eqref{eq:NLSG} is locally well-posed in $\mathcal{C}([0,T^*),H^k_G)$. Moreover, for $u_0\in H^k_G$, the maximal time $T^*$ satisfies $T^* \gtrsim \|u_0\|_{L^{\infty}_G}^{-2}$. 
\item The blow-up criterion can be refined:
\[T^{*}<\infty \Longrightarrow \|u(t)\|_{H^{\frac{3}{2}}_G} \underset{t \to T^*}{\longrightarrow}\infty\,.\]
\end{enumerate}
\end{proposition}

\begin{proof}[Sketch of proof for Proposition~\ref{prop:localTheory}] 
We only give formal arguments, which are easily converted into fully rigorous proofs by standard means. 

(\textit{i}) Let $k>\frac{3}{2}$. Applying the operator $(-\Delta_G)^{\frac{k}{2}}$ to~\eqref{eq:NLSG}, then multiplying by $(-\Delta_G)^{\frac{k}{2}}u$ and integrating by parts in space, we compute that:
\[
    \frac{\mathrm{d}}{\mathrm{d}t}\|u(t)\|_{H^k_G}^2 \lesssim \|(-\Delta_G)^{\frac{k}{2}}(|u|^2u)(-\Delta_G)^{\frac{k}{2}}\bar u\|_{L^1_G}\,.
\]
Applying the Hölder inequality, the algebra property of Lemma~\ref{lem.productLaw} and the Hölder inequality again, we finally arrive at the estimate 
\[
    \frac{\mathrm{d}}{\mathrm{d}t}\|u(t)\|_{H^k_G}^2 \lesssim \|u\|_{L^{\infty}_G}^2\|u\|_{H^k_G}^2\,,
\]
which by the Sobolev embedding and the Grönwall inequality gives an \textit{a priori} estimate in $H^k$ and implies the local theory.

(\textit{ii}) This follows from the energy estimate and inequality~\eqref{eq:trudinger2}, which give
\begin{equation*}
    \frac{\mathrm{d}}{\mathrm{d}t}\|u(t)\|_{H^{k}_G}^2 
     \lesssim \|u\|_{L^{\infty}_G}^2\|u\|_{H^{k}_G}^2
     \lesssim \|u\|_{H^{\frac{3}{2}}}^2\|u\|_{H_G^k}^2\log\left(1 + \frac{\|u\|_{H^{k}_G}}{\|u\|_{H^{\frac{3}{2}}}}\right)\,,
\end{equation*}
which implies the result after an application of Grönwall's inequality. 
\end{proof}

Note that a similar argument to that of (\textit{ii}), which relies on the inequality~\eqref{eq:trudinger1} can be used to prove that, if solutions exists in $H^{\frac{3}{2}}_G$, they are unique. This argument goes back to Yudovich~\cite{Yudovich1963}, then has been used by Vladimirov~\cite{Vladimirov1984} in the context of Schrödinger equations. Namely, if $u_1$, $u_2$ are two $H^{\frac{3}{2}}_G$ solutions in $L^{\infty}([0,T),H^{\frac{3}{2}}_G$), introduce $\phi(t)=\|u_1(t)-u_2(t)\|^2_{L^2_G}$ and fix $T_1 <T$, we prove that $\phi (t)=0$ for all $t\in [0,T_1]$. Denote by $w(t)=u_1(t)-u_2(t)$. Then 
\[
    \phi'(t)=2\int_{\mathbb{R}^2} w'(t)\bar{w}(t) \,\mathrm{d}x=2i\int_{\mathbb{R}^2}\Delta w(t) \bar{w}(t) \,\mathrm{d}x - 2i\int_{\mathbb{R}^2}(|u_1|^2u_1-|u_2|^2u_2)\bar{w}(t)\,\mathrm{d}x  
\]
Since the first terms equals $-2i \|\nabla w\|^2_{L^2} \in i\mathbb{R}$ and $\phi$ is real-valued, we have 
\[
    \phi'(t) \leqslant 2\left\vert \int_{\mathbb{R}^2}(|u_1|^2u_1-|u_2|^2u_2)\bar{w}(t)\,\mathrm{d}x \right\vert \lesssim \int_{\mathbb{R}^2} |w(t)|^2\left(|u_1(t)|^2+|u_2(t)|^2\right)\,\mathrm{d}x\,.
\]
Then for all $t \leqslant T_1$, we have
\begin{align*}
    \phi '(t) & \lesssim \int_{\mathbb{R}^2} |(u_1-u_2)(t)|^{2(1-1/p)}(|u_1(t)|^{2(1+1/p)}+|u_2(t)|^{2(1+1/p)})\,\mathrm{d}x \\
    & \lesssim \sqrt{p}\phi(t)^{1-\frac{1}{p}}\left(\|u_1\|_{H^{\frac{3}{2}}_G}^{2(1+1/p)} + \|u_2\|^{2(1+1/p)}_{H^{\frac{3}{2}}_G}\right)\,, \end{align*}
where we use~\eqref{eq:trudinger1} in the last step. Since $\|u_1(t)\|_{L^{\infty}([0,T_1],H^{\frac{3}{2}})} = C(T_1)$ and $2(1+\frac{1}{p}) \leqslant 3$, we obtain 
\[
    \phi '(t) \lesssim \sqrt{p}\phi(t)^{1-\frac{1}{p}}
\]
The, we integrate on $[0,t]$ to get
\[
    \phi (t) \leqslant \left(\frac{t}{\sqrt{p}}\right)^p\,,
\]
which goes to $0$ as $p \to \infty$, hence $\phi =0$.  

\subsection{An interpolation lemma}\label{appendix:C}
The aim of this appendix is to prove that if inequality~\eqref{eq:term_stein}
    \begin{equation*}
        \left\|\sum_{2B>A}(\operatorname{Id}-\Delta_G)^{{\ell}/2}(u_{A}(P_{>A}v_{B}))\right\|_{L^2_G}^2
	    \lesssim \sum_{\delta\in D_1}\|(u_{A})^{\delta_1} \|_{L^{\infty}_G}^2\|v\|_{H^{\ell}_G}^2
    \end{equation*}
holds for $\ell=0$ and $\ell=2$, then this inequality holds for all $\ell\in[0,2]$ by interpolation.
Writing $w=(\operatorname{Id}-\Delta_G)^{{\ell}/2}v$, we have
    \[
        \left\|\sum_{B>A}(\operatorname{Id}-\Delta_G)^{{\ell}/2} (u_{A}(P_{>A}v_{B}))\right\|_{L^2_G}
	    =\left\| (\operatorname{Id}-\Delta_G)^{{\ell}/2} \left(u_A \sum_{2B> A} (\operatorname{Id}-\Delta_G )^{-{\ell}/2}(P_{>A}w)_B\right)\right\|_{L^2_G}\,.\qedhere
    \] 
 
\begin{lemma}[Interpolation lemma]\label{lem.interpol} 
Let ${\ell}\in[0,1]$, then for any $w \in L^2_G$, there holds
\begin{equation*}
    \| (\operatorname{Id}-\Delta_G)^{{\ell}}\left(u_A\chi_{> A} (\operatorname{Id}-\Delta_G)^{-{\ell}} w\right)\|_{L^2_G}
    \lesssim \left(\sum_{\delta\in D_1}\|(u_{A})^{\delta_1}\|_{L^{\infty}_G}^2\right)^{1/2}\|w\|_{L^2_G}\,.
\end{equation*}
\end{lemma}

The proof of this lemma is an application of Stein's interpolation theorem in the case when the operators are bounded by the same constant.

\begin{theorem}[Stein interpolation theorem~\cite{Stein1956interpolation}, Theorem 1]\label{th.stein} Let $(\Omega _i, \Sigma_i, \mu_i)$, $i=0, 1$ be two measured spaces, $p_i, q_i \in [1,\infty]$, $S:= \{z \in \mathbb{C}\mid 0 <Re(z)<1\}$ and $(T_z)_{z \in \bar S}$ a family of operators from simple functions in $L^1(\mu_1)$ to $\mu_2$-measurable functions. Assume that there exists $c<\pi$ such that the following holds.
\begin{enumerate}[label=(\textit{\roman*})]
    \item For any fixed simple functions $f$ and $g$ on $\Omega_0$ and $\Omega_1$ respectively, $z \in\bar S\mapsto \int_{\Omega_1} T_z(f)g\d\mu_1$ is continuous on $\bar S$ and holomorphic in $S$, and satisfies
    \begin{equation*}
        \sup_{z \in S} e^{-c|\mathrm{Im}(z)|}\log \left\vert\int_{\Omega_1} T_z (f) g\d\mu_1\right\vert < \infty\,.
    \end{equation*}
    \item The operator $T_z : L^{p_0}(\Omega_0) \to L^{q_0}(\Omega_1)$ is continuous whenever $Re(z)=0$: there exists $C_0$ such that for all $f\in L^{p_0}(\Omega_0)$ and $r\in\mathbb{R}$,
    \[
    \|T_{0+ir}f\|_{L^{q_0}}\leq C_0(r)\|f\|_{L^{p_0}},
    \]
    and similarly, whenever $Re(z)=1$, the operator $T_z : L^{p_1}(\Omega_0) \to L^{q_1}(\Omega_1)$ is continuous: there exists $C_1$ such that for all $f\in L^{p_1}(\Omega_0)$ and $r\in\mathbb{R}$,
    \[
    \|T_{1+ir}f\|_{L^{q_1}}\leq C_1(r)\|f\|_{L^{p_1}}.
    \]
    Moreover, for $i\in\{0,1\}$,
    \begin{equation*}
        \sup_{r\in\mathbb{R}} e^{-c|r|} \log |C_i(r)|< \infty\,.
    \end{equation*}
\end{enumerate}
Then for $\theta \in [0,1]$, setting
\(\frac{1}{p_{\theta}} = \frac{1-\theta}{p_0}+\frac{\theta}{p_1} \text{ and } \frac{1}{p_{\theta}} = \frac{1-\theta}{q_0}+\frac{\theta}{q_1},\)
the operator $T_{\theta} : L^{p_{\theta}}(\Omega_0)  \to L^{q_{\theta}}(\Omega_1) $ is bounded. More precisely, there exists $C_\theta=C(\theta,C_0,C_1)$ such that
\[
\|T_{\theta}f\|_{L^{q_\theta}}
	\leq C_\theta\|f\|_{L^{p_\theta}},
    \quad f\in L^{p_\theta}(\Omega_0), \quad r\in\mathbb{R}.
\]
\end{theorem}

\begin{proof}[Proof of Lemma~\ref{lem.interpol}] We apply Theorem~\ref{th.stein} to the operators $\left(\sum_{\delta\in D_1}\|(u_{A})^{\delta_1}\|_{L^{\infty}_G}^2\right)^{-1/2} T_z$, with
\[T_zf
	:= (\operatorname{Id}-\Delta_G)^{z}\left(u_AP_{> A}(\operatorname{Id}-\Delta_G)^{-z}f\right)\,,\]
and $(\operatorname{Id}-\Delta_G) ^{z}=\exp \left(z\log (\operatorname{Id}-\Delta_G)\right)$.
We make the choice $\Omega_0=\Omega_1=\mathbb{R}^2$, $\mu_0=\mu_1$ is the Lebesgue measure $\lambda$, and $p_0=p_1=p_{\theta}=q_{\theta}=q_0=q_1=2$. 

Observe that for any real $r$, the operator $(\operatorname{Id}-\Delta_G)^{ir} $ acts by rotation of the Fourier coefficients in $L^2_G$ thus it is bounded in $L^2$ with norm $1$. Indeed, the decomposition $\mathcal{F}(w)(x,\eta)=\sum_{m\in\mathbb{N}}f_m(\eta)h_m(|\eta|^{\frac{1}{2}}x)$ leads to
\[
\mathcal{F}((\operatorname{Id}-\Delta_G)^{ir}w)(x,\eta)=\sum_{m\in\mathbb{N}}(1+(2m+1)|\eta|)^{ir}f_m(\eta)h_m(|\eta|^{\frac{1}{2}}x)\,,
\]
so that we can conclude by orthogonality of the decomposition. This implies that for all $w\in L^2(\mathbb{R}^2)$ and $z\in \bar S$, we have
\begin{equation*}
\|T_zf\|_{L^2(\mathbb{R}^2)}
	=\|T_{\mathrm{Re}(z)}((\operatorname{Id}-\Delta_G)^{-ir}f)\|_{L^2(\mathbb{R}^2)}\,,
\end{equation*}
with $\|(\operatorname{Id}-\Delta_G)^{-ir}w\|_{L^2(\mathbb{R}^2)}=\|w\|_{L^2(\mathbb{R}^2)}$. Thanks to the fact that inequality~\eqref{eq:term_stein} holds in the case ${\ell}=0$ and ${\ell}=2$, assumption (\textit{ii}) from Theorem~\ref{th.stein} then holds with some constant functions $C_0=C_1$ independent of $r=\mathrm{Im}(z)$.

We now establish assumption (\textit{i}). Fix two simple functions $f,g$ on $\mathbb{R}^2$. Then the map  $z \in\bar S\mapsto \int_{\mathbb{R}^2} T_z(f)g\d\lambda$ is continuous and holomorphic. Moreover, for all $z\in S$, we have
\begin{align*}
\left\vert \int_{\mathbb{R}^2} T_z(f)g\d\lambda\right\vert
	&\leqslant \|T_zf\|_{L^2(\mathbb{R}^2)}\|g\|_{L^2(\mathbb{R}^2)}\,.
\end{align*}
We write $z = {\ell} + ir$ with $({\ell},r)\in[0,1]\times\mathbb{R}$.
When ${\ell}=0$, we simply write
\[
\|T_zf\|_{L^2_G}
	\lesssim \|u_A\|_{L^{\infty}_G}\|P_{> A}f\|_{L^2_G}\,,
\]
and observe that $\|u_A\|_{L^{\infty}_G}^2\lesssim \|u_A\|_{H^{\frac{3}{2}+\varepsilon}_G}^2\lesssim A^{\frac{3}{2}+\varepsilon}\|u_A\|_{L^2_G}^2<+\infty$. Otherwise, if ${\ell}>0$, start by an application of the product law in Proposition~\ref{lem.productLaw}, which gives
\begin{align*}
    \|T_zf\|_{L^2_G}&=\|u_A P_{>A}(\operatorname{Id}-\Delta_G)^{-z}f\|_{H^{{\ell}}_{G}} \\
    & \lesssim \|u_A\|_{L^{\infty}_G} \|P_{>A}(\operatorname{Id}-\Delta_G)^{-z}f\|_{H^{\ell}_G} + \|u_A\|_{H^h_G} \|P_{>A}(\operatorname{Id}-\Delta_G)^{-z}f\|_{L^{\infty}_G} \\
    & \lesssim \|u_A\|_{L^{\infty}_G}\|f\|_{L^2_G} + \|u_A\|_{H^{\ell}_G}\|P_{>A}(\operatorname{Id}-\Delta_G)^{-ir}f\|_{W^{-{\ell},\infty}_G}\,.
\end{align*}
Then we observe that $\|u_A\|_{H^{\ell}_G} \lesssim A^{\frac{{\ell}}{2}}\|u_A\|_{L^2_G}< \infty$, $\|u_A\|_{L^{\infty}_G}<\infty$, as well as  $\|f\|_{L^2_G}< \infty$. 
It remains to study $\|P_{>A}(\operatorname{Id}-\Delta_G)^{-ir}f\|_{W^{-{\ell},\infty}_G}$. We first use the dual Sobolev embedding $L^p_G\hookrightarrow W^{-{\ell},\infty}_G$ where $\frac{1}{p}-\frac{{\ell}}{3}=0$ (so that $3\leq p<\infty$):
\[
    \|P_{>A}(\operatorname{Id}-\Delta_G)^{-ir}f\|_{W^{-{\ell},\infty}_G} 
    \lesssim \|P_{>A}(\operatorname{Id}-\Delta_G)^{-ir}f\|_{L^p_G}\,.
\]
Now, we conclude by using the continuity of $P_{> A}(\operatorname{Id}-\Delta_G)^{-ir}$ on $L^p_G$.

Indeed, $P_{> A}(\operatorname{Id}-\Delta_G)^{-ir}=F(-\Delta_G)$, where
\[
F(\lambda)=\left(1-\chi\left(\frac{1+\lambda}{A}\right)\right)(1+\lambda)^{ir}.
\]
Since $\chi\in\mathcal{C}_c^{\infty}[0,1)$, one knows that $F\in W^{2,\infty}(\mathbb{R})$. Now, Theorem~1 in~\cite{MartiniSikora2012} implies that for all $p\in(1,\infty)$, $F(-\Delta_G)$ is bounded in $\mathcal{L}(L^p(\mathbb{R}^2))$.
\end{proof}

\bibliographystyle{alpha}
\bibliography{biblio}

\end{document}